\documentclass[12pt,a4paper]{report}
\usepackage{t1enc} 
\usepackage[latin1]{inputenc}
\usepackage[english]{babel}
\usepackage{amsfonts} 
\usepackage{amsthm}
\def\M{\mathcal{M}}
\def\e{\varepsilon}
\def\R{\mathbb{R}}

\def\Cent{\hookrightarrow^c}
\def\A{A_{\varepsilon}}
\def\G{\Gamma_{\varepsilon}}
\def\Ae{{\cal{A}}_{\varepsilon}}

\def\cB{{\cal{B}}}
\def\F1{{\cal F}}
\def\pf{\begin{proof}[{\bf Proof:}]\thst}
\newcommand{\Hm}[1]{\mathcal{H}^{#1}}
\newcommand{\cl}[1]{\mathcal{#1}}
\newcommand{\p}[1]{\theta_{#1,\cdot}}
\newcommand{\T}[1]{\theta_{{#1}}}

\def\thst{\mbox{}\newline}
\newtheorem{def1}{Definition}[chapter]
\newtheorem{thm}{Theorem}[chapter]
\newtheorem{lem}{Lemma}[chapter]
\newtheorem{cor}{Corollary}[chapter]
\newtheorem{conj}{Conjecture}[chapter]
\newtheorem{prop}{Proposition}[chapter]
\newtheorem{cons}{Construction}[chapter]
\newtheorem{ques}{Question}[chapter]

\begin{document}
\title{Approximately $j$-dimensional Koch type sets are potentially minimal surface singularities.}
\author{Amos Koeller} 
\date{}
\maketitle
\tableofcontents

\chapter{Introduction}

In the fields of Geometric Measure Theory and Differential Geometry we find that the study of surfaces (Minimal Surfaces, 
Stationary Surfaces, Energy Minimising Surfaces, etc.) and the flows of surfaces (Mean Curvature flow, Ricci flow, 
Brakke flow, etc.) play a central role. Such objects are not in general well behaved in that they have initially, or 
develop in finite time, singularities. Simply speaking, these can be thought of as holes, edges, corners, or in general 
points that around which no neighbourhood can be described by a graph. It is natural then that an understanding of the 
structure of such sets would be desired.
\newline \newline
In the present state of knowledge surprisingly little is known about these sets. Although, particularly in weak formations, 
regularity theorems are relatively standard in studies of these objects (See in particular Allard, White, Simon, Brakke, 
Ecker), this tells us more about how much of the surface we may consider as being smooth (or regular) than about the 
structure or measure of the singularity set itself.
\newline \newline
Some important results on the structure of the singularity sets themselves are due to White, whose stratification results 
show that the dimension of the singularity set is at least 1 less than that of the surface, and Simon, who has shown that 
in a particular class of minimal surfaces the singularity set is always a finite union of countably rectifiable sets in 
the dimension of the singularity.
\newline \newline
What is not known is anything at all about the shape of a singularity sets. We do indeed have examples of singularity sets 
but they are all simple, (i.e. the subset of a line, or a point) which leaves a lot of space between examples and 
generally provable results.
\newline \newline
In his paper showing the rectifiability of singularity sets of a certain class of minimal surfaces, Simon shows that 
singularity sets can be approximated by planes in the dimension in which they occur. In mean curvature flow, Huisken 
and Sinestrari have shown that blow ups around singularity points lead to eventually bounding the singularity set 
(blown up) in a cylinder. This, when considering the axis of the cylinder as a plane in the appropriate dimension 
is again an approximation to a plane in the dimension of the singularity.
\newline \newline
This tells us that the properties of sets that are approximately planes of some dimension are worth considering 
to see what properties we can get "for free" and what sort of potential problems does one need to be wary of when 
considering the singularity sets.
\newline \newline
As a model for what is meant when we say that a set is approximately a $j$-dimensional plane or indeed that a set is 
approximately $j$-dimensional we use the 'plane like' properties shown by Simon to be possessed by singularity set 
approximations.
\newline \newline
We isolate these properties to construct an ordering of eight strengths of $j$-dimensional plane approximation of which the 
property combination specifically used by Simon is the fourth. We classify these definitions in terms of whether or not 
they ensure actual $j$-dimensionality and whether or not they ensure locally $\Hm{j}$-finite measure in either a strong or 
a weak sense.
\newline \newline
The definitions allow for the full spectrum of possibilities. The strongest definition implying that the set is locally 
a finite union of Lipschitz graphs and the weakest two do not even ensure that the set be $j$-dimensional.
\newline \newline
The most interesting case, however, is that of the complications of our fourth definition, intriguingly the same as 
that arising in Simons work. This definition ensures $j$-dimensionality, but what makes this case interesting is that while 
locally finite $j$-dimensional measure is not ensured, any counter examples are necessarily exotic. We show that while 
satisfying 'approximately $j$-dimensional' properties such sets have points of infinite $\Hm{j}$-density 
but that no piece of any Lipschitz graph may pass through such a point. This rules out any  vaguely well behaved sets 
(or countable unions of vaguely well behaved sets) from both satisfying our fourth definition and failing to have 
locally finite $\Hm{j}$-measure.
\newline \newline
Since our classification is complete it follows that we can (and indeed do) provide a set satisfying this fourth 
definition that also does not have locally finite $j$-dimensional measure. The set is a variation on the fractal known 
as the Koch set. Since all singularity sets are closed we go on to show that a closed version of this counter example 
exists which implies that in principle singularity sets could be as terribly behaved as the counter example.
\newline \newline
Especially since, at least in the minimal surface case, singularity sets are known to be finite unions of countably 
$j$-rectifiable sets (see $\cite{simon2}$) the question of whether such sets as these counter examples are 
finite unions of countably 
$j$-rectifiable sets (and so continue to, potentially, be singularity sets) becomes of interest.
\newline \newline
The answer to this question for the particular examples initially given turns out to be no, they are not rectifiable 
without considering measure conditions and so cannot be finite unions of countably $j$-rectifiable sets. However, 
since the explicitly constructed counter examples are members of a family of constructions this by no means rules out 
the possibility of very poorly behaved singularity sets.
\newline \newline
The second part of the work then defines generalisations of the construction of the constructed counter examples. We call 
these sets, due to their similarity to the Koch sets, Koch-type sets. We then concentrate on giving dimension, measure 
and rectifiability conditions for these generalised sets.
\newline \newline
We find, encouragingly for the study of singularity sets that should such a set be first of all rectifiable then it can 
also be written as a single Lipschitz graph.
\newline \newline
This would immediately imply, since we need to remove the 'corners' of the sets in order to satisfy our fourth definition 
that any singularity set that may be of a Koch type set form should also be a subset of a single Lipschitz graph.
\newline \newline
The structure of the work is as follows:
\newline \newline
In chapter 2 we present a more precise formulation of the motivating mathematics
including some particularly relevant standard general geometric measure theoretic definitions and results, 
provide the list of definitions as well as the results already known in terms of our classification aims and results from 
which looked for classification results are a short corollary. 
\newline \newline
In chapter 3 we construct the specific counter examples that will be used in our classifications including the explicit 
examples of Koch type sets mentioned above. We go on to prove some important properties of these sets. Some properties, 
for example dimension, follows from some relatively general previous results of Hutchinson (see $\cite{hutch}$). Since 
it is often instructive to see the direct proof for explicit examples we provide direct proofs for these results as well. 
\newline \newline
Before moving on to show that the counter examples do indeed satisfy the definitions that they are counter examples to, a 
by no means trivial task, we show in chapter 4 that the complexity of the counter examples constructed is indeed necessary; 
in that no 'simple' example could possibly suffice. Further, we show the path to showing that singularity sets have locally 
finite measure is shorter than was previously thought, in that we need only show that the set is graph possesing at 
all points of infinite density. This is shorter than previously thought since such a property is so very weak. It does not 
even require that the set be weakly locally countably rectifiable.
\newline \newline
In chapter 5 we fit the counter examples to their respective definitions and complete the task of classifying the 
definitions.
\newline \newline
Chapter 6 gathers a few other miscellaneous relevant results and describes dimension generalisation of the explicit 
counter examples which are constructed to satisfy approximations to dimension 1 (though, of course, some are 
actually of fractal dimension between 1 and 2.) 
\newline \newline
Finally, in Chapters 7 and 8, we deal with the question of dimension, measure and rectifiability of the family of sets 
that are the generalised form of the explicit counter examples given. These generalisations are divided into two levels 
of generalisation, first and second degree variation. We keep the two levels of generalisation distinct since, 
although first degree variation generalisations are also second degree variation generalisations, they allow for 
stronger results. This is because much more can 'go wrong' in the second degree variation case.

\chapter{Background, Definition and Existing Results}
\section{Preliminary Geometric Measure Theory}
We start straight of with some relevant measure theoretic background. 
The standard references are of course $\cite{simon1}$ and $\cite{federer}$. We assume basic familiarity with 
general measure theory and we use the usual symbol for $r$-dimensional Hausdorff measure $\Hm{r}$ for 
$r \in \R$. Also, we denote the Hausdorff volume of the unit $n$-ball by $\omega_n$. \newline \newline
As mentioned, a major part of our investigation regards dimension, for which we are interested in Hausdorff 
dimension which we define as follows.
\begin{def1}\label{defdim} \thst
Set $A \subset \R^n$ for some $n \in \mathbb{N}$. Then the {\bf Hausdorff dimension} of $A$ is defined as 
\begin{eqnarray}
dimA & := & \inf\{r\in\R :\Hm{r}(A)=0\} \nonumber \\
& = & \sup\{r\in\R :\Hm{r}(A)=+\infty\} \nonumber
\end{eqnarray}
\end{def1} \noindent
Another important quantity that we will be using  is density, and indeed $n$-dimensional density.
\begin{def1}\label{defden}\thst
Let $(X,\cl{B},\mu)$ be a measure space.Then for any subset $A$ of $X$, and any point $x \in X$, we define the 
$n$-dimensional upper and lower $n$-dimensional densities $\Theta^{*n}(\mu,A,x)$, $\Theta_*^n(\mu,A,x)$ respectively by 
$$\Theta^{*n}(\mu,A,x) = \limsup_{\rho \rightarrow 0}{{\mu(A \cap B_{\rho}(x))}\over{\omega_n\rho^n}}$$ and 
$$\Theta^{n}_*(\mu,A,x) = \liminf_{\rho \rightarrow 0}{{\mu(A \cap B_{\rho}(x))}\over{\omega_n\rho^n}}.$$
In the case that the two quantities are equal we call the common quantity the $n$-dimensional $\mu$-density of $A$ at $x$ 
denoted by $\Theta^n(\mu,A,x)$.
\end{def1} \noindent
{\bf Remark:} \newline 
Depending on which quantities are from the context understood, we will also use the terms density of $A$ at $x$ or 
simply the density at $x$. \newline 
The $\sigma$-algebra $\cl{B}$here is mentioned for formality but is unimportant in the definition.
\newline \newline
Also fundamental to our considerations is the concept of rectifiability. We will need several 
forms of the definition of rectifiability. Their equivalences are well presented in $\cite{simon1}$. We shall not 
here be interested in general rectifiable sets, so we restrict ourselves immediately to countably rectifiable sets.
Firstly and most basically we have the following definition.
\begin{def1}\label{defrec1}\thst
A set $M \subset \R^{n+k}$ is said to be countably $n$-rectifiable if 
$$M \subset M_0 \cup \bigcup_{j=1}^{\infty}F_j(\R^n)$$
where $F_j:\R^n \rightarrow \R^{n+k}$ are Lipschitz functions and $\Hm{n}(M_0)=0$
\end{def1}\noindent
{\bf Remark} By standard Lipschitz extension results we know that we can also write 
$$M = M_0 \cup \bigcup_{j=1}^{\infty}F_j(A_j)$$ 
for subsets $A_j \subset \R^n$.
\newline \newline
Notice that we have not required that the sets by measurable, which is occasionally required in definitions of 
rectifiable sets. It is however not necessary since, as we will see, all of the relevant sets we will be considering 
are in any case measurable since they can be shown to be expressable as countable unions and intersections of Borel sets 
in the appropriate Euclidean space. 
\newline \newline
From this basic definition it is known that the following expression for rectifiable sets holds.
\begin{lem}\label{lemrec1}\thst
$M$ is countably $n$-rectifiable if and only if 
$$M \subset \bigcup_{j=0}^{\infty}N_j,$$
where $\Hm{N}(N_0)=0$ and where each $N_j$, $j \geq 1$, is an $n$-dimensional embedded $C^1$ submanifold of $\R^{n+k}$.
\end{lem}\noindent
To introduce the final representation that we need we first need the following definitions.
\begin{def1}\label{def6} \thst 
We let the time fixed blow-up function be denoted by $\eta$, that is for any subeset $A \subset \R^n$
$$\eta_{y,\rho}(A) = \rho^{-1}(A-y).$$
Let $L$ be an subspace of $\R^n$ and $\rho \in \R$, $\rho>0$, then 
$$L^{\rho} = \{x \in \R^n : |x-y|<\rho \hbox{ for some } y \in L\}.$$
\end{def1} \noindent
and
\begin{def1}\label{defats}\thst
If $M$ is an $\Hm{n}$-measurable subset of $\R^{n+k}$ and $\theta$ is a positive locally $\Hm{n}$-integrable function on 
$M$, then we say that a given $n$-dimensional subspace $P$ of $\R^{n+k}$ is the approximate tangent space for $M$ with 
respect to $\theta$ if 
\begin{eqnarray}
\lim_{\lambda \rightarrow 0}\int_{\eta_{x,\lambda}M}f(y)\theta(x+\lambda y)d\Hm{n}(y) 
& := & \lim_{\lambda \rightarrow 0}\lambda^{-n}\int_Mf(\lambda^{-1}(z-x))\theta(z)d\Hm{n}(z) \nonumber \\
& = & \theta(x)\int_Pf(y)d\Hm{n}(y) \nonumber
\end{eqnarray}
for all $f \in C_C^0(\R^{n+k})$. The function $\theta$ is called the multiplicity function of $M$.
\end{def1}\noindent
We will in general consider sets with the multiplicity function set to 1.
\newline \newline
Our final definition of countably $n$-rectifiable sets is now stated in the form of the following theorem. 
\begin{thm}\label{thmrec1}\thst
Suppose $M$ is $\Hm{n}$-measurable. Then $M$ is countably $n$-rectifiable if and only if there is a positive locally 
$\Hm{n}$-integrable function $\theta$ on $M$ with respect to which the approximate tangent space $T_XM$ exists for 
$\Hm{n}$-a.e. $x \in M$.
\end{thm}\noindent
{\bf Remark:} We note that, for example in $\cite{simon2}$ it is often required that the total or 
$\Hm{n}$ measure of a set $M$ or at least $\Hm{n}(M \cap K)$ be finite for each compact set $K$. We do not, a priori 
make this assumption.
\newline \newline
Rectifiability can be seen as the weakest form of structure that a set can posses. However, we can explore parts of 
even unrectifiable sets in the case that they contain rectifiable parts. This fact will be useful to us. Particularly 
in chapter 4. For these reason we also define purely unrectifiable sets.
\begin{def1}\label{defpurunrec}\thst
A set $P$ is said to be purely $n$-unrectifiable if it contains no countably $n$-rectifiable subsets of $\Hm{n}$ positive 
measure.
\end{def1}\noindent
We note to this definition that for any set in $\R^{n+k}$, $A$, $A$ can always be decomposed into the 
disjoint union of two sets $A = R \cup P$ where $R$ is countably $R$ rectifiable and $P$ is purely $n$-unrectifiable.
\section{Motivation of the Classification}
We move now onto the motivation and construction of the problem at hand, previous results and results that follow 
more or less trivially from the literature. \newline \newline
An additional motivation to that mentioned in the introduction to this work 
was to perhaps uncover a way to attack the local $\Hm{j}$-finality of 
singularity sets for minimal surfaces or surfaces moving by their mean curvature. This is supported by the mentioned 
results in Leon Simons' 
$\cite{simon2}$ paper on the rectifiability of minimal surfaces, and recent work by Huisken and Sinistrari that 
shows that estimates on the shape of singularity sets is heading in the direction of satisfying the properties 
of the definitions under consideration. In particular, in Simon $\cite{simon2}$ 
a Lemma (the same one as has been previously discussed) shows that at least parts of the singularity
sets of particular types of minimal surfaces exactly satisfy one of the approximation properties. 
\newline \newline
We state this Lemma 
(after appropriate definitions) as a motivational starting point and also as it highlights some of the interesting 
points of the results. We then state the definitions mentioned in the introduction that 
we wish to classify and provide more fully a discussion of what it is we want 
to classify in these definitions. We also provide here a summary of the classification that is the central classification 
of the work. 
\begin{def1}\label{def2} \thst
By a {\bf Multiplicity one class} of minimal surfaces, 
$\M$, we will mean a set of smooth (i.e. infinitely differentiable) $n$-dimensional minimal submanifolds. 
Each $M \in \M$ is assumed to be properly embedded in $\R^{n+k}$ in the sense that for each 
$x \in M$ there is a $\sigma >0$ such that $M \cap \overline{B}_{\sigma}(x)$ is a compact connected embedded smooth 
manifold 
with boundary contained in $\partial B_{\sigma}(x)$. We also assume that for each $M \in \M$ there is a corresponding 
open set $U_M \supset  M$ such that $\Hm{n}(M \cap K) < \infty$ for each compact $K \subset U_M$, and such that $M$ is 
stationary in $U_M$ in the sense that
$$\int_M div_M\Phi d\mu = 0.$$
whenever $\Phi = (\Phi^1, ..., \Phi^{n+k}):U_M \rightarrow \R^{n+k}$ is a $C^{\infty}$ vector field with compact support 
in $U_m$. Where we have used $\mu = \Hm{n}|_M$.
We also require that the multiplicity one class of submanifolds are closed with respect to 
sequential compactness, orthogonal transformations and homotheties, that is:
\begin{enumerate}
\item $M \in \M \Rightarrow q \circ \eta_{x, \rho}M \in \M$ and $q \circ \eta_{x, \rho}U_M = U_{q \circ \eta_{x,\rho}M}$ 
for each $\rho \in (0,1]$, and for each orthogonal transformation $q$ of $\R^{n+k}$.
\item If $\{M_j\}_j \subset \M$, $U \subset \R^{n+k}$ with $U \subset U_M$, for all sufficiently large $j$, and 
$\sup_{j\geq 1}\Hm{n}(M_j \cap K) < \infty$ for each compact $K \subset U$, then there is a subsequence $M_{j^{\prime}}$ 
and an $M \in \M$ such that $U_M \supset U$ and $M_{j^{\prime}} \rightarrow M$ in $U$ in the sense that 
$$\int_{M_{j^{\prime}}}f d\Hm{n} \rightarrow \int_Mfd\Hm{n}$$ for any $f \in C^0_C(U, \R)$.
\end{enumerate}
\end{def1} \noindent
We assume here that the $M \in \M$ have no removable singularities: thus if $x \in \overline{M} \cap U_M$ and there is 
a $\sigma >0$ such that $\overline{M} \cap \overline{B}_{\sigma}(z)$ is a smooth connected embedded $n$ -dimensional
submanifold with boundry contained in $\partial B_{\sigma}(z)$, then $z \in M$. Subject to this agreement we 
can make the following definition 
\begin{def1} \label{def3} \thst
Suppose that $\M$ is as above and that $M \in \M$ then the {\bf (interior) singular set of $M$} 
(relative to $U_M$) is defined by 
$$\hbox{sing}M=U_M \cap \overline{M}\sim M$$
and the {\it regular set of $M$} is simply $M$ itself, that is
$$\hbox{reg}M=M.$$
\end{def1} \noindent
With these definitions we can now state the motivating Lemma due to Simon $\cite{simon2}$:
\begin{lem} \label{lem1} \thst
If $\M$ is a multiplicity one class of minimal surfaces, $M \in \M$, 
$$m := \max\{dim singM: M \in \M\}$$
$$z_0 \in singM$$
and 
$$S_+(z_0):=\{z \in M : \Theta^m(M,z) \geq \Theta^m(M,z_0)\}$$
Then for each $\e > 0$ there is a $p = p(\e,z_0,M)>0$ such that $S_+(z_0)$ has the following approximation property 
in $\overline{B_{p}(z_0)}$:
\newline \newline
For each $\sigma \in (0,p]$ and $z \in S_+(z_0) \cap \overline{B_p(z_0)}$ there is an $m$-dimensional affine 
subspace $L_{z,\sigma}$ containing $z$  with 
$$S_+ \cap B_{\sigma}(z) \subset \hbox{ the } (\e\sigma)-\hbox{nhood of } L_{z,\sigma}.$$
\end{lem}
\noindent
We note that in the case of Mean Curvature Flows, the singularity set can also be defined as follows:
\begin{def1}\label{def4} \thst
We say that a solution of Mean Curvature Flow $(M_t)_{t<t_0}$ reaches $x_0 \in \R^{n+1}$ at time $t_0$ if 
there exists a sequence $(x_j,t_j)$ with $t_j\nearrow t_0$ so that $x_j \in M_j$ and $x_j \rightarrow x_0$.
\end{def1} \noindent
\begin{def1}\label{def5} \thst
Let $\M = (M_t)$ be a smooth solution of mean curvature flow in $U \times (t_1, t_0)$. We say that $x_0 \in U$ is 
a {\bf singular point} of the solution at time $t_0$ if $\M$ reaches $x_0$ at time $t_0$ and has no smooth extension 
beyond time $t_0$ in any neighbourhood of $x_0$. All other points are called {\bf regular points}. The {\bf singular set} 
at time $t_0$ will be denoted by $sing_{t_0}\M$ and the regular set by $reg_{t_0}\M$.
\end{def1} \noindent  
As singularity sets are the motivation rather than the subject of our investigation, the properties of singular sets 
are used very little. 
However, in determining how applicable our results may be to singular sets we find that it is important to 
note that singular sets (from either definintion) are closed.
\begin{prop}\label{prop1} \thst
Singular sets as defined in either Definition $\ref{def3}$ or Definition $\ref{def5}$ are closed. 
\pf 
Suppose that the statement is not true, then there is a point $x \in regM$ such that for all $r > 0$ 
$B_r(x) \cap singM \not= \emptyset$. In particular since $x \in reg M$ there is a radius $\rho_x >0$ such that 
$\bar{M} \cap \overline{B_{\rho_x}(x)}$ 
is "smooth" (either in the infinitely differentiable in space time sense for mean curvature flow, 
or the sense outlined in Definition $\ref{def2}$, depending on whether we are proving the result for Definition 
$\ref{def3}$ or $\ref{def5}$) and such that 
$B_{\rho_x} \cap sing M \not=0$. Thus there is a $z \in singM$ and $\rho_z >0$ such that 
$B_{\rho_z}(z) \subset B_{\rho_x}(x)$. It follows that $\bar{M} \cap \bar{B_{rho_z}(z)}$ is "smooth" and thus $z \in regM$.
This contradiction shows such a point $x$ cannot be found which completes the proof. 
\end{proof}
\end{prop} \noindent
We now construct the properties that we will be investigating. 
We will always be considering sets being approximated by $j$-dimensional affine spaces that are subspaces of $\R^n$. 
We will identify $\R \times \{0\}$ with $\R$ and denote the projection onto $\R$ by $\pi_x$. Further, if $L$ is a 
$1$-dimensional affine space in $\R^2$ we will denote the projection onto $L$ by $\pi_L$.
\newtheorem*{defA}{Definition A}
\begin{defA}\label{defa} \thst
Let $A \subset \R^n$ be an arbitrary set and $\delta > 0$; then
\newline \newline
(i) $A$ has the weak $j$-dimensional $\delta$-approximation property if for all $y \in A$ there is $\rho_y > 0$ such that 
for all $\rho \in (0,\rho_y]$, $B_{\rho}(y) \cap A \subset$ the $\delta \rho$-neighbourhood of some $j$-dimensional
affine space $L_{y, \rho}$ containing $y$.
\newline \newline
(ii) $A$ has the weak $j$-dimensional $\delta$-approximation property with local $\rho_y$-uniformity if for all 
$y \in A$ there is a $\rho_y > 0$ such that for all $\rho \in (0,\rho_y]$ and all $x \in B_{\rho_y}(y)$
$$B_{\rho_y}(x) \cap A \subset L_{x,\rho}^{\delta \rho}$$
for some $j$-dimensional affine space $L_{x,\rho}$.
\newline \newline
(iii) $A$ is said to have the fine weak $j$-dimensional $\delta$-approximation property if for all $\delta >0$ $A$ 
has the weak $j$-dimensional $\delta$-approximation property with respect to $\delta$.
\newline \newline
(iv) $A$ has the fine weak $j$-dimensional approximation property with local $\rho_y$-uniformity if $A$ satisfies $(ii)$ 
for all $\delta > 0$.
\newline \newline
(v) The property (i) is said to be $\rho_0$-uniform, if $A$ is contained in some ball of radius $\rho_0$ and if, for 
every $y \in A$ and every $\rho \in (0,\rho_0]$, $B_{\rho}(y) \cap A \subset$ the $\delta \rho$-nhood of some 
$j$-dimensional affine space $L_{y,\rho}$ containing $y$.
\newline \newline
(vi) $A$ has the strong $j$-dimensional $\delta$-approximation property if for each $y \in A$ there is a $j$-dimensional 
affine space $L_y$ containing $y$ such that the definition (i) holds with $L_{y, \rho} = L_y$ for every 
$\rho \in (0,\rho_y]$.
\newline \newline
(vii) $A$ has the strong $j$-dimensional $\delta$-approximation property with local $\rho_y$-uniformity 
if for all $y \in A$ there exists a $\rho_y > 0$ 
and an affine space $L_y$ such that for all $x \in B_{\rho_y}(y)$ and all $\rho \in (0,\rho_y]$
$$B_{\rho}(x) \cap A \subset L_y^{\delta \rho}.$$
(viii) The property in (vi) is said to be $\rho_0$-uniform if $A$ is contained in some ball of radius $\rho_0$ and if 
for each $y \in A$ there is a $j$-dimensional affine space $L_y$ containing $y$ such that 
$B_{\rho}(y) \cap A \subset$ the $\delta \rho$-nhood of $L_y$ for each $\rho \in (0,\rho_0]$.
\newline \newline
Due to the long names of the properties, they will be henceforth referred to only by their number.
\end{defA} \noindent
Our classification is to gat a simple yes or no answer for each of the eight definitions with respect to two questions.
\begin{ques}\label{ques1}\thst 
We wish to classify the definitions in Definition 1 with respect to the following questions: 
\begin{enumerate}
\item if the set will be of dimension $j$ (or rather $\leq j$), and 
\item if the set will have some locally finite Hausdorff measure property.
\end{enumerate}
\end{ques} \noindent
With these questions in mind we will concern ourselves with asking about the answer to (1) or (2) with respect to 
a certain definition, for example the answer to (i) (1) is no.
\newline \newline
As we are generally probing here for 'free information' about singularity sets, and the use of more than 
one definition of the terms about which we are asking in the literature we remain open as to 
which definition it is that we are making classifications with respect to.
We therefore allow 
for two strengths of locally finite $\Hm{j}$ measure. In only one case do find that the answer as to 
possesing locally finite $\Hm{j}$ measure is affected by the choice of strength of definition, that is 
for $(vii)$ where the definition ensures satisfaction of the weaker but not the stronger. The 
definitions are:
\begin{def1}\label{def7}\thst
A subset $A \subset \R^n$ is said to have {\bf locally finite $\Hm{j}$ measure} 
(or local $\Hm{j}$-finality) if for all compact subsets $K \subset \R^n$,
$$\Hm{j}(K \cap A) <\infty ,$$
or equivalently, if for all $y \in \R^n$ there exists a radius $\rho_y > 0$ such that 
$$\Hm{j}(B_{\rho_y}(y) \cap A) <\infty.$$
A subset is said to have {\bf weak locally finite $\Hm{j}$ measure}
(or weak local $\Hm{j}$-finality) if for each $y \in A$ there exists a radius 
$\rho_y > 0$ such that 
$$\Hm{j}(B_{\rho_y}(y) \cap A ) < \infty.$$
\end{def1} \noindent
An example of the difference is that 
$$\cl{N}:=\bigcup_{n=1}^{\infty}\R \times\left\{{{1}\over{n}}\right\}$$
has weak local $\Hm{j}$-finality but not local $\Hm{j}$-finality. The use of allowing the weak definition is that in 
some cases, such as the one just given a set with weak local $\Hm{j}$-finality will be the finite union of a collection 
of sets with local $\Hm{j}$-finality. Which still counld be understood as having reasonably behaved local measure 
when the structure giving the locally infinite measure is known.
\newline \newline
As we will see, and has been hinted at, we do not necessarily get very much information for free. Paritulcarly as 
we get a "no" to answer the definition corresponding Simons' Lemma. However, as mentioned in the introduction. In 
this case we do show that in order for something to go wrong the set does have to be truly badly behaved which should 
be helpful. We now note formally that the condition in Simons' Lemma is definition (iv). 
\begin{prop}\label{prop2}\thst
The $S_+(z_0)$ sets introduced in Lemma $\ref{lem1}$ are (iv).
\pf
Direct comparison between the property shown in Lemma $\ref{lem1}$ and (iv) shows that this is exactly what is shown in 
Lemma $\ref{lem1}$.
\end{proof}
\end{prop} \noindent
\section{Results Following from the Literature}
Although the problem we are looking at has not previously been systematically investigated, a few of the 
results follow easily from results already in the literature for which 
proofs can be found, for example in Simon $\cite{simon3}$. 
Excepting a counter example, the relevant results can be convieniently stated in the following Lemma.
\begin{lem}\label{lem2}\thst
(i) There is a function $\beta:[0,\infty) \rightarrow [0,\infty)$ with $\lim_{\delta \searrow 0}\beta(\delta)=0$ 
such that if $A \subset \R^n$ has the $j$-dimensional weak $\delta$-approximation property for some given 
$\delta \in (0,1]$, then $\Hm{j+\beta(\delta)}(A)=0$. (In particular if $A$ has the $j$-dimensional weak $\delta$-
approximation property for each $\delta > 0$, then dim$A \leq j$.)
\newline \newline
(ii)If $A \subset \R^n$ has the strong $j$-dimensional $\delta$-approximation property for some $\delta \in (0,1]$, 
then $A \subset \cup_{k=1}^{\infty}G_k$, where each $G_k$ is the graph of some Lipschitz function over some 
$j$-dimensional subspace of $\R^n$.
\newline \newline
(iii) If $A \subset \R^n$ has the $\rho_0$ uniform 
strong $j$-dimensional $\delta$-approximation property for some $\delta \in (0,1]$, 
then $A \subset \cup_{k=1}^QG_k$, where $G_k$ is the graph of some Lipschitz function over some $j$-dimensional subspace 
of $\R^n$, $L$.
\end{lem}\noindent
We show in the following Corollary that the above Lemma allows us to answer yes to properties  
(vi) (1), (viii) (1) and (2), (iii) (1), (iv) (1) and (vii) (1) and (2), although we answer 
yes to (vii) (2) only with weak local $\Hm{j}$-finality, to local $\Hm{j}$-finality we answer no.
\begin{cor}\label{cor1}\thst
The answer to the following Definitions is yes:
\begin{list}{}{}
\item (vi) (1), 
\item (viii) (1),
\item (viii) (2),
\item (iii) (1),
\item (iv) (1),
\item (vii) (1), and 
\item (vii) (2)
\end{list}
\pf
(iii) (1) follows from Lemma $\ref{lem2}$ $(i)$ isince 
\newline \newline
"{\it In particular if $A$ has the $j$-dimensional weak $\delta$-
approximation property for each $\delta > 0$, then dim$A \leq j$.}"
\newline \newline
means that should $A$ satisfy (iii), then dim$A \leq j$ which proves that the answer to 
(iii) (1) is yes. Further, since (iv) (1) is a strengthening 
of (iii), sets satisfying the properties of (iv) must further 
satisfy any properties following from sets satisfying (iii), thus the answer to 
(iv) (1) must also be yes.\newline \newline
Any graph of a Lipschitz function over a $j$-dimensional affine space clearly has dimension 
less than or equal to $j$. It follows then 
that any countable union of such graphs will also have dimension bounded above by $j$. It thus follows from Lemma 
$\ref{lem2}$ (ii) and (iii) that the answers to (vi) (1) and (viii) (1) are yes. Similarly to the 
preceeding paragraph, the fact that (vii) is a strengthening of (vi) that the 
answer to (vii) (1) is yes. 
\newline \newline
Further concerning (viii), suppose that we have a set $A$ satifying 
the conditions of property $(viii)$. Suppose also that $ x \in \R^n$ and $\rho >0$. Then we know that 
$$A \cap B_{\rho}(x) \subset \bigcup_{k=1}^Qg_k(\pi_{L_k}(B_{\rho}(x)))$$
where $L_k$ are the $j$ dimensional affine spaces that Lemma $\ref{lem2}$ ensures exist and the $g_k$ are the Lipschitz 
functions over the $L_k$ that combined contain $A$. Thus 
$$\Hm{j}(A \cap B_{\rho}(x)) \leq \sum_{k=1}^Q\Hm{j}(g_{k}(\pi_{L_k}(B_{\rho}(x)))).$$
Since $card(\{g_k\}_{k=1}^Q)=Q < \infty$ there exists a 
$$M = \max_k \hbox{ Lip}g_k < \infty$$
so that by the Area formula 
\begin{eqnarray}
\Hm{j}(A \cap B_{\rho}(x)) & \leq & \sum_{k=1}^Q\Hm{j}(g_{k}(\pi_{L_k}(B_{\rho}(x)))) \nonumber \\
& \leq & \sum_{k=1}^Q M\Hm{j}(\pi_{L_k}(B_{\rho}(x))) \nonumber \\
& = & \sum_{k=1}^QM\omega_j \rho^j \nonumber \\
& = & QM\omega_j\rho_j \nonumber \\
& < & \infty. \nonumber
\end{eqnarray}
We thus have that property $(viii)$ does ensure locally finite measure, and thus we have shown that the answer to 
(viii) (2) is yes.
\newline \newline
Finally we note that should we have a set satisfying (vii), then, by definition, for each $y \in A$ 
there is a $\rho_y > 0$ and an affine space $L_y$ such that for all $x \in B_{\rho_y}(y)$ and all $\rho \in (0,\rho_y]$
$$B_{\rho}(x) \cap A \subset L_y^{\delta \rho}.$$
It follow that $B_{\rho_y}(y) \cap A$ satisfies (viii), thus $\Hm{j}(K \cap A) < \infty$ for each 
compact $K \subset \R^n$, That is 
\begin{eqnarray}
\Hm{j}(B_{\rho_y}(y) \cap A) & \leq & \Hm{j}(\bar{B_{\rho_y}(y)} \cap A) \nonumber \\
& < & \infty. \nonumber
\end{eqnarray}
Thus giving weak local $\Hm{j}$-finality, and thus allowing us to answer (vii) (2) with yes.
\end{proof}
\end{cor} \noindent
{\bf Remark: }
We note that the proof as written is also optimal in that we cannot get better than weak local $\Hm{j}$-finality 
for (vii) as seen in the already given example of $\cl{N}$. For each $y \in \cl{N}$ we can 
find a $\rho_y > 0$ such that $B_{\rho_y}(y) \cap \cl{N} \subset \R \times \{1/n\}$ for some $n \in \mathbb{N}$, 
and by setting $L_y$ as this affine space for each $y$ it is clear that $\cl{N}$ satisfies (vii), 
However, for each $r>0$
$$\Hm{j}(B_r((0,0)) \cap \cl{N}) = \infty$$
so that $\cl{N}$ does not have locally finite $\Hm{j}$-measure.
\newline \newline
Another contribution that comes from Simon $\cite{simon3}$ is a set that is similar in form to the main and 
most interesting counter example that is presented here. Its actual construction and properties will be 
discussed in the following section, however, in noting results that have already been essentially shown, we 
acknowledge its existence and that it was known to satisfy one of the definitions.
\begin{lem}\label{lem3}\thst
There is a set, $\Gamma_{\e}$ that satisfies (i) for $j = 1$ that has dimension greater than~$1$.
\end{lem}\noindent
In later chapters dimension of $\Gamma_{\e}$ and related sets will be discussed. The original proofs that we present 
will be based on the knowledge of how to calculate the dimension of $\G$. The proof of the relevant formula will, 
however, not be presented, as it also already exists in the literature. The proof can be found in $\cite{hutch}$.
\begin{cor}\label{cor2}\thst
The answer to (i) (1) and (i) (2) is no.
\pf
The set $\Gamma_{\e}$ of Lemma 3 constructed in the following section provides a counter example to the answer 
to (i) (1) being yes. Since the dimension of a set satisfying (i) with $j=1$ could be 
greater than $1$, there is clearly no gaurantee of any form of finite $\Hm{1} = \Hm{j}$ measure. Thus the answer 
to (i) (2) is also no.
\end{proof}
\end{cor}\noindent
This completes the survey of the results that were already known, or rather, at least already almost known. 
So that the complete classification of all the definitions is presented in a convenient easily digested way somewhere 
we complete this chapter with a table of the complete classification that we prove in this thesis.
\newline \newline
The classification of the definitions in Definition 1 with respect to the questions presented in Questions 1 is as 
follows:
\begin{eqnarray}
(i) & &(1) \ \  no \ \ (2) \ \  no \nonumber \\
(ii) & &(1)  \ \ no \ \  (2) \ \  no \nonumber \\
(iii) & &(1) \ \  yes \ \  (2) \ \ no \nonumber \\
(iv) & &(1) \ \ yes \ \ (2) \ \ no \nonumber \\
(v) & &(1) \ \ no \ \  (2) \ \ no \nonumber \\
(vi) & & (1) \ \ yes \ \  (2) \ \  no \nonumber \\
(vii) & & (1) \ \  yes \ \  (2) \ \  yes(weak) / no(strong)  \nonumber \\
(viii) & & (1) \ \  yes \ \ (2) \ \ yes. \nonumber  
\end{eqnarray}
We note that those definitions classified as yes have all already been answered. It remains only to show that the 
classification of the remaining definitions is no.

\chapter{Construction of the Counter Examples}
Having answered all the questions that will be answered with yes, we now turn our attention to providing counter 
examples for the remaining questions so as to answer no to each of these. For those with a didactic turn of mind, 
of course these counter examples were constucted in association 
with answering our questions and not constructed before hand, only to be quite coincidently successfully used later.
\newline \newline
The sets being considered are not all trivial sets to construct or to understand. At least not at first sight. We 
therefore provide only the constructions and some intrinsic properties of the sets, leaving the proofs that 
they actually satisfy the definitions that they are respectively intended to be counterexamples to until later. 
For the more complicated sets, particularly $\A$, there is more than one method to construct the set. Some of 
these will be discussed further in Chapters 7 and 8. For now, however, we satisfy ourselves with the definitions 
most easily used to fit the constructed sets to the relevant definitions and thus complete the classification.
\newline \newline
In this chapter we construct 3 sets and 3 1-parameter families of sets. Of the latter three the first 
is our own construction of a known set, the same that appears in Lemma 3, which we provide since the necessary 
properties for our purposes 
are more easily proven with our construction method. The latter two are then variations of the 
same set allowing for important extra properties by adding another point of variation.
For the sets with a variable there is a range of values of the parameter (independent of which set) for which 
each resultant specific example is appropriate for our purposes. We will, however, calculate with the 
parameter left arbitrary since it provides more generality and makes no difference to the proofs of the results 
that we want to prove with these sets. 
\newline \newline
The three simpler sets are of little interest apart from the fact that 
they are appropriate counter examples to particular definitions. The other three are of independent interest. 
As well as allowing us to show that some good behaviour is ensured by the approximate $j$-dimensionality of the sets 
if not as definite as we had hoped, they provide a range of interesting results on dimension, rectifiability and measure 
density. Alot of the general proofs concerning properties of these sets are included 
in the discussion of generalised Koch Type sets 
(the generalisations of these three 'specific' examples) in Chapters 7 and 8. We include in any case the direct proofs 
of the properties that we are interested that are relevant to the classification work. That is we include direct proofs 
that for each definition for which there is a counter example there is a closed counter example (important, since as 
we have shown in Proposition $\ref{prop1}$ singularity sets are all closed) and that the sets of integer dimension are 
shown directly to have their respective dimensions. 
\newline \newline
We construct firstly the three simpler sets. We then construct $\G$ which will be a counter example to (i) (1) followed 
by a property of $\G$ important to our study. We then construct the second more complicated set $\A$ which is a 
counter example to (iv) (2). Since $\A$ is not closed and is therefore not possibly a singularity set we make 
the third construction $\Ae$, which is a subset of the second, constructed to be closed but retain the necessary properties. 
We then prove some necessary proprties of $\A$ and $\Ae$.
\section{Simple and Known Sets}
The first set has already been defined, and is:
$$\cl{N}:=\bigcup_{n=1}^{\infty}\R \times\left\{{{1}\over{n}}\right\}$$
Note that we will henceforth identify $R^n \times [0]^{N-n}$ with $\R^n$ in $\R^N$ for each choice of $n,N \in \mathbb{N}$ 
with $n <N$.
The other simple ones, are used in a similar way to $\cl{N}$ but need differing levels of fineness approximation with 
bad properties at one point. Being a collection of flat sheets, $\cl{N}$ does not have this property, 
we therefore define the subset of $\R^2$ defined for each $\delta > 0$ as 
$$\Lambda_{\delta} = \bigcup_{n=1}^{\infty} \hbox{graph}\left\{{{\delta x}\over{n}}\right\} 
\cup \bigcup_{n=1}^{\infty} \hbox{graph}\left\{{{-\delta x}\over{n}}\right\}$$
and the subset of $\R^2$ defined as
$$\Lambda^2 = \bigcup_{n=1}^{\infty} \hbox{graph}\left\{{{x^2}\over{n}}\right\} 
\cup \bigcup_{n=1}^{\infty} \hbox{graph}\left\{{{-x^2}\over{n}}\right\}.$$
\newline \newline
We now construct the more complicated examples, they are both based on the "Koch Curve" which was originally 
constructed as a fractal set being of dimension between $1$ and $2$. The first we construct is the set $\Gamma$ 
given by Simon in $\cite{simon3}$, on which the remaining sets are based. The second set, which is 
actually a function from $\R^+$ into $2^{\R^n}$ (that is, the set is constructed with respect to a variable 
$\varepsilon \in \R^+$) will be denoted $A_{\varepsilon}$, and is used as a counter example 
to (iv) (2). Although $\G$ was actually constructed as a fixed set, we will allow the set to 
be constructed with respect to a variable $\varepsilon$, which will later allow us to find appropriate counter examples 
with respect to (i) (1) for any given $\delta$. The variable set will then be denoted 
$\Gamma_{\varepsilon}$ 
\newline \newline
These constructions rely heavily on the use of triangles so we first make the following definition.
\begin{def1}\label{def8}\thst
Let $L = (a,b) = ((a_1,a_2),(b_1,b_2))$ be a line in $\R^2$. A {\bf $\varepsilon$-triangular cap} or, when 
the context is clear, simply a {\bf cap} will be the triangle , $T$, with vertices $a,b$ and $c + (a+b)/2$  
(we write $c$ also as $(c_1,c_2)$), where $c$ is chosen such that 
\begin{list}{}{}
\item $|c| = \varepsilon$ and
\item $<c,b-a> = 0$.
\end{list}
Further to ensure that the cap is well defined we choose $c$ from the two remaining possible points in $\R^2$ as follows. 
Should $L$ be an edge of a previously constructed triangular cap, $T_0$, then $c$ is chosen such that $T \subset T_0$.
Otherwise, if $c_1 \not\equiv (a_1 +b_1)/2$ (regardless of which of the two possbilities) then $c$ is 
chosen such that  $c_1 > (a_1 +b_1)/2$, otherwise we choose $c$ such that $c_2 > (a_2+b_2)/2$.
\end{def1} \noindent
\begin{cons}\label{cons1}
\begin{proof}[] \thst
We construct the set $\Gamma_{\e}$ as follows.
\newline \newline
Let $\varepsilon > 0$. 
We begin with a $\varepsilon$-triangular cap, $T_0$, constructed over the line $A_{0,1} := ((0,0),(1,0))$. We then 
name the two new edges $A_{1,j}$, $j = 1,2$. We denote the first "approximation", which is $T_0$, as $A_0$.
We note that $l:=\Hm{1}(A_{1,j}) < \Hm{1}(A_{0,1})$, $j=1,2$. 
We then construct $l\varepsilon$-triangular caps $T_{1,j}$ on $A_{1,j}$. We name the four new edges 
$A_{2,j}$, $j \in \{1,2,3,4\}$. We denote this second "approximation", $\cup_{j=1}^2T_{1,j}$ by $A_1$. We note 
that $A_{2,j}$, $j=1,...,4$ are the $2^2$ shortest edges of length $l^2$. We note also that $A_1$ can also be 
constructed by the appropriately rotated union of two copies of $A_0$
\newline \newline
We now continue inductively, suppose that we have a set $A_n$ consisting of $2^n$ triangular caps, $T_{n,j}$ 
with base length $l^{n-1}$ and altogether $2^{n+1}$ "shortest sides", $A_{n+1,j}$ of length $l^{n+1}$. On 
each $A_{n+1,j}$ we construct a $l^{n+1}\varepsilon$-triangular cap $T_{n+1,j}$. We set 
$$A_{n+1} := \bigcup_{j=1}^{2^{n+1}}T_{n+1,j}.$$
This $A_{n+1}$ will then have all of the same properties as $A_n$ with $n$ replaced by $n+1$.
We note also that with the numbering of the caps, we always count from "left" to "right" so that 
$T_{n+1,2j-1} \cup T_{n+1,2j}\subset T_{n,j}$.
\newline \newline
We then define 
$$\Gamma_{\varepsilon} = \bigcap_{n=0}^{\infty}A_n$$
where the dependence on $\varepsilon$ comes from the initial choice of $\varepsilon$. 
\end{proof}
\end{cons} \noindent
One property of $\Gamma_{\varepsilon}$ that should be noted now, as it is particularly intrinsic to the construction 
is that $\G$ is essentially the union of two scaled copies of itself. We show this after the following definitions.
\begin{def1}\label{def9}\thst
We denote the end points of a line of finite length, $A$, as $E(A)$, and call them the 
{\bf edge points of} $A$. Let $T_{n,i}$ be a triangular cap. $T_{n,i}$ will then have 3 vertices which will be called the
{\bf edge points of} $T$. Let $A_n$ be a stage in Construction 1 or 2 (we will see that the definition applies to definition 
as well to Construction 2) then the {\bf edge points of} $A_n$ are
$$E(A_n) := \bigcup_{i=1}^{2^n}E(T_{n,i})$$
and the {\bf edge points of} $\Gamma_{\e}$ are 
$$E({\Gamma}_{\e}) := \bigcup_{n=1}^{\infty}E(A_n).$$
Also, as we will see, the same definition applies to $A_{\e}$. That is we can and do 
define the {\bf edge points of } $A_{\e}$ to be 
$$E(A_{\e}) := \bigcup_{n=1}^{\infty}E(A_n).$$
\end{def1} \noindent
We see the the edge points are all of the corners that appear in the constructions of $\Gamma_{\e}$ and $A_{\e}$.
\begin{def1}\label{def10}\thst
We define the {\bf edgeless} $\Gamma_{\e}$ as 
$$\Gamma_{\e}^E := \Gamma_{\e} \sim E(\Gamma_{\e}).$$
\end{def1} \noindent
\begin{prop}\label{prop3}\thst
There are contraction mappings, $S_1$ and $S_2$, and an open set, $O$, such that 
$$\Gamma_{\e}^E \subset O,$$
$$S_1(O) \cup S_2(O) \subset O,$$
$$S_1(\Gamma_{\e}^E) \cup S_2(\Gamma_{\e}^E) = \Gamma_{\e}^E$$
and
$$S_1(O) \cap S_2(O) = \emptyset.$$
Further
\begin{eqnarray}
Lip S_1 = Lip S_2 & = & l \nonumber \\
& := & (1/4 + \e^2)^{1/2} \nonumber
\end{eqnarray}
\pf
It is not too difficult to check that the contraction mappings of Lipschitz constants $l$ defined by 
\newline \newline
$$S_i(x,y) =  \left( \begin{array}{cc}
cos((-1)^itan^{-1}(\e)-\pi) & -sin((-1)^itan^{-1}(\e)-\pi) \\
sin((-1)^itan^{-1}(\e)-\pi) & cos((-1)^itan^{-1}(\e)-\pi)  \end{array} \right)v(x,y)$$
where
$$v(x,y) = \left( l \left( \left( \begin{array}{c} x \\ y \end{array} \right)-\left( \begin{array}{c}1/2 \\ \e /2 \end{array} \right)\right)
+\left( \begin{array}{c} 1/2 \\ \e /2 \end{array} \right) \right)$$
are such that 
$$S_1(T_0) = T_{1,2},$$ 
$$S_2(T_0) = T_{1,1}$$
and thus 
$$S_1(A_0) \cup S_2(A_0) = S_1(T_0) \cup S_2(T_0) = T_{1,1} \cup T_{1,2} = A_1.$$
Further, by setting $O$ to be the open quadrilateral with vertices 
$$\{(0,0),(1/2,3\e /2),(1,0),(1/2,-\e /2\}$$ 
we see 
$$\Gamma_{\e}^E \subset T_0 = A_0 \subset O$$ 
and that we have $S_1(O)$ is the quadrilateral of vertices 
$$\{(1/2,\e),(1,0),l((1,0)-(1/2,3\e /2))+(1/2, 3\e /2),l((1,0)-(1/2,-\e /2))+(1/2, -\e /2)\}$$
and $S_2(O)$ is the quadrilateral of vertices
$$\{(0,0),l(1/2,3\e /2),(1/2,\e),l(1/2,-\e /2)\}.$$
It follows that 
$$S_1(O) \cup S_2(O) \subset O$$ 
and 
$$S_1(O) \cap S_2(O) = \emptyset.$$
Note specifically that since the proceedure, $P$, of taking two triangular caps on the shorter sides of a union of 
isosceles triangles is clearly invariant under orthogonal transformation (since chosing the new cap to be 
within the previous triangle is independent of orientation) and homothety, that is $P(R(T)) = O(R(T))$ where $T$ is 
an isosceles triangles and $R$ is any orthogonal transformation on $\R^2$, and if 
$l \in \R$, $P(lT) = lP(T)$. Since $S_1$ and $S_2$ are indeed just combinations of homothety and orthogonal transformation 
we have $P(S_i(T)) = S_i(P(T))$ for $i=1,2$.
\newline \newline
We claim that 
$$A_n = S_1(A_{n-1}) \cup S_2(A_{n-1})$$ 
for each $n \in \mathbb{N}$. \newline \newline
We already have a starting point ($n=1$). Now, supposing that 
$$A_n = S_1(A_{n-1}) \cup S_2(A_{n-1})$$ 
for some $n \in \mathbb{N}$,
we then have
\begin{eqnarray}
A_{n+1} & = & P(A_n) \nonumber \\
& = & P(S(A_{n-1}) \cup S_2(A_{n-1}) \nonumber \\
& = & S_1(PA_{n-1}) \cup S_2(PA_{n-1}) \nonumber \\
& = & S_1(A_n) \cup S_2(A_n) \nonumber
\end{eqnarray}
Completing the induction. Then, since $A_1 \subset A_0$, we then have 
\begin{eqnarray}
\Gamma_{\e}^E & = & \bigcap_{n=0}^{\infty}A_n \sim E \nonumber \\
& = & \bigcap_{n=1}^{\infty}(A_n \sim E) \nonumber \\
& = & \bigcap_{n=1}^{\infty}S_1(A_{n-1} \sim E) \cup S_2(A_{n-1} \sim E) \nonumber \\
& = & \bigcap_{n=0}^{\infty}S_1(A_n \sim E) \cup S_2(A_n \sim E) \nonumber \\
& = & S_1\left(\bigcap_{n=1}^{\infty}A_n \sim E \right) \cup S_2\left(\bigcap_{n=1}^{\infty}A_n \sim E \right) \nonumber \\
& = & S_1(\Gamma_{\e}^E \sim E) \cup S_2(\Gamma_{\e}^E \sim E). \ \ \hbox{ } \nonumber
\end{eqnarray}
\end{proof}
\end{prop}\noindent
\section{Pseudo-Fractal Sets}
We now construct the "strangest" sets. These are similar to $\Gamma_{\varepsilon}$ in construction, however, 
as we noted in Proposition $\ref{prop3}$, the construction for $\Gamma_{\varepsilon}$ retains the basic shape 
of the triangular caps. This will not be sufficient for the cases when we want to prove properties for the 
case where approximations should hold for all $\delta > 0$. We therefore allow the relative height of the 
triangular caps to shrink, so that the "angles" involved in the triangles approach zero as we look at 
smaller and smaller sections of the triangles. As we will see later, even this adjustment is not 
sufficient. We therefore remove all of the interior at each stage, 
take, in a sense, a limit, remove the approximating sets and the edges. We make the specific constructions below in 
Constructions $\ref{cons2}$ and $\ref{cons3}$. Note that the heuristic path to our set just given was not the one 
that originally led 
to its construction, but rather, it is the result of being the embodiment of the worst case allowed in a 
failed attempt to prove that the answer to Definition A (iv) (2) was yes. As has been mentioned, the third 
set is then a carefully selected subset of this chosen in such a way as to ensure that it is closed.
\begin{cons}\label{cons2}
\begin{proof}[] \thst
We construct the set, as previously, as a subset of $\R^2$. We start with 
$$A_0 := [(0,0),(1,0)].$$
We then denote by $T_0$ the $2\varepsilon$-triangular cap on $A_0$.
\newline \newline
We now set
$$A_1:=\overline{(\partial T_{0,1} \sim A_0)},$$
which is the union of two lines (namely the two shorter edge lines of $T_{0,1}$), we name the two lines 
$A_{1,i}$, $i=1,2$. 
To continue, we denote by $T_{1,i}$ the $\varepsilon$-triangular cap constructed on $A_{1,i}$, 
considered as an edge of $T_{0,1}$ for each $i$. 
\newline \newline
We then set 
$$A_2 := \overline{\left(\partial \left(\bigcup_{i=1}^2T_{1,i}\right) \sim A_1\right)},$$
which will be a union of $4$ lines $A_{2,i}$, $i=1,2,3,4$. Each an edge of a triangular cap $T_{1,i}$.
\newline \newline
We continue the construction inductively.
Assuming we have $\cl{A}_n$, a union of $2^n$ lines, 
$\{A_{n,i}\}_{i=1}^{2^n}$ that lie on the boundary of 
$2^{n-1}$ triangular caps $\{T_{n-1,i}\}_{i=1}^{2^{n-1}}$, 
(and $A_n$ a union of $2^n$ triangular cups), we construct $2^n$ 
$2^{1-n}\varepsilon$-triangular caps, $\{T_{n,i}\}_{i=1}^{2^n}$,
on each of the $2^n$ lines. As previously we number from "left" to "right" so that 
$T_{n+1,2j-1} \cup T_{n+1,2j}\subset T_{n,j}$. We then set 
$$A_{n+1} := \overline{\left(\partial \left(\bigcup_{i=1}^{2^n}T_{n,i}\right) \sim A_n\right)},$$ 
Finally, we define  
\begin{eqnarray}
\A & := & \bigcap_{i=1}^{\infty}\left(\overline{\bigcup_{n=1}^{\infty}A_n} \sim \bigcup_{n=1}^{i}(A_n \sim E)\right) 
\sim E \nonumber \\
& = & \overline{\bigcup_{n=1}^{\infty}A_n} \sim \bigcup_{n=1}^{\infty}A_n, \nonumber
\end{eqnarray}
where
$$ E = \bigcup_{n=1}^{\infty}\bigcup_{i=1}^{2^n}E(A_{n,i}),$$
and $E(A_{n,i})$ denotes the endpoints of the line $A_{n,i}$. As previously, the $\varepsilon$ refers to the 
arbitraryily chosen $\varepsilon>0$ at the begining of the construction, which may, of course, be chosen, 
as small as is necessary. 
\end{proof}
\end{cons}\noindent
{\bf Remarks: } \newline
{\bf (1) } 
The removal of the endpoints is very important for the example. With the endpoints, there are of course points in 
the set with a fixed angle that must be squeezed into a $\delta$ approximation for exery $\delta>0$. This is 
not possible. With the endpoints missing we can, for each element of the set choose, for any given angle greater 
than zero, avoid all "corners" of angle greater than or equal to the given one, so as to make the set 
flater than the given angle in that neighbourhood. Since we are asking questions of measure, it is also important 
to note that the union of all the endpoints, that is 
$$\bigcup_{n=0}^{\infty}\bigcup_{i=1}^{2^n}E(A_{n,i}),$$
is countable and therefore of zero $\Hm{1}$ measure, thus having no effect on any $\Hm{1}$ measure properties that 
we are looking at. 
\newline \newline
{\bf (2) }
A remark on both Construction $\ref{cons1}$ and Construction $\ref{cons2}$ and indeed on 
Definition $\ref{def8}$ is that a triangular cap 
constructed on the edge of a previously constructed triangular cap may not be well defined in that it may 
not be a subset of the previous triangular cap. Another problem is that, as we often do, constructing 
triangular caps on both of the sides of identical length on an isoceles triangle may lead to the two 
new triangular caps indeed being subsets of the previous cap, but intersecting with each other. Choosing 
the vertical heights prevents this problem, and indeed, should the initial vertical heigth, $h$ be less than 
$1/4$ of the base length $b$, then provided the new vetical height is less than or equal to 
$h(\sqrt{h^2 + b^2/4})/b$ we will 
encounter no problems, such a proceedure cannot lead to a new vertical height being more than $1/4$ the base 
length, and furthermore in this situation, $h(\sqrt{h^2 + b^2/4})/b = l > h/2$ (where the $l$ is 
the scaling factor in Construction $\ref{cons1}$) so that no problems later in the induction can occur in 
Constructions $\ref{cons1}$, $\ref{cons2}$ and $\ref{cons3}$. 
We will always assume that the appropriate conditions on the vertical height have been satisfied. This is no problem 
as we want our caps to be very flat in any case.
\begin{def1}\label{def11}\thst
For each $n \in \mathbb{N} \cup \{0\}$ and each $i \in \{1,...,2^n\}$ there is a triangular cap $T_{n,i}$ 
constructed on $A_{n,i}$, we denote the vertice of $T_{n,i}$ that is not in $A_{n,i}$ (that is, the new vertice created) by
$a_{n,i}$.
\end{def1} \noindent
\begin{cons}\label{cons3}
\begin{proof}[]\thst
As previously mentioned we will be looking at a subset of $A_{\e}$. We have already noted that the edge points of 
$A_{\e}$ are countable, we now give them an ordering that will prove important later. We take 
$$e_1 = (0,0),$$
$$e_2 = (1,0),$$
$$e_3 = a_{0,1}$$
and then in general
$$e_{2+i+\sum_{j=0}^{n-1}2^j} =a_{n,i}.$$
We set 
$$\rho_1 = {{1}\over{4}}2^{-7}(1+7\cdot 16\e^2)^{1/2}$$
and for $n \in \mathbb{N}$, we set 
$$\rho_n = 4^{1-n}r_1 < 2^{-6-n}(1+(n+6)16\e^2)^{1/2}.$$
We now define a set of radii. We set 
$$r_1 = \rho_1$$
and 
$$r_2 = \rho_2.$$
Then for $i \geq 2$ there is a unique $n \in \mathbb{N} \cup \{0\}$ such that 
$$i \in \{2+\sum_{j=0}^n2^n, ..., 1+\sum_{j=0}^{n+1}2^j\}$$
so that we can define
$$r_i = \min\{\rho_i,d(e_i,A_{n-1} \cup A_{n-2})/2\}.$$
We then define 
$$\cal{B} := \bigcup_{\hbox{$i=1$}}^{\infty}\hbox{$B_{r_i}(e_i)$} $$
Note that 
$$E = \{e_i\}_{i=1}^{\infty} \subset \bigcup_{i=1}^{\infty}B_{r_i}(e_i) = \cal{B}.$$
Finally we define 
$$\Ae = A_{\e} \sim \cal{B}.$$
We note that this can also be written as 
$$\Ae = (A_{\e} \cup E) \sim \cal{B}.$$
\end{proof}
\end{cons}\noindent
There are three points concerning $A_{\e}$ and $\Ae$ that are important that should be noted. Firstly, the entire 
purpose of altering $A_{\e}$ to $\cal{A}_{\e}$ was that $\Ae$ should be closed. We therefore prove that this 
important property indeed holds. Secondly, although we will show that $A_{\e}$ and $\Ae$ have property (iv) 
with respect to $j=1$ and thus have dimension 1, the sets have some interesting properties in and of themselves. For 
this reason and as support for the consistency of the results here we provide a direct proof that the dimension of 
$A_{\e}$ and $\Ae$ is 1. Finally, as we will show in chapter 4,
the exotic counterexamples of $\A$ and $\Ae$ are necessary. 
Further, to support the idea that counter examples to (iv) 2 need necessarily be badly behaved, we note that 
$\A$, $\Ae$ are not rectifiable. As substantial preparation is necessary and since the fact is not necessary for 
our classification, we present the proof in Chapter 7 along with the generalisations of the sets. A direct 
proof for these specific examples is also given.
\begin{lem}\label{lem4} \thst
$\Ae$ is closed.
\pf
We first show that $A_{\e} + E$ is closed.
\newline \newline
Consider a convergent sequence of points $\{x_n\} \subset A_{\e} +E$. We must show that 
$$x:=\lim_{n \rightarrow \infty}x_n \in A_{\e} + E.$$
If 
$$x \in E$$
we are finished, so assume that this is not the case.
Now, for each $x_n$, either $x_n \in E$ 
or $x_n \in A_{\e}$. 
\newline \newline
In the first case $x_n = e_i \in E$ and there is an $n_0 \in \mathbb{N}$ such that 
$e_i \in A_n$ for each $n > n_0$. By taking $\{x_{n,j} \}_{j=1}^{\infty}$ such that $x_{n,j} = x_n$ for each $j$ then 
$x_{n,j} \in A_m$ for some $m \geq j$ for each $j$ and $\lim_{j \rightarrow \infty}x_{n,j}=x_n$. 
\newline \newline
In the second case 
$$x_n \in A_{\e} = \overline{\left(\bigcup_{m=1}^{\infty}A_m\right)} \sim \bigcup_{m=1}^{\infty}A_m.$$
Thus there exists a sequence $x_{n,j}$ such that $|x_{n,j} - x_n| < 1/j$ so that $\lim_{j \rightarrow \infty}x_{n,j} = x_n$ 
and $\{x_{n,j}\}_{j=1}^{\infty} \subset \bigcup_{m=1}^{\infty}A_m$.
Now assume that there is a finite number $q$ such that $\{x_{n,j}\}_j \subset \cup_{m=1}^qA_m$. Then since 
$\cup_{m=1}^qA_m$ is a finite union of closed lines it is closed so that $\lim_jx_{n,j} \in \cup_{m=1}^qA_m$ 
and thus $x_n \in \cup_{m=1}^q A_m$. However, since
$x \not\in E$ and $\cup_{m=1}^qE(A_m)$ is finite, $d(x,\cup_{m=1}^qE(A_m) > 0$. Thus in this case 
$x_n \in \cup_{m=1}^qA_m \sim E$. It follows then that we would have 
\begin{eqnarray}
x_n & \not\in & \overline{\left(\bigcup_{m=1}^{\infty}A_m\right)} \sim \bigcup_{m=1}^{\infty}A_m +E \nonumber \\
& =& A_{\e} \cup E. \nonumber
\end{eqnarray}
We can therefore take a subsequence and relabel to assume that $x_{n,j} \in A_m$ for some $m \geq j$ for 
each $j \in \mathbb{N}$. 
\newline \newline
We now take the sequence $\{x_m\}_{m=1}^{\infty}$ given by
$$x_m = x_{m,m},$$
and note that $\{x_m\} \subset \bigcup_{n=1}^{\infty}A_n$ so that 
$\lim_{m \rightarrow \infty} x_m \in \overline{\cup_{n=1}^{\infty}A_n}$. By the condition that $|x_{n,j} -x_n| < 1/j$, 
this diagonal selection gives us
\begin{eqnarray}
x & = & \lim_{m \rightarrow \infty} x_m \nonumber \\
& \in & \overline{\bigcup_{n=1}^{\infty}A_n}. \nonumber
\end{eqnarray}
Since, following from construction 2, 
for each $n \in \mathbb{N}$ and each $y \in A_n \sim E$ there is a radius $r>0$ such that 
$d(y,\cup_{m=n+1}^{\infty}A_m)>r$ it follows that for each $n \in \mathbb{N}$ $x \not\in A_n \sim E$.
Thus 
\begin{eqnarray}
x & \in & \overline{\left(\bigcup_{m=1}^{\infty}A_m\right)} \sim \bigcup_{m=1}^{\infty}(A_m \sim E) \nonumber \\
& = & \overline{\left(\bigcup_{m=1}^{\infty}A_m\right)} \sim \bigcup_{m=1}^{\infty}A_m \cup E \nonumber \\
& = & A_{\e} \cup E. \nonumber 
\end{eqnarray}
We therefore have that $A_{\e}$ is closed.
\newline \newline
Now since $\cal{B}$ is the countable union of open balls it is also open. Since $E \subset \cal{B}$ we can write 
$$\Ae = A_{\e} \sim \cB = \A \cup E \sim \cB$$
which is a closed set without an open set and thus is closed, proving the Lemma.
\end{proof}
\end{lem} 
\section{Properties of $\A$ and $\Ae$}
We now look at some direct properties of $\A$ and $\Ae$ that will be important to us later. Some of the properties, for 
example the dimension of $\A$ and $\Ae$ follow from more general Theorems that we shall use. However, since the direct 
proof is more instructive as to the properties of the sets and is not particularly longer, we present the direct proof 
here.
\begin{lem}\label{lem5} \thst
Let $\varepsilon > 0$ be such that $A_{\varepsilon}$ is well defined. Then 
for each $n \in \mathbb{N}$ and each $j \in \{1,2,...,2^n\}$,
the base length of a triangular cap $T_{n,j}$ in the construction of $A_{\varepsilon}$ has length
$$\Hm{1}(A_{n,j}) = {{(1+n16\varepsilon^2)^{1/2}}\over{2^n}},$$
and thus 
$$\Hm{1}(A_n) = (1+n16\e^2)^{1/2}$$
for each $n \in \mathbb{N}$.
\pf
Clearly $\Hm{1}(A_0) = 1$. 
$$\Hm{1}(A) \geq \liminf_{n \rightarrow \infty}\Hm{1}(A_n).$$
Then $\Hm{1}(A_1)$ is the sum of two hypotheses of triangles $(1/2)\Hm{1}(A_0)$ base length and $2\varepsilon$ height. 
that is 
\begin{eqnarray}
\Hm{1}(A_1) & = & 2\left(\left({{1}\over{2}}\right)^2 + (2\varepsilon )\right)^{1/2} \nonumber \\
& = & (1 + 16\varepsilon^2 )^{1/2} \nonumber
\end{eqnarray}
Having that it is true for $n=0,1$ I now claim that 
$$\Hm{1}(A_n) = (1+n16\varepsilon^2)^{1/2}.$$
Assuming it is true for $n$ we note that $\Hm{1}(A_{n+1})$ is the sum of $2^{n+1}$ hypotheses of triangles of
base length $\Hm{1}(A_n)/2^{n+1}$ and height $2^{2-n}\varepsilon$. That is 
\begin{eqnarray}
\Hm{1}(A_{n+1}) & = & 2^{n+1}\left(\left({{\Hm{1}(A_n)}\over{2^{n+1}}}\right)^2 + (2^{2-(n+1)}\varepsilon)^2\right)^{1/2} 
\nonumber \\
& = & ((\Hm{1}(A_n))^2 + 2^{2-2n + 2 + 2n}\varepsilon )^{1/2} \nonumber \\
& = & (1 + n16\varepsilon^2 + 16\varepsilon^2)^{1/2} \nonumber \\
& = & (1 + (n+1)16\varepsilon^2)^{1/2}, \nonumber 
\end{eqnarray}
proving the inductive claim. Then for each $n \in \mathbb{N}$ and $j \in \{1,2,...,2^n\}$ the base length of a 
triangular cap $T_{n,i}$ is one equal $2^n$-th part of the length of $A_n$. That is 
$$\Hm{1}(A_{n,j}) = {{(1+n16\varepsilon^2)^{1/2}}\over{2^n}}.$$
\end{proof}
\end{lem}\noindent
\begin{def1}\label{defpi}\thst
We denote the projection of a space onto a subset, $S$, whenever thye concept of projection makes sense for $S$ by 
$\pi_S$. An exception to this rule is the projection of $\R^2$ onto the $x$-axis identified with $\R$. This projection 
is denoted by $\pi_x$.
\end{def1}
\begin{thm}\label{thm1}\thst
dim$A_{\e} = $ dim$\Ae =$ $1$.
\pf
First note that 
$$\Hm{1}(A_{\e}) \geq \Hm{1}(\pi_x(A_{\e}))$$ 
and similarly
$$\Hm{1}(\Ae) \geq \Hm{1}(\pi_x(\Ae))$$
First, since $E$ is countable we can consider $x \in [0,1] \sim \pi_x(E) \not= \emptyset$. Since each $A_n$ can 
be considered as a connected path joining $(0,0)$ and $(0,1)$ there is an $x_n \in \pi_x^{-1}(x) \cap A_n$. 
Then we have $\{x_n\}_n$ a subsequence of $\cup_{n=1}^{\infty}A_n$. Since this sequence is in a bounded set 
($[0,1] \times [0,2\e]$) there is a convergent subsequence. Since for all $n \in \mathbb{N}$, 
$\pi_x(x_n) \not\in \pi_x(E)$, it follows that $x_0 = \lim x_n \not\in E$. Similarly to in the previous Lemma, this 
also implies that $x \not\in A_m$ for each $m \in \mathbb{N}$. Therefore 
$$x_0 \in \overline{\left(\cup_{n=1}^{\infty}A_n \right)}\sim \cup_{n=1}^{\infty A_n} \sim E = A_{\e}.$$
It follows that 
\begin{eqnarray}
\Hm{1}(\pi_x(A_{\e}) & \geq & \Hm{1}([0,1] \sim \pi_x(E)) \nonumber \\
& \geq & \Hm{1}([0,1]) - \Hm(\pi_x(E)) \nonumber \\
& = & 1 \nonumber \\
& > & 0. \nonumber
\end{eqnarray}
Now, we note that $\Ae = A_{\e} - \cB$ and that 
$$r_i = {{1}\over{4}}2^{-7}(1+7\cdot 16 \e^2)^{1/2} < 2^{-7}$$ 
(since we are in any case always taking $\e < 0.01$). It follows that 
\begin{eqnarray}
\Hm{1}(\pi_x (\Ae)) & = & \Hm{1}(\pi_x(A_{\e}\sim \pi_x(\cB)) \nonumber \\
& \geq & \Hm{1}(A_{\e}) - \sum_{i=1}^{\infty}r_i \nonumber \\
& > & 1 - 2^{-7} \sum_{i=1}^{\infty}4^{-i} \nonumber \\
& > & 1 - 2^{-7} \nonumber \\
& > & 0. \nonumber
\end{eqnarray}
It follows that 
\begin{equation}
\hbox{dim}A_{\e} \geq \hbox{dim}\Ae \geq \hbox{$1$}
\label{e:dimae1}
\end{equation}
\newline \newline
Now let $s > 0$ and $\delta>0$. Then for any given $\e>0$ 
there is an $n \in \mathbb{N}$ such that 
\begin{equation}
\delta \in (2^{2-n}\e,2^{3-n}\e]. 
\label{e:epsdel1}
\end{equation}
We note that the vertical height of 
the trianglular caps in the $n$-th stage of construction of $A_{\e}$ is $2^{1-n}$ so that $\delta \geq 2$ times the 
vertical height of the triangular caps in the $n$th construction stage. Since 
$$A_{\e} \subset \bigcup_{i=1}^{2^n}T_{n,i}$$
any cover of $\cup_{i=1}^{2^n}$ is also a cover of $A_{\e}$. By taking balls of radius $\delta$ with centers in 
$A_{n}$ we note that we can take these balls along an $A_{n,i}$ such that the overlaps ensure that 
$A_{n,i}^{\delta/\sqrt{2}}$ is covered. By taking such a cover of $A_{n,i}$ for each $i$ we have a cover consisting 
of balls of radius $\delta$, $\cB_{\delta} = \{B_{\delta}\} $ such that
$$\bigcup_{B_{\delta} \in \cal{B}_{\delta}} B_{\delta} \supset A_n^{\delta/\sqrt{2}} \supset A_n^{2^{1-n}}.$$
Since 
\begin{eqnarray}
A_n^{2^{1-n}} & \supset & \bigcup_{i=1}^{2^n}T_{n,i} \nonumber \\
& \supset & A_{\e} \nonumber
\end{eqnarray}
we also have that $\cB_{\delta}$ is a cover of $A_{\e}$.
Since with such a cover no more than $\delta/\sqrt{2}$ of the radius of a ball in $\cB_{\delta}$ will uniquely 
contribute to the cover of $A_n$, and since the inefficiencies of taking $A_{n,i}$'s that meet at non-uniform angles 
can not do any worse than forcing us to cover $A_n$ twice it follows that 
$$ \sum_{B_{\delta} \in \cB_{\delta}}\delta \leq 2\sqrt{2} \Hm{1}(A_n)$$
so that from Lemma 5 we have 
$$ \sum_{B_{\delta} \in \cB_{\delta}}\delta \leq 2\sqrt{2} (1+n16\e^2)^{1/2}.$$
Thus from ($\ref{e:epsdel1}$) we have 
$$ \sum_{B_{\delta} \in \cB_{\delta}}\delta^{1+s} \leq (2\e)^s 2\sqrt{2} {{(1+n16\e^2)^{1/2}}\over{2^{ns}}}.$$
Since $\cB_{\delta}$ is a cover of $A_{\e}$ this means 
$$\Hm{1+s}_{\delta}(A_{\e}) \leq (2\e)^s 2\sqrt{2} {{(1+n16\e^2)^{1/2}}\over{2^{ns}}}$$
so that we have
\begin{eqnarray}
\Hm{1+s}(A_{\e}) & = & \lim_{\delta \rightarrow 0} \Hm{1+s}_{\delta} \nonumber \\
& \leq & \lim_{n \rightarrow \infty} (2\e)^s 2\sqrt{2} {{(1+n16\e^2)^{1/2}}\over{2^{ns}}} \nonumber \\
& = & 0. \nonumber
\end{eqnarray}
Since this is true for all $s>0$ it follows that dim$A_{\e}$ $\leq 1$ and since $\Ae \subset A_{\e}$ 
that dim$\Ae$ $\leq 1$. Combining with ($\ref{e:dimae1}$) gives the result. 
\end{proof}
\end{thm}

\chapter{The limited Potency of Simple Examples and Weak Requirements for Locally Finite Measure}
\section{Limits on Approximately $j$-Dimensional Sets Entering and Exiting on the Same Side}
As we have already mentioned, several of the questions we are asking must be answered in the negative. To show this, 
clearly 
we need counter examples. Some of the counter examples, such as $\cal{N}$, $\Lambda_{\delta}$ and $\Lambda^2$ are 
relatively simple in that they are countable collections of nicely behaved functions whose relevant properties are 
clear. $\Gamma_{\e}$ is not so transparent as the sets already mentioned. It is, however, relatively clear that we need 
something a bit more complex to satisfy a $j$-dimensional approximation with a set that is not $j$-dimensional so 
as to provide a counter example for those properties not ensuring $j$-dimensionality. 
\newline \newline
$A_{\e}$ and $\Ae$, 
however are another matter, being "pseudo-fractal" sets (in the sense that every magnification of $A_{\e_1}$ looks 
like $A_{\e_2}$ for some $\e_2 <\e_1$ so that $A_{\e}$ is semi-selfsimilar.) that are in fact $j$-dimensional (where 
$j=1$ in this case). The obvious question is to ask if we could find a tricky way of putting nicely behaved functions 
together to get a different counter example to (iv) (2). (iv) is particularly important as we actually know that some 
singularity sets with 
a relationship to this property. We answer this question with an encouraging "no". This is encouraging as it means 
that to show that singularity sets have any sort of nice properties would then directly imply locally finite $\Hm{j}$-
measure. In fact, as mentioned previously, we can show that $\A$ and $\Ae$ are not countably $j$-rectifiable for the 
$j$ used in property (iv) which, since $\A$ and $\Ae$ are the only known counterexamples certainly supports the 
assertion that such sets must be poorly behaved.
\newline \newline
We find that any counter example must in fact be very poorly behaved in that for any point of locally infinite measure 
(where the essential part of a counter example is) cannot possibly have any part of the set (no matter how small) going 
through it that could be almost everywhere described by a Lipschitz function under some rotation and still satisfy property 
(iv). That is the set has to be a broken non-function at all critical points at all magnifications. 
\newline \newline
Conversely this means, to ensure a singularity set satisfying (iv) is locally $\Hm{j}$ finite we would expect only 
to need to show that no point on the singularity set has a neighbourhood in which the singularity set is purely 
unrectifiable.
\newline \newline
This section proves these assumptions. The key idea is that to have a function of infinite measure in a small neighbourhood 
means that at some point it has to be sharply folded on itself at all levels of magnification which will prevent the 
set from having property (iv). We make a couple of necessary definitions, then prove a Lemma proving an important special 
case which we use in the Theorem proving our claim.
\begin{def1}\label{def12} \thst
Let $u:\R \rightarrow \R$ be a function and let
$$\hbox{graph}u \cap B_{\rho}(y) \subset L_y^{\delta}$$
for some affine space $L_y$ and some $\delta \in (0,1/4)$. Then $u$ is said to {\bf enter and exit the same side of}
$B_{\rho}(y)$ {\bf with respect to } $L_y^{\delta}$ if 
$$\max\{|z-x|:y,x \in \hbox{graph}u \cap \partial B_{\rho}(y)\} < {{\pi \rho}\over{2}}.$$
\end{def1} \noindent
We note then that for a ball $B_{\rho}(y)$ and an affine space $L_y \ni y$ 
$$L_y^{\delta} \cap \partial B_{\rho}(y) = \Psi_1 \cup \Psi_2$$
for some arcs $\Psi_1$ and $\Psi_2$ in $\R^2$. We can therefore make the following definition.
\begin{def1}\label{def13}\thst
Suppose a function $u$ enters and exits $B_{\rho(y)}$ on the same side with respect to $L_y^{\delta}$. Then 
$$L_y^{\delta} \cap \partial B_{\rho}(y) = \Psi_1 \cup \Psi_2$$
for some arcs $\Psi_1$ and $\Psi_2$ in $\R^2$. Further graph$u \cap \Psi_i \not= \emptyset$ for exactly one 
$i =i(u) \in \{1,2\}$. We denote this $\Psi_{i(u)}$ by $\Psi_u$ and the other by $\Psi^u$.
\end{def1} \noindent
\begin{lem}\label{lem6}\thst
Suppose $u:\R \rightarrow \R$ is continuous and graph$u$ $\subset A \subset \R^2$.
Suppose that $A$ has property (iv) and that for some $y \in A$ and $\delta \in (0,1/4)$ $\rho_y$ is an 
appropriate radius at $y$ with respect to $\delta$. If $u$ enters and exits $B_{\rho_y}(y)$ on the same side. Then 
$$\max\{d(\Psi_u,y):y \in \hbox{graph}u \cap B_{\rho_y}(y)\} < 4\delta \rho_y.$$
\pf
We first show that 
$$\hbox{graph}u \cap B_{\rho_y}(y) \subset \hbox{graph}u(I_{u,y,\rho_y})$$
where 
$$I_{u,y,\rho_y} := [\inf\{\pi_x(\hbox{graph}u \cap \partial B_{\rho_y}(y))\},
\sup\{\pi_x(\hbox{graph}u \cap \partial B_{\rho_y}(y))\}].$$
Suppose that this were not to be the case, then there is a $z \in B_{\rho_y}(y) \subset \R^2$ with 
$z \in $ graph$u$ (and thus $u(\pi_x(z))=z$) and such that either 
$$\pi_x(z) > \sup\{\pi_x(\hbox{graph}u \cap \partial B_{\rho_y}(y))\} \hbox{or}$$
$$\pi_x(z)< \inf\{\pi_x(\hbox{graph}u \cap \partial B_{\rho_y}(y))\}.$$
without loss of generality we consider the case 
$$\pi_x(z) > \sup\{\pi_x(\hbox{graph}u \cap \partial B_{\rho_y}(y))\}$$
the other case follows similarly.
Since $u$ is a continuous function graph$u$ is connected and by the choice of $z$ 
$$\max\{\pi_x(B_{\rho_y}(y))\}>\max\{\pi_x(\hbox{graph}u \cap \partial B_{\rho_y}(y)\}$$
Thus the path 
$$P:=u([\max\{\pi_x(\hbox{graph}u \cap \partial B_{\rho_y}(y))\},\max\{\pi_x(B_{\rho_y}(y))\}+1])$$
intersects $B_{\rho_y}(y)$ only at its starting point on the boundary of $B_{\rho_y}(y)$. That is 
$$P \cap B_{\rho_y}(y) = u(\max\{\pi_x(\hbox{graph}u \cap \partial B_{\rho_y}(y)\})$$
(Otherwise $u(x) \cap \partial B_{\rho_y}(y) \not= \emptyset$ for some 
$x>\max\{\pi_x(\hbox{graph}u \cap \partial B_{\rho_y}(y))\}$ (in order for the connected path, $P$, to leave the ball) 
contradicting the choice of $\max\{\pi_x(\hbox{ graph}u \cap \partial B_{\rho_y}(y))\}$.)
\newline \newline
Thus 
$$\pi_x(z) \in \left[\max\{\pi_x(\hbox{ graph}u \cap \partial B_{\rho_y}(y))\},\max\{\pi_x(B_{\rho_y}(y)\}+1\right]$$
which implies 
$$u(\pi_x(z))\not\in B_{\rho_y}(y).$$
This contradiction means that $z \not\in \hbox{ graph}u$.
\newline \newline
For $z \in \hbox{ graph}u \cap B_{\rho_y}(y)$ Let 
$$z_{\partial} :=\pi_x^{-1}(\pi_x(z)) \cap \Psi_u$$
which will be a unique point. Now assume 
$$\max\{d(\Psi_u,z):z \in \hbox{graph}u \cap B_{\rho_y}(y)\} \geq 4\delta \rho_y$$
Then there is a $z \in \hbox{graph}u \cap B_{\rho_y}(y)$ such that 
\begin{eqnarray}
|\pi_y(z)-\pi_y(z_{\partial})| & > & d(z,z_{\partial}) \nonumber \\
& > & 4\delta \rho_y, \nonumber
\end{eqnarray}
since for all $a \in \Psi_u$, $|a-z_{\partial}| < 2\delta \rho_y$ and thus 
$|\pi_y(a) - \pi_y(z_{\partial})|<2 \delta \rho_y$.
\newline \newline
This implies 
$$\inf\{|\pi_y(z) - \pi_y(a)|:a \in \Psi_u\} > 2 \delta \rho_y.$$
w.l.o.g. assume that $\pi_y(z) > \sup\{\pi_y(a):a \in \Psi_u\}$.
\newline \newline
Then, as $u$ is continuous, there exist two connected paths $P_1$, $P_2$ such that 
\begin{list}{}{}
\item $\pi_x(P_1) \leq \pi_x(z)$,
\item $\pi_x(P_2) \geq \pi_x(z)$ and
\item $P_1$ and $P_2$ are connected to $\Psi_u$.
\end{list}
Thus
$$P_1 \cap \pi_y^{-1}(\pi_y(z) - 2\delta \rho_y) \not= \emptyset$$
and
$$P_2 \cap \pi_y^{-1}(\pi_y(z) - 2\delta \rho_y) \not= \emptyset.$$
Let
$$z_1 \in P_1 \cap \pi_y^{-1}(\pi_y(z) - 2\delta \rho_y)$$
and
$$ z_2 \in P_2 \cap \pi_y^{-1}(\pi_y(z) - 2\delta \rho_y).$$
Without loss of generality assume $|\pi_x(z_1)-\pi_x(z)| \leq |\pi_x(z_2) - \pi_x(z)|$
This choice implies that 
\begin{eqnarray}
|\pi_x(z_1) - \pi_x(z)| & \leq & 1/2 \sup\{|\pi_x(a_1)- \pi_x(a_2)|: a_1,a_2 \in \Psi_u\} \nonumber \\
& \leq & \delta \rho_y. \nonumber
\end{eqnarray}
Then notice  
\begin{eqnarray}
\rho_z & := & |z_2 - z_1| \nonumber \\
& \leq & \sup\{|\pi_x(a_1) - \pi_x(a_2)|:a_1,a_2 \in \Psi_u\} \nonumber \\
& = & 2\delta \rho_y \nonumber \\
& \leq & 1/2 \rho_y \nonumber
\end{eqnarray}
so we consider $B_{5\rho_z /4}(z_1)$. 
\newline \newline
Notice also that 
$|\pi_x(z) - \pi_x(z_1)| < |\pi_x(z_2) - \pi_x(z)|$ implies 
$$|\pi_x(z) - \pi_x(z_1)| \leq {{1}\over{2}}\rho_z.$$
Now call the subpath of $P_1 \subset \hbox{graph}u$ connecting $z_1$ to $z$ $P_{z_1}$. Note 
$$\pi_x(P_{z_1}) \subset [\pi_x(z_1),\pi_x(z)] \hbox{ and}$$
$$z \not\in B_{\rho_z}(z_1).$$
which implies 
$$P_{z_1} \cap \partial B_{\rho_z}(z_1) \not= \emptyset$$
and for all $$w \in |\pi_x(w)-\pi_x(z_1)|<{{1}\over{2}}\rho_z \hbox{ and}$$
$$d(w,z_1) = \rho_z$$
which implies 
$$|\pi_y(w) - \pi_y(z_1)| > {{\sqrt{3}}\over{4}}\rho_z.$$
However, for any choice of $L_{z_1,\rho_z}^{\delta \rho_z}$ we must have 
$$\sup\{|\pi_y(l)-\pi_y(z_1)|:l \in L_{z_1,\rho_z}^{\delta \rho_z}\} < {{9}\over{4}}\delta \rho_z.$$
Since $\delta < {{1}\over{16}}$ we note 
$${{\sqrt{3}}\over{4}}\rho_z > {{\rho_z}\over{4}} > {{9\rho_z}\over{64}} > {{9}\over{4}}\delta \rho_z.$$
Thus it is impossible to choose a $L_{z,\rho_z}$ such that 
$$A \cap B_{\rho_z}(z_1) \subset L_{z,\rho_z}^{\delta \rho_z}.$$
This would imply $A$ does not have property (iv). This contradiction proves the Lemma. 
\end{proof}
\end{lem} \noindent
\section{Set Constraints for Dually Approximately $j$-Dimensionality and Infinite Density}
We now prove the main theorem of this chapter by showing that we can reduce the problem to an application of the above 
lemma.
\begin{thm}\label{thm2} \thst
Suppose $A \subset \R^2$ and that there exists a $y \in A$ such that 
$$\Hm{1}(y \cap B_{\rho}(A)) = \infty \hbox{ for all } \rho > 0$$
and for some $\rho_1>0$,
\begin{list}{}{}
\item $y \in G_y^{-1}($graph$u) \cap B_{\rho_1}(y)$  and
\item $\overline{B_{\rho_y}(y)A \cap G_y^{-1}(\hbox{graph}u)} = G_y^{-1}($graph$u) \cap B_{\rho_y}(y) $
\end{list}
where $u$ is Lipschitz, $G_y \in G(1,2)$ and 
$G_y(\cdot):\R^2 \rightarrow \R^2$
is defined as the rotation such that $G_y(G_y) = \R$. \newline \newline
Then $A$ does not have property (iv) for $j=1$.
\pf
By the invariance of the relevant quantities under orthogonal transformations we can assume that $y = (0,0)$ and 
$G_y = \R$.
\newline \newline
Assume that $A$ does satisfy satisfy property (iv).
\newline \newline
Then for a given $\delta < 1/8$ there is a $\rho_y = \rho_y(y) \in (0,\rho_1)$ such that there exists an affine space 
$L_{y,\rho_y}$ such that 
$$A \cap B_{\rho_y}(y) \subset L_{y,\rho_y}^{\delta \rho_y}$$
and furthermore, for each $x \in A \cap B_{\rho_y}(y)$ and $\rho \in (0,\rho_y]$ there is an affine space $L_{x,\rho}$ 
such that 
$$A \cap B_{\rho}(x) \subset L_{x,\rho}^{\delta \rho}.$$
Noting that $y \in $graph$u$ and that clearly 
\begin{eqnarray}
d(y, \partial B_{\rho_y}(y)) & = & \rho_y \nonumber \\
& > & 4\delta \rho_y \nonumber 
\end{eqnarray}
it follows that 
$$\max\{d(\Psi_u,y):y \in \hbox{graph}u \cap B_{\rho_y}(y)\} < 4\delta \rho_y$$
and thus by Lemma $\ref{lem6}$ $u$ cannot enter and exit $B_{\rho_y}$ on the same side with respect to any affine space. 
\newline \newline
In particalar for each $w \in L_{y, \rho_y}$
$$\hbox{ graph}u \cap \pi_{L_{y,\rho_y}}^{-1}(w) \cap L_{y, \rho_y}^{\delta \rho_y} \not= \emptyset.$$
Also, if 
$$A \cap B_{\rho_y/2}(y) \subset \hbox{ graph}u$$
then
$$\Hm{1}(A \cap B_{\rho_y/2}(y)) \leq {{\rho_y}\over{2}} \cdot \omega_1 \cdot \hbox{ Lip}u < \infty,$$
a contradiction to our assumptions on the measure of balls around $y$.
\newline \newline
It follows that there exists an $x \in A \cap B_{\rho_y/2}(y)$ such that $x \not\in $graph$u$.
\newline \newline
Note that $\pi_{L_{y, \rho_y}}^{-1}(x) \cap $graph$u \not= \emptyset$ which implies 
$$d(x,\hbox{ graph}u) \leq 2\delta \rho_y < {{1}\over{2}}\rho_y.$$
Now select $z \in $graph$u$ such that
\begin{eqnarray}
d(z,x) & < & {{9}\over{8}}\inf\{d(w,x):w \in \hbox{ graph}u\} \nonumber \\
& =: & {{9}\over{8}}d \nonumber \\
& < & \rho_y. \nonumber
\end{eqnarray}
By the hypotheses there is an $z_1 \in $graph$u \cap A \cap B_{(1/16)d}(z)$. We now consider 
$B_{\rho_x}(z_1) \ni x$.
\newline \newline
Note that for any choice of $L_{z_1,\rho_x}$
$$L_{z_1,\rho_x}^{\delta \rho_x} \cap \partial B_{\rho_x}(z_1) = \Psi_1 \cup \Psi_2,$$
a union of two arcs as considered in Definition $\ref{def13}$ and that 
$$d(x,\partial B_{\rho_x}(z_1)) < {{1}\over{4}}d.$$
This implies that for some $i = i(x) \in \{1,2\}$
\begin{eqnarray}
\Psi_i & \subset & B_{(1/4)d + 2\delta \rho_x}(x) \nonumber \\
& = & B_{(1/4)d + 2\delta(9/8)d}(x). \nonumber 
\end{eqnarray}
Since $\delta$ was chosen such that $\delta < 1/8$
\begin{eqnarray}
{{1}\over{4}}d + 2\delta {{5}\over{4}}d & < & {{4}\over{16}} + {{5}\over{16}}d \nonumber \\
& < & {{15}\over{16}}d \nonumber
\end{eqnarray}
which implies 
$$\hbox{ graph}u \cap \Psi_{i(x)} = \emptyset.$$
This in turn implies that $u$ enters and exits $B_{\rho_x}(z_1)$ on the same side with respect to any 
affine space possibly allowing property (iv) to hold.
\newline \newline
Since $z_1 \in $ graph$u$
\begin{eqnarray}
\max\{d(w,\partial B_{\rho_x}(z_1)):w \in \hbox{ graph}u\}) & = & \rho_x \nonumber \\
& > & 4\delta \rho_x. \nonumber
\end{eqnarray}
This implies, by Lemma $\ref{lem6}$, 
that $A$ does not have property (iv). This contradiction completes the proof of the Theorem.
\end{proof}
\end{thm}\noindent
In order to more definitely relate what has previously been discussed to this result, I observe the following trivial 
corollaries.
\begin{cor}\label{cor3} \thst
Suppose $A \subset \R^2$ and that there exists a $y \in A$ such that 
$$\Hm{1}(y \cap B_{\rho}(A)) = \infty \hbox{ for all } \rho > 0$$
and for some $\rho_1>0$,
$$A \cap B_{\rho_1}(y) = G_y^{-1}\left(\bigcup_{n=1}^{Q}\hbox{ graph}u_n\right) \cap B_{\rho_1}(y)$$
for some $Q \in \mathbb{N} \cup \{\infty\}$
where $u_n$ is Lipschitz for each $n$, $G_y \in G(1,2)$ and 
$G_y(\cdot):\R^2 \rightarrow \R^2$
is defined as the rotation such that $G_y(G_y) = \R$. \newline \newline
Then $A$ does not have property (iv).
\pf
Since 
$$y \in A \cap B_{\rho_1}(y) = G_y^{-1}\left(\bigcup_{n=1}^{Q}\hbox{ graph}u_n\right) \cap B_{\rho_1}(y)$$
$y \in $ graph$u_{n_0}$ for some $1 \leq n_0 \leq Q$. With $u = u_{n_0}$ the conditions of Theorem 2 are then satisfied 
from which the conclusion follows. 
\end{proof}
\end{cor} \noindent
\begin{cor}\label{cor4}\thst 
$\cal{N}$, $\Lambda_{\delta}$ and $\Lambda^2$ are not counter examples to (iv) (2).
\pf
Let $\Xi = \cal{N}$ or $\Lambda_{\delta}$. Then since $\Xi$ is a countable 
union of Lipschitz graphs, any point of infinite density in $\Xi$ satisfies Theorem 2.
\newline \newline
For $\Lambda^2$ we note that the only point of density is $(0,0)$. Note that restricted to $[-1,1]$ the functions 
making up $\Lambda^2$, ($u_n = x^2/n$) are Lipschitz. Thus taking $\rho_1 = 1/2$ and $y = (0,0)$ in Theorem 2 the 
conditions of Theorem 2 are satisfied so that $\Lambda^2$ does not satisfy property (iv).
\end{proof}
\end{cor} \noindent
{\bf Remark} \newline
We note that in Lemma 6 and Theorem 2 we only used $\delta < 1/8$. Thus the full power of property (iv) has not been used. 
It is therefore possible and in fact likely that we could force any potential counter examples to (iv) (2) to be 
even stranger than what we have forced here. Even without using the $\delta$-fine property I believe that an 
improvement to Theorem 2 could be made in the form of the following conjecture.
\begin{conj}\label{conj1}\thst
Suppose $A \subset \R^2$ and that there exists a $y \in A$ such that 
$$\Hm{1}(y \cap B_{\rho}(A)) = \infty \hbox{ for all } \rho > 0$$
and for some $\rho_1>0$,
\begin{list}{}{}
\item $y \in G_y^{-1}($graph$u) \cap B_{\rho_1}(y)$  and
\item $\overline{A \cap G_y^{-1}(\hbox{graph}u)} = G_y^{-1}($graph$u) \subset A$
\end{list}
where $u \in C^0(\R;\R)$, $G_y \in G(1,2)$ and 
$G_y(\cdot):\R^2 \rightarrow \R^2$
is defined as the rotation such that $G_y(G_y) = \R$. \newline \newline
Then $A$ does not have property (iv).
\end{conj}\noindent
The idea being that although in this case the full infinite measure could all be produced from the one function, in 
the case where all the measure does come from the single function it must fold on itself sufficiently tightly and 
densely to either create a maximum or minimum somewhere we we could apply Lemma 6, or where essentially parallel lines 
would appear in which case choosing the correct size ball would mean that the approximating affine space would be 
essentially one of the lines and the intersection with the neighbouring line would then provide a contradiction to 
$A$ having property (iv). 
\newline \newline
More quantitatively, we note that there are several methods of attacking the proof and "almost getting there". One 
method, using Lemma 6, reduces the proof to the following.
\begin{conj}\label{conj2}\thst
Suppose $I_1, I_2$ are compact subintervals of $\R$ and 
$$u: I_1 \rightarrow I_2.$$
Suppose further that for all $x_1, x_2 \in I_1$ such that $u(x_1) = u(x_2)$
$$\sup\{|u(y) - u(x_1)|:y \in [x_1,x_2] \} < |x_1 - x_2|$$
Then, for any $\delta >0$ there exists a partition $P = \{p_1, ..., p_Q\}$ of $I_1$ with 
$$\max\{|p_i-p_{i-1}|:2\leq i \leq Q\} < \delta\}$$
and 
$$\sum_{i=2}^Q|u(p_i) - u(p_{i-1})| < C < \infty .$$
\end{conj} \noindent
Having discussed the non-simplicity of counter examples to (iv) (2), we reform what we have shown in how it is written to 
emphasise that a set thus need only be (iv) and posses at every point of infinite density a piece of graph 
to be sure that we have locally finite measure. 
This is an improvement on previous theory since such sets need not even be weak locally rectifiable. Clearly, 
we must first give a formal definition of these types of sets.
\begin{def1} \label{def14} \thst
Let $\mu$ be a measure on $\R^{n+k}$. Then $A \subset \R^{n+k}$ is said to posses a piece of Lipschitz graph at $x \in A$ 
if there exists an $r>0$, $G \in G(n,n+k)$ and a Lipschitz function 
$u:G \mapsto G^{\perp}$ such that 
$$x\in graph u$$
and
$$\Hm{n}((graph u \sim A) \cap B_r(x)) = 0.$$
\end{def1}\noindent
\begin{def1}\label{def15} \thst
A set $A$ is called weak locally countably $n$-rectifiable if for all $x \in A$ there exists $r>0$ such that 
$$A \cap B_r(x)$$
is countably $n$-rectifiable.
\end{def1}\noindent
It is clear that Definition $\ref{def14}$ is the same condition as that given in Theorem $\ref{thm2}$ so that the 
claim that this condition together with (iv) leads to locally finite measure follows from the same theorem.
\newline \newline
The claim that this is a lesser task to showing rectifiability follows from the existence of n-unrectifiable sets of 
$\Hm{n}$ finite measure for any $n$. Thus any set satisfying Definition $\ref{def14}$ in union with any $n$-unrectifiable 
set continues to satisfy the conditions of Definition $\ref{def14}$.

\chapter{Fitting the Counter Examples}
We mentioned in Chapter 2 that only questions with the answer "no" remain to be shown. 
In this section we show these results by appropriately fitting counter examples. For us this means showing 
firstly that the set in fact satisfies the definition that we claim it does and secondly that the set either 
has the wrong dimension (i.e. dimension greater than 1) or does not have locally finite 
$\Hm{1}$-measure depending on which property it is to which we wish to answer "no". 
As mentioned in the introduction, the higher dimensional cases will be discussed the following chapter. The 
reason the general dimension is not dealt with here is that they in any case reduce to the $1$-dimensional case as 
we shall see. \newline \newline
There is in fact, in terms of classifying the properties of our defintions, little that remains to be shown. 
What remains, however, is technical and non-trivial. \newline \newline
Fitting counter the counter example to (iv) (2) in particular shows that a non-rectifiable set (we show that 
$\A$ and $\Ae$ are non-recitifiable later) spiralling at all points and magnifications does not spiral too tightly 
around any given point. \newline \newline
The structure of the Chapter is that we show that $\Lambda_{\delta}$ satisfies 
(vi) which will answer (vi) (2) in the negative. We do the same with $\Lambda^2$ for (iii). 
$\Ae$ is then shown to satisfy 
(iv) (actually via first showing that $\A$ satisfies (iv)), from which (iv) (2) is answered in the negative, 
and as a corollary therefore 
(iii) (2) is also answered in the negative. 
Finally $\G$ is shown to satisfy (v), from 
which it follows that (v) (1) is answered with a no, and therefore as a corollary, the remaining questions: 
(v) (2), (ii) (1) and (ii) (2) are also answered with no.
\newline \newline
The proofs that the sets satisfy the definitions are mainly geometric and will actually mostly involve fitting sets in 
cones and then considering an appropriate neighbourhood of the center point. For this we need to develop notation 
to describe the cones we are using. As we will also find sets that should be covered by a cone meeting at a point, 
notation and theory also need to be developed for angles between sets. The appropriate definitions will be made as 
(or shortly before) they are used.
\begin{def1}\label{def16}\thst
Let $A$ be a $1$-dimensional affine subspace of $\R^2$, $\delta > 0$ and $x \in \R^n$, 
then $A$ is said to be a subset of the {\bf $\delta$-cone} at $x$, $C_{\delta}(x)$, if 
$$A\subset \left\{y = (y_1,y_2) \in \R^n: {{y_2}\over{y_1}} < \delta \right\} + x =:C_{\delta}(x).$$
More generally, if $L$ is a $1$-dimensional affine space in $\R^2$, $x \in A \cap L$ and $\phi$ 
is the orthogonal transformation such that 
$$\phi(L) = \R$$
and
$$\phi(x) = 0$$
then we say that $A$ is a subset of the {\bf $\delta$-cone around $L$ at $x$}, $C_{\delta,L}(x)$ if 
$$A \subset \phi^{-1}\left(\left\{y = (y_1,y_2) \in \R^n: {{y_2}\over{y_1}} < \delta \right\}\right) =:C_{\delta,L}(x).$$
\end{def1} \noindent
\section{Simple Counter Examples}
We now present the relevant classification results following from the simpler counter examples.
\begin{prop} \label{prop4}\thst
$\Lambda_{\delta}$ satisfies (vi), and further does not have weak locally finite $\Hm{1}$ measure 
so that the answer to (vi) (2) (weakly locally finite measure) is no.
\pf
There are two types of points to consider. If $x \not= (0,0)$, then if $x = (x_1, x_2)$
$$x \in \hbox{ graph}\left({{sgn(x_1)sgn(x_2) \delta x}\over{n}}\right) $$ 
for some $n \in \mathbb{N}$. 
Then for 
$$r_x = {{|x|\delta}\over{4(n+1)}},$$
\begin{eqnarray}
B_{r_x}(x) \cap \Lambda_{\delta} & \subset & \hbox{ graph}\left({{sgn(x_1)sgn(x_2) \delta x}\over{n}}\right) \nonumber \\
& \subset & G_{\delta /n,x}^{\delta r}, \nonumber
\end{eqnarray}
where $G_{\delta /n} \in G(1,2)$ is the affine space defined by graph$((sgn(x_1)sgn(x_2) \delta x) /n)$,
for each $r \in (0,r_x]$. 
Thus, by setting $L_x = G_{\delta /n,x}$, $x$ is an acceptable point with respect to (vi).
\newline \newline
If $x = (0,0)$, then by construction, we may choose $L_x = \R$ and note that 
$$G_{\delta /n,x} \subset C_{\delta}(x)$$
for each $n \in \mathbb{N}$, so that 
$$\Lambda_{\delta} \subset C_{\delta}(x).$$
It follows that 
$$\Lambda_{\delta} \subset \R^{\delta \rho} = L_x^{\delta \rho}$$
for each $\rho>0$. Thus choosing a $r_x>0$ at random we have 
$$\Lambda_{\delta} \subset L_x^{\delta r}$$ 
for each $r \in (0,r_x]$.
\newline \newline
It follows that $\Lambda_{\delta}$ satisfies (vi).
\newline \newline
Note, however, that due to the fact that there are countably infinitly many lines of length $2r$ going through any 
ball of radius $r$ around $(0,0)$, it follows that for all $r >0 $
$$\Hm{1}(\Lambda_{\delta} \cap B_r((0,0))) = \infty$$
so that $\Lambda_{\delta}$ is not weak locally $\Hm{1}$ finite. It follows that the answer to (vi) (2) 
is no.  
\end{proof}
\end{prop}\noindent
\begin{prop}\label{prop5}\thst 
$\Lambda^2$ satisfies (iii), and further does not have weak locally finite $\Hm{1}$ measure 
so that the answer to (iii) is no.
\pf
There are two types of points to consider. If $x \not= (0,0)$, then if $x = (x_1, x_2)$
$$x \in \hbox{ graph}\left({{sgn(x_1)sgn(x_2) \delta x^2}\over{n}}\right) $$ 
for some $n \in \mathbb{N}$. 
Then for 
$$r_x = {{|x|^2\delta}\over{4(n+1)}},$$
$$B_{r_x}(x) \cap \Lambda^2  \subset  \hbox{ graph}\left({{sgn(x_1)sgn(x_2) \delta x^2}\over{n}}\right)$$
Since also $x^2$ is differentiable there is a tangent line $L_{x}$ to $sgn(x_1)sgn(x_2)x^2/n$ at $x$ and a radius 
that can be chosen to be smaller than $r_x$, $r_{x_1} = r_{x_1}(\delta)>0$,  such that for all 
$$ y  \in \hbox{ graph}{{sgn(x_1)sgn(x_2)x^2}\over{n}} \cap B_{r_{x_1}}(x)$$
$$|\pi_{L_x^{\perp}}(y) - \pi_{L_x^{\perp}}(x)| < \delta |\pi_{L_x}(y) - \pi_{L_x}(x)|$$
so that 
$$B_{r}(x) \cap \Lambda^2 \subset L_{x}^{\delta r}$$
for each $r \in (0,r_{x_1}]$.
Thus $x$ is an acceptable point with respect to (vi).
\newline \newline
If $x = (0,0)$, then by construction, we may choose $L_x = \R$ and note that 
for $|x| < \delta$
\begin{eqnarray}
{{|x^2|}\over{n}} & = & {{|x||x|}\over{n}} \nonumber \\
& < & |x| \delta \nonumber
\end{eqnarray}
for each $n \in \mathbb{N}$. Thus it follows that for each $r \in (0,r_x=\delta]$
$$\Lambda^2 \cap B_{r}((0,0)) \subset L_x^{r\delta}.$$
It follows that $\Lambda^2$ satisfies (vi).
\newline \newline
Note, however, that due to the fact that there are countably infinitly many lines of length greater than 
or equal to $2r$ going through any 
ball of radius $r$ around $(0,0)$, it follows that for all $r >0 $
$$\Hm{1}(\Lambda^2 \cap B_r((0,0))) = \infty$$
so that $\Lambda^2$ is not weak locally $\Hm{1}$ finite. It follows that the answer to (vi) (2) 
is no. 
\end{proof}
\end{prop}\noindent
\section{Spiralling}
For $A_{\e}$ and $\Ae$ we show that the required measure properties hold first. That is that both of the 
sets are not weak locally $\Hm{1}$-finite. After that we then demonstrate 
that the set indeed satisfies (iv). Indeed, we have to work quite hard to get the results for 
$\Gamma_{\varepsilon}$ and $A_{\varepsilon}$. This arises from the fact, as has been mentioned and as will 
be shown in the next chapter, that $\G$ and $\A$ develope spirals in the set. 
In order to show the required properties we need to show that these spirals are not too tight. 
We now prove a technical lemma 
showing that we can find a "spiral free" view of our sets $\G$ and $\A$. We can then discuss the measure properties of 
$\A$ and $\Ae$.
\newline \newline
In order to discuss spiralling, we clearly need to discuss angles. For us, most essential will be 
the angle between two sets, particularly the angle between two trianglular caps. As simply saying the angle 
between two sets is unclear, we make a definition that will be sufficient for our needs.
\begin{def1}\label{def17} \thst
Let $A$ and $B$ be two sets with a single common point $z$ that can be divided by some $G \in G(1,2)$ in a 
sense that is explained below. 
Then the angle between the two sets $\psi^A_B$ is defined by 
$$\psi^A_B = \min \{\theta: C_{\theta}(z) \supset G(A \cup B) \hbox{ for some } G \in G(1,2) 
\hbox{ dividing } A \hbox{ and }B\}$$
where as usual $G(1,2)$ is the grassman manifold, $G(\cdot)$ denotes the rotation that takes $G \in G(1,2)$ to $\R_x$, 
and $G$ divides $A$ and $B$ if for all $X \in A$, $\pi_x(X) \leq 0$ and for all $Y \in B$, $\pi_x(Y) \geq 0$
\end{def1} \noindent
{\bf Remarks:} 
Clearly if $A_1 \subset A$, and $B_1 \subset B$ are such that $A_1 \cap B_1 = A \cap B = \{z\}$ then 
$\psi^{A_1}_{B_1} \leq \psi^A_B$. Note that the order is important due to the dividing of $A$ and $B$. The 
notation $\psi^A_B$ will always denote that $A$ is in the "left cone half" (i.e. $\pi_x(G(A)) \subset \R_x^{-}$) 
and $B$ is in the "right cone half" (i.e. $\pi_x(G(B)) \subset \R_x^{+}$) for the $G$ giving the 
minimum. 
We note that $\psi^{(\cdot)}_{(\cdot)}$ is subadditive in the sense that, if $A,B$ and $C$ are sets for which the 
definition makes sense for the pairings $\{A,B\}$ and $\{B,C\}$ with $z_1 = A \cap B$ and $z_2 \in B \cap C$, 
then 
$$\psi^A_{C - \{z_2-z_1\}} \leq \psi^A_B + \psi^B_C,$$
provided that such a value is less than $\pi/2$ (to ensure the dividing of the sets continues to make sense). Note  
that $\psi^{(\cdot)}_{(\cdot)}$ is translation and rotation invariant.
We note also particularly that in considering the angle between sets $A$ and $B$, 
if there is an affine space $L$ such that $A \cap L = \{z,z_a\}$ (i.e. contains the 
point common with $B$, $z$, and another point), then $\psi^L_B \leq \psi^A_B$ otherwise it would be impossible 
to contain $z_a$ and $B$ in a cone of angle $\psi^A_B$ around $z$.
\newline \newline
We also need to consider the angles that are actually intrinsic to the triangular caps.
\begin{def1}\label{def18}\thst
Let $n \in \mathbb{N} \cup \{0\}$ and $j \in \{1,2, ..., 2^n\}$, then we see from Constructions 1,2 and 3 that 
the triangular cap $T_{n,j}$ is an isosceles triangle. We denote the angles of $T_{n,j}$ as $\theta_{n,j}^A$ and 
$\pi - 2\theta_{n,j}^A$ 
where
$$\theta_{n,j}^A = \tan^{-1}\left({{2^{2-n}\varepsilon}\over{{{(1+n16\varepsilon^2)^{1/2}}\over{2^{n+1}}}}}\right)$$
and where the $\e$ is that associated with the construction of $\A$. Should the set $A$ be understood we will simply write 
$\theta_{n,j}$. Further, as in this chapter, should the $\T(n,j)$ be independent of $j$ for the understood set $A$; 
$\T(n,j)$ will be written $\T(n,\cdot)$. \newline \newline
Also, suppose that $L$ is an $1$-dimensional affine subspace (i.e. a line) of $\R^2$ of finite length 
(so that it has a middle point $l$), then we use $O_L$ to denote the orthogonal transformation such that 
$$O_L:L \rightarrow \R$$
and 
$$O_L(l) = (0,0).$$
\end{def1} \noindent
{\bf Remark:} At the present time the angles $\theta_{n,j}^A$ are independent of the index $j$. However, in Chapters $7$ 
and $8$ when we look at general forms of the construction of $\A$, the angles will be allowed to vary dependent on $n$ and 
$j$. For uniformity and simplicity later in the work, we introduce the symbol for the more general needs immediately.
\newline \newline
{\bf Note:} We note that from here on we take $\psi(0,\e) < \pi/32$. 
Thus we need $\e$ such that 
$$tan^{-1}\left({{8\varepsilon}\over{(1+16\varepsilon^2)^{1/2}}}\right) < {{\pi}\over{32}}$$
(coming from the definition of $\psi(n,\e)$.)
That is 
$${{8\varepsilon}\over{(1+16\varepsilon^2)^{1/2}}} < 0.09$$
so that taking 
$$0 < \varepsilon < {{1}\over{100}}$$
is sufficient. Since we in any case want to look at 
very small $\e$ and eventually will also be looking at $\e \rightarrow 0$, this presents us with no 
problems. We will therefore henceforth assume the $\e$ used to construct $\G$, $\A$, $\Ae$ and other similar sets 
is less than $0.01$. The reason for this assumption is that it is required for the spiralling Lemmas to work.
\begin{lem}\label{lem7}\thst
Suppose that $\A$, $\Ae$ and $\G$ are as defined in Constructions 1,2 and 3. Then 
\newline \newline
(1) \newline
should two neighbouring triangles, $T_{n,i}$ and $T_{n,i+1}$, be contained in another (necessarily earlier) triangular 
cap $T_{m,j(i)}$ ($m \leq n$) then
$$\psi^{T_n,i}_{T_{n,i+1}} \leq 2 \T{m,j(i)} \leq 2 \T{0,1}.$$
and \newline \newline
(2) \newline
the rectangle 
$$R_{n,i} = \pi_x \left( O_{A_{n,i}}\left(\cup_{j:|i-j|\leq 1}A_{n,j}\right) \right) \times 
[-2\Hm{1}(A_{n,i}),2\Hm{1}(A_{n,i})]$$
has the property
$$O_{A_{n,i}}^{-1}(R_{n,i}) \cap A \subset \bigcup_{j:|i-j|\leq 1}A_{n,j},$$
in the case of $\A$ and $\Ae$ and 
$$O_{A_{n,i}}^{-1}(R_{n,i}) \cap A \subset \bigcup_{j:|i-j|\leq 1}T_{n,j},$$
in the case of $\G$
\pf 
We write the proof for $\A$, from which the proofs for $\Ae$ and $\G$ follows. 
This is true for $\Ae$ since $\Ae \subset \A$ and it is true for $\G$ since we make all claims 
with respect to the triangular caps, and the second claim for $\A$ follows by noting that in $\T{n,j}$ 
only $A_{n,j}$ is in $A$ in any case. The only additional tool used is properties of $\T{n,j}$. However 
since the only property of $\T{n,j}$ from the construction of $\A$ that is used is that $\T{n,j} \leq \T{m,i}$ for 
$m \leq n$ and 
since $\T{n,j} \equiv \T{0,1}$ for all $n \in \mathbb{N}$, $j \in \{1,...,2^n\}$ 
in the construction of $\G$, all arguments involving $\T{\cdot,\cdot}$ also translate directly to $\G$.
\newline \newline
For (1), let $T_{n,i}$ and $T_{n,i+1}$ be two neighbouring triangular caps with common point $z$. Then, by the 
construction of $A_{\varepsilon}$, $z = z_{n_1+1,2i_1}$ is the vertex of a triangular cap 
$T_{n_1,i_1}$ for some $n_1 <n$ and some appropriate $i_1$. Further, since $z \in T_{m,j(i)}$ and 
$T_{n,i}, T_{n,i+1} \subset T_{m,j(i)}$ so that $z \not\in E(A_{m,j(i)})$ $m<n_1$ as otherwise the vertex $a_{n_1,i_1}$ 
cannot be in $T_{m,j(i)}$.
\newline \newline
Then by considering $G_{n_1,i_1} \in G(1,2)$ chosen such that $G_{n_1,i_1}||A_{n_1,i_1}$ 
we see that we can choose two "halves" (divided at $z_{n_1+1, 2i_1}$) of $G_{n_1,i_1}$, $G_{n_1,i_1}^-$ and 
$G_{n_1,i_1}^+$, such that 
$$\psi^{A_{n_1+1,2i_1}}_{G_{n_1,i_1^+}+z} \leq \T{{n_1},\cdot} \hbox{ and }
\psi_{A_{n_1+1,2i_1-1}}^{G_{n_1,i_1^-}+z} \leq \T{{n_1},\cdot}$$ 
so that, since in both cases in finding the minimum over cones, from which the definition of
$\psi^{A_{n_1+1,2i_1}}_{G_{n_1,i_1^+}+z}$ and 
$\psi_{A_{n_1+1,2i_1-1}}^{G_{n_1,i_1^-}+z}$ comes, we used the cone with respect to $G_{n_1,i-1}$, we have
$$\psi^{A_{n_1+1,2i_1}}_{A_{n_1+1,2i_1-1}} \leq \T{{n_1},\cdot}.$$
Since then $T_{n_1+1, 2i_1-1}$ and $T_{n_1+1, 2i_1}$ are constructed on the interior of $T_{n_1, i_1}$ with a base 
angle of $\T{{n_1}+1,\cdot}$, it follows similarly that 
$$\psi^{T_{n_1+1,2i_1}}_{G_{n_1,i_1^+}+z} \leq \T{n_1,\cdot}+\T{{n_1}+1,\cdot} \hbox{ and }
\psi_{T_{n_1+1,2i_1-1}}^{G_{n_1,i_1^-}+z} \leq \T{{n_1},\cdot} + \T{{n_1}+1,\cdot}$$
so that, since we have, as above, in both cases again made the statements about 
$\psi^{\cdot}_{\cdot}$ with respect to a cone around $G_{n_1,i_1}$
$$\psi^{T_{n_1+1,2i_1-1}}_{T_{n_1+1, 2i_1}} \leq \T{{n_1},\cdot} + \T{{n_1}+1,\cdot}.$$
Now, since $\T{n,\cdot} > \T{m,\cdot}$ for all $n<m$ it follows that $\T{{n_1},\cdot} \leq \T{m,\cdot} \leq \T{0,\cdot}$ 
and that 
$\T{{n_1}+1,\cdot} \leq \T{m,\cdot} \leq \p{0}$ so that 
$$\psi^{T_{n_1+1,2i_1-1}}_{T_{n_1+1, 2i_1}} \leq 2\p{m} \leq 2\p{0}.$$
Finally, we note that now, by construction (in that $A_{\varepsilon}$ is defined through intersection of the constructing 
levels) that $T_{n,i} \subset T_{n_1,i_1}$ and $T_{n,i+1} \subset T_{n_1,i_1+1}$ so that 
$$\psi^{T_{n,i}}_{T_{n,i+1}} \leq 2\p{m} \leq 2\p{0}.$$
This proves $(1)$. 
\newline \newline
For (2), note that since $\e < 1/100$, $\p{0} < \pi/32$.
\newline \newline
We first need to make a subclaim. 
\newline \newline
The claim is that if $T_{n,i}$ and $T_{n,j}$ are triangular caps with $2 \leq |i-j| \leq 3$ then 
$$\pi_x\left(O_{A_{n,i}}\left(\bigcup_{j:|i-j|<2}T_{n,j}\right) \right) \cap 
\pi_x(O_{A_{n,i}}(T_{n,j})-\{z_{n,i-2},z_{n,i+1}\}) = \emptyset$$
From this claim we will prove (2). As claimed above, we note that since 
$$\bigcup_{j:|i-j<2}A_{n,j} = A \cap \bigcup_{j:|i-j|<2}T_{n,j}$$ 
it is sufficient to prove that for any $T_{n,i}, T_{n,i+1}, T_{n,i+2}$ we have 
$$A \cap \pi_x\left(O_{A_{n,i+1}}\left(\bigcup_{j:|i+1-j|\leq 1}T_{n,j}\right)\right) \times 
[-2\Hm{1}(A_{n,\cdot}),2\Hm{1}(A_{n,\cdot})] \subset \bigcup_{j:|i+1-j|\leq 1}T_{n,j}.$$
We now consider our claim.
\newline \newline
We prove the case for $j-i > 0$, the other case following symmetrically. Note that we know from (1) that 
$$\psi^{T_{n,i}}_{T_{n,i+1}} \leq 2\p{0}$$
and that 
$$\psi^{T_{n,i+1}}_{T_{n,i+2}} \leq 2\p{0}$$
so that 
$$\psi^{T_{n,i}}_{T_{n,i+2}-(z_{n,i+1} - z_{n,i})} \leq 4\p{0}.$$
Indeed, since 
$$\psi^{T_{n,i+2}}_{T_{n,i+3}} \leq 2\p{0},$$
\begin{eqnarray}
\psi^{T_{n,i}}_{T_{n,i+2}-(z_{n,i+2} - z_{n,i})} & = & \psi^{T_{n,i}}_{T_{n,i+3}-(z_{n,i+2}-z_{n,i+1})-(z_{n,i+1}-z_{n,i})} 
\nonumber \\
& \leq & \psi^{T_{n,i}}_{T_{n,i+1}} + \psi^{T_{n,i+1}}_{T_{n,i+2}-(z_{n,i+2} - z_{n,i+1})} \nonumber \\
& \leq & \psi^{T_{n,i}}_{T_{n,i+1}} + \psi^{T_{n,i+1}}_{T_{n,i+2}} + \psi^{T_{n,i+2}}_{T_{n,i+3}} \nonumber \\
& \leq & 6\p{0}. \nonumber
\end{eqnarray}
It thus follows that $\psi^{A_{n,i}}_{T_{n,i+3} - (z_{n,i+2} - z_{n,i})} \leq 6\p{0}$. 
\newline \newline
Since $A_{n,i}$ is a line meeting the center of the cone 
$$C_{6\p{0}}(G(z_{n,i})) \supset G(A_{n,i} \cup (T_{n,i+3}-(z_{n,i+2}-z_{n,i})))$$ 
it follows that 
$$O_{A_{n,i}}(G^{-1}(C_{6\p{0}}(G(z_{n,i})))) \subset C_{12\p{0}}((0,\Hm{1}(A_{n,i})/2))$$
and thus that 
$$O_{A_{n,i}}(T_{n,i+3}-(z_{n,i+2}-z_{n,i})) \subset C^+_{12\p{0}}((0,\Hm{1}(A_{n,i})/2))$$
(where $C^+$ denotes the RHS of the cone), and therefore from translation invariance of the cone containing a set
$$O_{A_{n,i}}(T_{n,i+3}) \subset C^+_{12\p{0}, \R_x + z_{n,i+2}}(z_{n,i+2}).$$
This being the worse of the two possible $j$ cases
($j=i+1$ and $j = i+2$), an identical proceedure can be used to show that 
$$O_{A_{n,i}}(T_{n,i+2}) \subset C^+_{8\p{0}, \R_x + z_{n,i+1}}(z_{n,i+1}).$$
We note that 
$$8\p{0} < 12\p{0} < {{12\pi}\over{32}} < {{\pi}\over{2}}.$$
Thus
$$\pi_x(O_{A_{n,i}}(T_{n,i+2} \cup T_{n,i+3})) \subset [\pi_x(O_{A_{n,i}}(z_{n,i+1})),\infty)$$
and
$$\pi_x(O_{A_{n,i}}(T_{n,i+2} \cup T_{n,i+3})-\{z_{n,i+1},z_{n,i-2}\}) \subset (\pi_x(O_{A_{n,i}}(z_{n,i+1})),\infty).$$
We find that a similar argument to the above produces
$$O_{A_{n,i}}(T_{n,i+1}) \subset C^+_{4\p{0}, \R_x + z_{n,i}}(z_{n,i}).$$
so that since $4\p{0} < \pi/2-\p{0}$
\begin{eqnarray}
\max\{\pi_x(y):y \in O_{A_{n,i}}(T_{n,i+1})\} & = & \pi_x(O_{A_{n,i}}(z_{n,i+1})) \nonumber \\
& > & \pi_x(O_{A_{n,i}}(z_{n,i})) \nonumber \\
& = & \max\{\pi_x(y):y \in O_{A_{n,i}}(T_{n,i})\} \nonumber \\
& = & \pi_x(O_{A_{n,i}}(z_{n,i-1})) + \Hm{1}(A_{n,\cdot}) \nonumber \\
& \geq & \max\{\pi_x(y):y \in O_{A_{n,i}}(T_{n,i-1})\}. \nonumber 
\end{eqnarray}
Thus clearly
$$\pi_{x}\left(\bigcup_{j:|i-j|<2}O_{A_{n,i}}(T_{n,j})\right) \subset (-\infty,\pi_x(O_{A_{2,1}}(z_{2,3})],$$
so that
$$\pi_{x}\left(\bigcup_{j:|i-j|<2}O_{A_{n,i}}(T_{n,j})\right) \cap 
\pi_x(O_{A_{n,i}}(T_{n,i+2} \cup T_{n,i+3})-\{z_{n,i+1},z_{n,i-2}\}) = \emptyset$$
proving the claim.
\newline \newline
We now prove (2) by induction. We first note that for $A_0$ and $A_1$ it is obvious, as there are $1$ and $2$ 
triangular caps respectively, meaning that $A$ is clearly a subset of any "triple" (using " "  as it is actually 
impossible to choose a triple) of the form required. For $A_2$ there are four triangular caps, so that there is 
something to prove. However, we note that for any chosen $i$ every triangle is either in the "triple" around $i$ or has an 
index $j$ such that $2 \leq |i-j| \leq 3$. Since $A$ is a subset of the four triangles, the required result follows 
directly from the above proved claim.
\newline \newline
We now prove the inductive step. So we suppose that the inductive hypothesis (i.e. (2)) holds for all triples 
$\{T_{p,i-1}, T_{p,i},T_{p,i+1}\}$ for a given $p \in \mathbb{N}$ and show that it holds for an arbitrary triple 
$\{T_{p+1,i-1}, T_{p+1,i},T_{p+1,i+1}\}$. 
We set 
$$\mathcal{T} = \cup \{T_{p+1,i-1}, T_{p+1,i},T_{p+1,i+1}\}.$$
Note first that 
$$\bigcup_{j:|i-j|<2}T_{p+1,j} \subset \bigcup_{j:|i_1-j|<2}T_{p,j}$$
where $i_1 = (i/2)^{\sqcap}-1$ ($x^{\sqcap}$ is the smallest integer $q \geq x$), so that the triple is in fact 
a subset of a triple in the $p$th construction level. This triple in the $p$th construction level, by the induction 
hypothesis contains exactly $6$ trianglular caps in the ($p+1$)th construction level, namely 
$\{T_{p+1},j\}_{j=2i_1-3}^{2i_1+2}$ with $T_{p+1,i} \in \{T_{p+1,2i_1-1},T_{p+1, 2i_1}\}$. We also have by the inductive 
hypothesis that 
$$A \cap R_{p,i_1} \subset \bigcup_{j=2i_1-3}^{2i_1+2}T_{p+1,j}.$$
It follows that 
$$A \cap R_{p+1,i} \cap R_{p,i_1} \subset \bigcup_{j=2i_1-3}^{2i_1+2}T_{p+1,j}.$$
Now, since $i \in \{2i_1-1, 2i_1\}$ we see that for all $j \in \{2i_1-3, ..., 2i_1+2\}$, either $|i-j| < 2$ 
or $2 \leq |i-j| \leq 3$. From the above proven claim it follows that for each $j$ such that 
$2 \leq |i-j| \leq 3$, $(T_{p+1, j} \sim \mathcal{T}) \cap R_{p+1,i} = \emptyset$. Thus
$$A \cap R_{p+1,i} \cap R_{p,i_1} \subset \mathcal{T}.$$
The induction then follows in the case that $R_{p+1,i} \subset R_{p,i_1}$, as in this case 
$$A \cap R_{p+1,i} = A \cap R_{p+1,i} \cap R_{p,i_1} \subset \mathcal{T}.$$
We therefore prove that this is the case. It is clearly sufficient to show that 
$$O_{A_{p,i_1}}(R_{p+1,i}) \subset O_{A_{p,i_1}}(R_{p,i_1})$$
as in this case 
\begin{eqnarray}
R_{p+1,i} & = & O_{A_{p,i_1}} \circ O_{A_{p,i_1}}^{-1}(R_{p+1,i}) \nonumber \\
& \subset & O_{A_{p,i_1}} \circ O_{A_{p,i_1}}^{-1}(R_{p,i_1}) \nonumber \\
& = & R_{p,i_1}, \nonumber 
\end{eqnarray}
which is what we need.
\newline \newline
Without loss of generality we may assume that 
\begin{eqnarray}
O_{A_{p,i_1}}(T_{p+1,i}) & \subset & \triangle ((0,0),(-\Hm{1}(A_{p,j}/2,0),(0,\varepsilon \Hm{1}(A_{p,j}))) \nonumber \\
& \subset &\triangle ((0,0),(-\Hm{1}(A_{p,j}/2,0),(0,\Hm{1}(A_{p,j})/100)) \nonumber
\end{eqnarray}
where $\triangle(a,b,c)$ denotes the triangle in $\R^2$ with vertices $a,b$ and $c$. 
The other cases follow with symmetric arguments.
\newline \newline
We have
$$\pi_{O_{A_{p,i_1}}}\left(\bigcup_{j:|i-j|<2}O_{A_{p,i_1}}(T_p+1,j)\right) \subset $$
$$\left\{t\left(-\Hm{1}(A_{p,j},
-{{\Hm{1}(A_{p+1,j})}\over{100}}\right) + (1-t)\left({{\Hm{1}(A_{p,j})}\over{2}},{{2\Hm{1}(A_{p+1,j})}\over{100}}\right) 
:t \in [0,1] \right\};$$
so that 
$$O_{A_{p,i_1}}\left(\bigcup_{j:|i-j|<2}O_{A_{p,i_1}}(T_p+1,j)\right) \subset \{x=y+z\}$$
where 
$$ y \in \left\{t\left(-\Hm{1}(A_{p,j},
-{{\Hm{1}(A_{p+1,j})}\over{100}}\right) + (1-t)\left({{\Hm{1}(A_{p,j})}\over{2}},{{2\Hm{1}(A_{p+1,j})}\over{100}}\right) 
:t \in [0,1] \right\}$$
and 
$$z \in \left\{2s\Hm{1}(A_{p+1,j})\left({{-4}\over{100}},2\right):s \in [-1,1]\right\}.$$
That is $O_{A_{p,i_1}}(R_{p,i})$ is a subset of the quadrilateral with vertices 
$$V_1:=(-1.54\Hm{1}(A_{p,j}),2\Hm{1}(A_{p+1,j}))$$
$$V_2:=(0.96\Hm{1}(A_{p,j}),2.04\Hm{1}(A_{p+1,j}))$$
$$V_3:=(1.04\Hm{1}(A_{p,j}),-2\Hm{1}(A_{p+1,j}))$$
and
$$V_4:=(-1.46\Hm{1}(A_{p,j}),-2.04\Hm{1}(A_{p+1,j})).$$
Noting then that, due to the fact that $\p{0}<\pi/32$ and the general fact that $\psi^{T_{p,j}}_{T_{p,j+1}} < 2\p{0}$
(from (1)) we get 
\begin{eqnarray}
\Hm{1}(\pi_x(O_{A_{p,i_1}}(T_{p,j}))) & > & \cos\left({{\pi}\over{8}}\right)\Hm{1}(A_{p,j}) \nonumber \\
& > & 0.9\Hm{1}(A_{p,j}) \nonumber 
\end{eqnarray}
for all $j$ such that $|j-i_1| < 2$, and since 
\begin{eqnarray}
\Hm{1}(A_{p+1,j}) &=& {{1}\over{2}}(1 + 16\varepsilon^2)^{1/2}\Hm{1}(A_{p,k}) \nonumber \\
& < & 0.6\Hm{1}(A_{p,k}) \nonumber 
\end{eqnarray}
we have 
\begin{eqnarray}
R_{p,i_1} & = & O_{A_{p,i_1}}(R_{p,i_1}) \nonumber \\
& \supset & [-1.9\Hm{1}(A_{p,j}),1.9\Hm{1}(A_{p,j})] \times [-2\Hm{1}(A_{p,j}),2\Hm{1}(A_{p,j})] \nonumber \\
& \supset & [-3\Hm{1}(A_{p+1,j}),3\Hm{1}(A_{p+1,j})] \times [-3\Hm{1}(A_{p+1,j}),3\Hm{1}(A_{p+1,j})]. \nonumber
\end{eqnarray}
Since clearly
$$V_1,V_2,V_3,V_4 \in [-3\Hm{1}(A_{p+1,j}),3\Hm{1}(A_{p+1,j})] \times [-3\Hm{1}(A_{p+1,j}),3\Hm{1}(A_{p+1,j})]$$
it follows that 
$$O_{A_{p,i_1}}(R_{p+1,i}) \subset O_{A_{p,i_1}}(R_{p,i_1})$$ 
and thus that
$$R_{p+1,i} \subset R_{p,i_1}$$
completing the proof of (2). 
\end{proof}
\end{lem} \noindent
\section{Measure Properties of $\A$ and $\Ae$}
We now show that $\A$ and $\Ae$ are not weak locally $\Hm{1}$-finite. We start with the simpler: $\A$.
\begin{lem}\label{lem8}\thst
Let $\varepsilon > 0$ be such that $A_{\varepsilon}$ is well defined. Then 
$A_{\varepsilon}$ is not weak locally $\Hm{1}$ finite.
\pf
We note that for each $n_0 \in \mathbb{N}$, since $A_{\varepsilon}$ 
makes the lines in $A_{n_0}$ less straight, the refinement 
to $A_{\varepsilon}$ increases the measure of $A_{n_0}$. That is
$$\Hm{1}\left(\bigcap_{i=1}^{n_0}\bigcup_{n=i}^{\infty}A_n\right) \geq \Hm{1}(A_{n_0}).$$
Also, for arbitrary $n \in \mathbb{N}$, from Lemma $\ref{lem5}$ we have 
$$\Hm{1}(A_n) = (1+n16\varepsilon^2)^{1/2}.$$
So that 
\begin{eqnarray}
\Hm{1}(A_{\varepsilon}) & \geq & \liminf_{n \rightarrow \infty}(1 + n16\varepsilon^2)^{1/2} \nonumber \\
& = & + \infty \nonumber
\end{eqnarray}
Since for each $n\in \mathbb{N}$, 
$A_{\varepsilon} \subset \cap_{n=0}^{\infty}\cup_{i=1}^{2^n}T_{n,i}$ and $\lim_{n \rightarrow \infty}$ base length 
$T_{n,i}=0$ we know that for each $x \in A$, and for each $R>0$ there exists $n$ and $i$ such that 
$x \in T_{n,i} \subset B_{R}(x)$. Thus also $A_{n,i} \subset B_{R}(x)$.
\newline \newline
We now note that by the construction of $A_{\varepsilon}$ we actually have that the further construction of 
$A_{\varepsilon}$ on $A_{n,i}$ is 
the same as that for $A_{\varepsilon}$ except that we start with a base length $\Hm{1}(A_{n,i})$ of 
$\Hm{1}(A_n)/2^n$ instead of $1$. That is 
$A_{\varepsilon} \cap T_{n,i}$ is a version of $A_{\varepsilon}$ scaled by a factor of $\Hm{1}(A_{n,j})$. Thus 
\begin{eqnarray}
\Hm{1}(A_{\varepsilon} \cap B_{R}(x)) & \geq & \Hm{1}(T_{n,i} \cap B_{R}(x)) \nonumber \\
& \geq & {{\Hm{1}(A_n))}\over{2^n}}\Hm{1}(A_{\varepsilon}) \nonumber \\
& = & + \infty. \nonumber
\end{eqnarray}
Since this is true for each $x \in A_{\varepsilon}$ and each $R>0$ the conclusion follows. 
\end{proof}
\end{lem}\noindent
This result also leads to the following interesting result. Not only is it interesting in itself, showing that the set 
$\A$ has infinite density in its own dimension everywhere in the set. It is also useful in showing the 
nonrectifiability of $\A$ later on.
\begin{cor}\label{cor5}\thst
For each $y \in \overline{\A}$
$$\Theta^1(\Hm{1},\A,y) = \infty.$$
\pf
Let $y \in \overline{\A}$ and $\rho>0$. 
\newline \newline
Then there is a $y_1 \in \A \cap B_{\rho /2}(y)$ such that 
$B_{\rho /2}(y_1) \subset B_{\rho}(y)$. 
\newline \newline
Since $y_1 \in \A$ for each $n \in \mathbb{N}$ there is a triangular cap $T_{n,i(n,y)} \ni y$. Also, there is an 
$n_0 \in \mathbb{N}$ such that $\Hm{1}(A_{n,\cdot}) < \rho /4$ for each $n > n_0$ so that 
$T_{n,i(n,y)} \subset B_{\rho /2}(y_1)$ for each $n > n_0$. 
\newline \newline
Now From the symmetry of construction we see that 
$T_{n_0+1,i(n_0+1,y)}$ is a $\Hm{1}(A_{n_0+1,\cdot})$ scale copy of $A_{2^{-n_0}\e}$. However, from Lemma $\ref{lem8}$ 
we know \newline $\Hm{1}(A_{2^{-n_0\e}})~=~\infty$, thus 
\begin{eqnarray}
\Hm{1}(\A \cap B_{\rho}(y)) & \geq & \Hm{1}( \A \cap B_{\rho /2}(y_1)) \nonumber \\
& \geq & \Hm{1}(\A \cap T_{n_0+1,i(n_0+1,y)}) \nonumber \\
& = & \Hm{1}(A_{n_0+1,\cdot})\cdot \Hm{1}(A_{2^{-n_0\e}}) \nonumber \\
& = & \infty. \nonumber
\end{eqnarray}
It follows that 
\begin{eqnarray}
\Theta(\A,y) & = & \lim_{\rho \rightarrow \infty}{{\Hm{1}(B_{\rho}(y) \cap \A)}\over{\omega_1 \rho}} \nonumber \\
& = & \infty. \nonumber
\end{eqnarray}
\end{proof}
\end{cor}\noindent
Although having an infinitely dense point is not that uncommon, and infact having a set of $\Hm{1}$ positive 
measure of points of $\Hm{1}$ infinite density is not uncommon, that $\A$ is a set of positive $\Hm{1}$ measure 
that has infinite $\Hm{1}$ density at all points of its closure is less common, which makes $\A$ a set of 
peculiar interest in its own right without association to the properties that we are currently discussing. 
\newline \newline
Although it is possible that$\Ae$ has this same peculiar property, it is not easy to prove, and in fact we don't. 
We settle for finding one such point, however by removing small open balls around such points gauranteed by the proof 
that follows, we know that there must be at least countably many points in $\Ae$ of infinite density.
\newline \newline
Although the proof that $\Ae$ is not weakly locally $\Hm{1}$ finite is more involved than that for 
$\A$ it is similar. We find approximating sets (subsets of $A_n$) that we can take a limiting infimum of to 
bound the measure of $\Ae$ from below. We then show that this limiting infimum is infinite. The proof 
that there is a point of infinite density is then an indirect proof using a covering argument.
\begin{lem}\label{lem9}\thst
$\Hm{1}(\Ae) = \infty$.
\pf
Let $\delta >0$ and $\cB_{\delta}$ be a cover of $\Ae$ of balls of radius smaller than or equal to $\delta$. Then 
as $\Ae$ is compact we can find a finite subcover 
$$\cB_f = \{B_{w_i}(x_i)\}_{i=1}^Q$$
of balls of radius $w_i <\delta$. Further, since it is a finite collection we can define
$$w_m:=\{w_1,...,w_Q\}>0.$$
By appropriately selecting $8$ balls around each $B_{w_i}(x_i)$ of radius $w_i$ we get a new finite collection of 
$9Q$ balls with $$w_m \leq w_i \leq \delta$$ that we relabel
$$\cB_f = \{B_{w_i}(x_i)\}_{i=1}^{9Q}$$
such that 
\begin{equation}
\Ae^{w_m} \subset \bigcup_{i=1}^{9Q}B_{w_i}(x_i)
\label{e:l91}
\end{equation}
and
\begin{equation}
\Hm{1}_{\delta}(\Ae) = \inf\left\{{{1}\over{9}}\sum_{i=1}^{9Q}w_i \right\}
\label{e:l92}
\end{equation}
where the infimum is taken over all $\delta$-covers of $\Ae$.
\newline \newline
Now, suppose that $\gamma >0$ and that there is an $n_0 \in \mathbb{N}$ such that 
$$(A_n - \cB) \not\subset (\Ae)^{\gamma}$$
for all $n \geq n_0$.
\newline \newline
Then for each $\gamma > 0$ there is a sequence $\{x_n\}$ with $x_n \in A_n \sim \cB$ ($m \geq n$) such that 
$x_m \not\in \Ae^{\gamma}$. 
\newline \newline
Then, as $\{x_n\}$ is infinite in $[0,1]\times [0,2\e]$ there exists a convergent subsequence $\{y_n\}$ where 
for all $n_0 \in \mathbb{N}$ there is a $y_n$ such that 
\begin{equation}
y_n \in A_m \sim \cB \hbox{ for some } m \geq n_0.
\label{e:l93}
\end{equation}
We note that by construction 
$$\bigcup_{m=n+2}^{\infty}A_m \subset \bigcup_{i=1}^{2^{n+1}}T_{n+1,i}$$
and the boundaries (of the $T_{n+1,i}$'s) closest to $A_n$ are then the $(A_{n+2,i})$. Since then the angle 
between an $A_{n,\cdot}$ and an $A_{n+2,\cdot}$ that meet is $\p{n} - \p{n+1}$; that for all 
$B_{r_i}(x_i) \in \cB$ such that $A_n \cap B_{r_i}(x_i) \not= \emptyset$ we have
$$i \leq 2 + \sum_{j=0}^n2^j (<\infty)$$ 
and there exists a 
$$\min\left\{r_i:1\leq i \leq 2 + \sum_{j=0}^n2^j\right\}.$$
Since also, by Lemma $\ref{lem7}$ (2) the closest $A_{n+2,j}$ to an $A_{n,i}-\cup\{B_r(x):x \in E(A_{n,i})$ must be an 
$A_{n+2,j} \subset A_{n,i}$ we then have 
\begin{eqnarray}
d_n & := & d\left(A_n - \cB, \bigcup_{i=n+1}^{\infty}A_n \right) \nonumber \\
& \geq & d(A_n - \cB,A_{n+2} - \cB) \nonumber \\
& \geq & sin^{-1}(\p{n+1} - \p{n+2})\cdot\inf\left\{r_i:i \leq 2 + \sum_{j=0}^{n+2}2^j\right\} \nonumber \\
& > & 0. \label{e:l94}
\end{eqnarray}
Thus from ($\ref{e:l93}$) and ($\ref{e:l94}$) we have 
$$\lim_{n \rightarrow \infty}y_n \in \overline{\left(\bigcup_{n=1}^{\infty}(A_n - \cB)\right)}$$ 
and 
$$\lim_{n \rightarrow \infty}y_n \not\in A_m$$ for all $m \in \mathbb{N}$ which means 
$$\lim_{n \rightarrow \infty}y_n \in \Ae$$
but for all $n \in \mathbb{N}$, $y_n \not\in \Ae^{\gamma}$ so that 
$d(\Ae, y_n) > \gamma$ for all $n \in \mathbb{N}$ which implies 
$$d(\lim_{n \rightarrow \infty}y_n, \Ae) \geq \gamma.$$
This contradiction means that for all $\gamma > 0$ and all $n_o \in \mathbb{N}$ there is an $n > n_0$ such that 
$$(A_n - \cB) \subset \Ae^{\gamma}.$$
Thus for each $n_0 \in \mathbb{N}$ there is an $n . n_0$ such that 
$$(A_n - \cB) \subset \Ae^{w_m}$$
so that by ($\ref{e:l91}$) $\cB_f$ is a cover of balls of radius smaller than or equal to $\delta$ for $A_n - \cB$ and 
therefore 
$$\Hm{1}_{\delta}(A_n - \cB) \leq \sum_{i=1}^{9Q}w_i.$$
Since this is true for some $n \geq n_0$ for any $n_0 \in \mathbb{N}$ it follows that 
$$\liminf_{n \rightarrow \infty} \Hm{1}_{\delta}(A_n -\cB ) \leq \sum_{i=1}^{9Q}w_i.$$
Since this is true for any cover of $\Ae$ of balls with radius bounded above by $\delta$ it follows that 
$$\liminf_{n \rightarrow \infty}\Hm{1}_{\delta}(A_n - \cB) \leq \inf\left\{\sum_{i=1}^{9Q}w_i\right\}$$
where the infimum is taken over all $\delta$-covers of $\Ae$. Thus by ($\ref{e:l92}$)
$$\liminf_{n \rightarrow \infty}\Hm{1}_{\delta}(A_n - \cB) \leq 9\Hm{1}_{\delta}(\Ae).$$
Define 
$$L_n:=A_n - \cB.$$
We will show that 
$$\lim_{n \rightarrow \infty}\Hm{1}(L_n) = \infty.$$
Before we do this, however, we show how such a fact can be used to complete the proof.
\newline \newline
Assuming $\Hm{1}(L_n) = \infty$ we then know that for all $M \in \R$ there is an $n_0 \in \mathbb{N}$ such that 
$$\Hm{1}(L_n) > 2M$$ for all $n \geq n_0$. Let $n_0 \in \mathbb{N}$ be such a number and let 
$$L_{M,m} :=\{n \in \mathbb{N}: n \geq n_0 \hbox{ and } \Hm{1}_{1/m}(L_n) < M\}.$$
Since $\Hm{1}_{\delta}(A)$ is increasing as $\delta$ decreases for any $A \subset \R^2$, $L_{M,m_1} \subset L_{M,m_2}$ 
whenever $m_1 > m_2$.
\newline \newline
Then suppose there is no $m \in \mathbb{N}$ such that $L_{M,m} = \emptyset$, then 
$$\bigcap_{m=1}^{\infty}L_{M,m} \not= \emptyset$$
so that there is an $n \geq n_0$ such that 
$\Hm{1}_{1/m}(L_n) < M$ for all $m \in \mathbb{N}$. Thus 
\begin{eqnarray}
\lim_{\delta \rightarrow 0}\Hm{1}_{\delta}(L_n) & = & \lim_{m \rightarrow \infty}\Hm{1}_{1/m}(L_n) \nonumber \\
& < & M \nonumber 
\end{eqnarray}
contradicting $\Hm{1}(L_n) \geq 2M$. 
\newline \newline
It follows that there is a $\delta(M)>0$ such that $\Hm{1}_{\delta(M)}(L_n) \geq M$ for all $n$ such that 
$\Hm{1}(L_n) \geq 2M$. Thus 
$$\liminf_{n \rightarrow \infty}\Hm{1}_{\delta(M)}(L_n) \geq M.$$
Thus for all $M \in \R$ there is a $\delta(M) > 0$ such that for all $\delta < \delta(M) > 0$ 
\begin{eqnarray}
\Hm{1}_{\delta}(\Ae) & \geq & \Hm{1}_{\delta(M)}(\Ae) \nonumber \\
& \geq & {{1}\over{9}} \liminf_{n \rightarrow \infty}\Hm{1}_{\delta(M)}(L_n) \nonumber \\
& \geq & {{M}\over{9}}. \nonumber 
\end{eqnarray}
Thus 
\begin{eqnarray}
\Hm{1}(\Ae) & = & \lim_{\delta \rightarrow 0}\Hm{1}_{\delta}(\Ae) \nonumber \\
& = & \infty. \nonumber 
\end{eqnarray}
To complete the proof, we therefore now need to show that 
$$\lim_{n \rightarrow \infty}\Hm{1}(L_n) = \infty.$$
We consider first $L_7$. Note 
\begin{eqnarray}
\max\{r_i:B_{r_i}(x_i)(A_7) \not= \emptyset\} & = & r_1 \nonumber \\
& = & {{1}\over{4}}2^{-7}(1+7\cdot 16\e^2)^{1/2} \nonumber \\
& < & {{1}\over{2}}\Hm{1}(A_{7,\cdot}) \nonumber 
\end{eqnarray}
and thus from Lemma $\ref{lem7}$ (2) we know that this implies
\begin{eqnarray}
A_7 \cap B_{r_i}(x_i) & \subset & \cup\{T_{7,j}:x_1 \in E(T_{7,j})\} \nonumber \\ 
& = & T_{7,1}. \nonumber 
\end{eqnarray}
We remove this triangular cap from the measure that we count toward $L_7$ and note that 
\begin{eqnarray}
\sum_{B_{r_i}(x_i):B_{r_i}(x_i) \cap A_7 \not= \emptyset}r_i & \leq & \sum_{i=2}^{\infty}r_i \nonumber \\
& \leq & \sum_{i=1}^{\infty}r_i \nonumber \\
& < & \sum_{i=1}^{\infty}4^{-i}2^{-7}(1+7\cdot 16\e^2)^{1/2} \nonumber \\
& < & 2^{-6}. \nonumber
\end{eqnarray}
Each ball $B_{r_1}(x_i)$ such that $B_{r_i}(x_i) \cap A_7 \not= \emptyset$ (which by construction are those $B_{r_i}(x_i)$ 
such that $x_i = a_{n_i}$ for $n \leq 7$), again by Lemma $\ref{lem7}$ (2) meets only the two adjacent $A_{7,\cdot}$ 
so that by letting 
$$D_n:=\{\hbox{triangular caps discarded from the estimation of }L_n\}$$
we have $D_7 = \{T_{7,1}\}$ and thus 
$$\Hm{1}\left( \left(A_7 \sim \bigcup_{T_{7,i} \in D_7}T_{7,i}\right) \cap \cB\right) \leq 2^{-6}$$
and since each triangle $T_{7,i}$ gives the same value from 
$$\Hm{1}(A_7 \cap T_{7,i})$$
we have with 
$$N_n := \hbox{card}(D_n)$$
$$\Hm{1}\left(A_7 \sim \bigcup_{T_{7,i}\in D_7}T_{7,i}\right) = {{2^7 - N_7}\over{2^7}}\Hm{1}(A_7)$$
and thus 
$$\Hm{1}(L_7) \geq {{2^7 - N_7}\over{2^7}}\Hm{1}(A_7) - 2^{-6}.$$
We note in particular that $N_7 = 1 < 2$ and we make the following inductive claims.
\newline \newline
For each $n \geq 7$, by removing triangular caps intersecting 
$$\cup\{B_{r_i}(x_i):i \leq n-6\}$$ we have $N_n \leq 2^{n-5}-2$ so that 
\begin{eqnarray}
M_n & := & \{\hbox{ triangular caps remaining at the }n\hbox{th stage}\} \nonumber \\
& \geq & 2^{-n}-(2^{n-5}-2).\nonumber 
\end{eqnarray}
Further, we have that 
$$\sup\{r_1:r > n-6\} < {{1}\over{2}}\Hm{1}(A_{n,\cdot})$$
and that 
$$\Hm{1}\left(\bigcup_{i>n-6}B_{r_i}(x_i) \cap A_n\right) \leq \sum_{i=1}^{\infty} r_i < 2^{-6} $$
so that 
\begin{eqnarray}
\Hm{1}(L_n)&>&\Hm{1}\left(A_n\sim\bigcup_{T_{n,i}\in D_n}T_{n,i}\sim\left(A_n \cap \bigcup_{i>n-6}B_{r_i}(x_i)\right)\right)
\nonumber \\
& \geq & \Hm{1}\left(A_n\sim\bigcup_{T_{n,i}\in D_n}T_{n,i}\right) - 2^{-6} \nonumber \\
& \geq & {{2^n - (2^{n-5}-2)}\over{2^n}}\Hm{1}(A_n) - 2^{-6}.\nonumber
\end{eqnarray}
We know that these conditions hold for $n = 7$. Now we assume that they hold for some $n \geq 7$ and show the inductive step 
to show that they hold for $n+1$.
\newline \newline
First, we know that $N_n \leq 2^{n-5}-2$. This is the number of triangular caps that we have removed due to the intersection 
with balls $\{B_{r_i}(x_i)\}_{i=1}^{n-6}$. Thus there exist no more than the triangular caps $T_{n+1,j}$ such that 
$$T_{n+1,j} \subset T_{n,i} \in D_n$$
for some $i$. For each such $T_{n,i} \in D_n$ there are $2$ such triangular caps $T_{n+1,j}$.
\newline \newline
Then, as 
$$r_{n+1-6} < 4^{-(n+1-6}\Hm{1}(A_{7,\cdot}) < {{1}\over{2}}\Hm{1}(A_{n+1-6,\cdot})$$
it follows from Lemma $\ref{lem7}$ (2) that 
$B_{r_{n+1-6}}(x_{n+1-6})$
intersects at most $2$ triangular caps $T_{n+1,i}$. Thus removing these triangular caps means that removal of triangular 
caps of the $n+1$th level due to intersections with $\{B_{r_i}(x_i)\}_{i=1}^{n+1-6}$ has led to the removal of 
$$2N_n+2 \leq 2(2^{n-5}-2)+2 = 2^{(n+1)-5}-2$$
triangular cpas at the $n+1$th level. It follows that 
$$M_{n+1} = 2^{n+1}-N_{n+1} \geq 2^{n+1}-(2^{(n+1)-5}-2)$$
as required.
\newline \newline
Now 
\begin{eqnarray}
\sup\{r_i:i>n+1-6\} & < & \left({{1}\over{4}}\right)^{n+1-6}\Hm{1}(A_{7,\cdot}) \nonumber \\
& < & {{1}\over{2}}\Hm{1}(A_{n+1,\cdot}). \nonumber
\end{eqnarray}
Since, by construction 
$$B_{r_i}(x_i) \cap A_{n+1} = \emptyset \hbox{ for each }i>2+\sum_{j=0}^{n+1}2^j$$
and since for each $i \leq 2+\sum_{j=0}^{n+1}2^j$, 
$$x-i \in E(A_{n+1})$$ 
we know that apart from 
$$\bigcup_{i=1}^{n+1-6}B_{r_i}(x_i),$$
for which we have already removed the relevant triangular caps, $B_{r_i}(x_i)$ is a ball around an edge point of 
$A_{n+1}$ with $r_i < 1/2\Hm{1}(A_{n+1,\cdot})$. Thus by Lemma $\ref{lem7}$ (2) $B_{r_i}(x_i)$ 
can only meet the two triangular caps $T_{n+1,i}$ with intersection point $x_i$, it follows that for 
$$n+1-6 \leq i \leq 2 + \sum_{j=0}^{n+1}2^j$$
$B_{r_i}(x_i)$ consists of two straight lines of length $r_i$. Thus 
\begin{eqnarray}
\Hm{1}\left(\bigcup_{i=n+1-6}^{\infty}B_{r_i}(x_i) \cap A_{n+1} \right) 
& = & \Hm{1}\left(\bigcup_{i=n+1-6}^{2+\sum_{j=0}^{n+1}2^j}B_{r_i}(x_i) \cap A_{n+1} \right) \nonumber \\
& \leq & \sum_{i=n+1-6}^{2+\sum_{j=0}^{n+1}2^j}2r_i \nonumber \\
& < & \sum_{i=n+1-6}^{\infty}2^{-i} \nonumber \\
& < & 2^{-6}. \nonumber
\end{eqnarray}
Therefore
\begin{eqnarray}
\Hm{1}(L_{n+1}) & = &\Hm{1}(A_{n+1} \sim \cB) \nonumber \\
& \geq & \Hm{1}\left(A_{n+1} \sim \bigcup_{T_{n+1,i}\in D_n}T_{n+1,i}\sim\bigcup_{i=n+1-6}^{\infty}B_{r_i}(x_i)\right)
\nonumber \\
& \geq & \Hm{1}\left(A_{n+1} \sim \bigcup_{T_{n+1,i}\in D_n}T_{n+1,i}\right) - 
\Hm{1}\left(A_{n+1} \cap \bigcup_{i=n+1-6}^{\infty}B_{r_i}(x_i)\right) \nonumber \\ 
& \geq & \Hm{1}\left(A_{n+1} \sim \bigcup_{T_{n+1,i}\in D_n}T_{n+1,i}\right) - 2^{-6}. \nonumber
\end{eqnarray}
Since $\Hm{1}(A_{n+1} \cap T_{n+1,i})$ is constant over $i$ we know 
$$\Hm{1}\left(A_{n+1} \sim \bigcup_{T_{n+1,i}\in D_n}T_{n+1,i}\right) \geq {{M_{n+1}}\over{2^{n+1}}}\Hm{1}(A_{n+1})$$
so that 
$$\Hm{1}(L_{n+1}) \geq {{2^{n+1}-(2^{n+1-5}-2)}\over{2^{n+1}}}\Hm{1}(A_{n+1}) - 2^{-6}$$
completing the inductive step.
\newline \newline
We thus, most importantly have 
\begin{eqnarray}
\Hm{1}(L_n) & \geq &{{2^{n}-(2^{n-5}-2)}\over{2^{n}}}\Hm{1}(A_{n}) - 2^{-6} \nonumber \\
& \geq & (1-2^{-4})\Hm{1}(A_n) - 2^{-6} \nonumber 
\end{eqnarray}
for each $n \geq 7$, so that 
\begin{eqnarray}
\liminf_{n \rightarrow \infty}\Hm{1}(L_n) & \geq & \liminf_{n \rightarrow \infty}(1-2^{-4})\Hm{1}(A_n)-2^{-6} \nonumber \\
& = &-2^{-6} + (1-2^{-4})\liminf_{n \rightarrow \infty}\Hm{1}(A_n) \nonumber \\
& = & \infty. \nonumber 
\end{eqnarray}
Thus we can use the word limit and write 
$$\lim_{n \rightarrow \infty}\Hm{1}(L_n) = \infty$$
which is what was required to complete the proof. 
\end{proof}
\end{lem}\noindent
\begin{cor}\label{cor6}\thst
$\A$ and $\Ae$ are not weak locally $\Hm{1}$-finite. 
\pf
For $\A$ this follows directly from Corollary $\ref{cor5}$. \newline \newline
Now, suppose that $\Ae$ is weakly locally $\Hm{1}$-finite. Then for each $y \in \Ae$ there is a radius 
$\rho_y >0$ such that 
$$\Hm{1}(\Ae \cap B_{\rho_y}(y)) <\infty.$$
$\{B_{\rho_y}(y)\}_{y \in \Ae}$ is an open cover of $\Ae$ so that since $\Ae$ is compact there must exist 
a finite subcover $\{B_{\rho_{y_n}}(y_n)\}_{n=1}^Q$ of $\Ae$ and further we know that 
$$M:= \max\{\Hm{1}(\Ae \cap B_{\rho_{y_n}}(y_n)):1\leq n \leq Q\} < \infty.$$
It follows that 
\begin{eqnarray}
\Hm{1}(\Ae) & \leq & \Hm{1}\left(\Ae \cap \bigcup_{n=1}^QB_{\rho_{y_n}}(y_n)\right) \nonumber \\
& \leq & \sum_{n=1}^Q\Hm{1}(\Ae \cap B_{\rho_{y_n}}(y_n)) \nonumber \\
& \leq & QM \nonumber \\
& < & \infty. \nonumber 
\end{eqnarray}
This contradiction implies that there must exist a $y \in \Ae$ such that for each $\rho >0$ 
$$\Hm{1}(\Ae \cap B_{\rho}(y)) = \infty$$ 
and therefore that $\Ae$ is not weak locally $\Hm{1}$-finite. 
\end{proof}
\end{cor} \noindent
\section{Approximate $j$-Dimensionality of $\A$ and $\Ae$}
Having shown the measure theoretic properties of $\A$ and $\Ae$ that are required for them to be appropriate counter 
examples to (iv) (2), we now go on to show that $\A$ and $\Ae$ actually do satisfy the requirements of the Definition 
of (iv).
\begin{lem}\label{lem10}\thst
$\A$ and $\Ae$ satisfy property (iv).
\pf
Since $\Ae \subset \A$, proving that $\A$ satisfies (iv) is sufficient to prove the Lemma. We therefore 
proceed to prove that $\A$ is (iv).
\newline \newline
We first consider an arbitrary triangular cap, $T_{n,i}$ from somewhere in our construction. From the 
construction it is clear that it must be isosceles. From Lemma $\ref{lem7}$ and Construction $\ref{cons2}$ (particularly 
the constructed vertical heights, and Lemma $\ref{lem7}$ (1)) we see that it must have the two sorts of 
angles, $\psi(n,\varepsilon)$ and $\pi - 2\p{n}$, where, as in Definition $\ref{def18}$
$$\p{n} = \tan^{-1}\left({{2^{2-n}\varepsilon}\over{{{(1+n16\varepsilon^2)^{1/2}}\over{2^{n+1}}}}}\right).$$ 
so that we have
\begin{eqnarray}
\lim_{n \rightarrow \infty}\p{n} & = & \lim_{n \rightarrow \infty}\tan^{-1}
\left({{2^{2-n}\varepsilon}\over{{{(1+n16\varepsilon^2)^{1/2}}\over{2^{n+1}}}}}\right) \nonumber \\
& = & \tan^{-1}\left({{2^{2-n+n+1}\varepsilon}\over{(1+n16\varepsilon^2)^{1/2}}}\right) \nonumber \\
& = & 0. \label{e:l7p3psi0}
\end{eqnarray}
We now choose arbitrarily some $\delta > 0$ and $x \in A_{\varepsilon}$. 
We show that there is a $r_x$ such that for each $r \in (0,r_x]$ 
$$A_{\varepsilon} \cap B_{r}(x) \subset L_{x,r}^{\delta r}.$$
Since the endpoints of $A_{n,i}$ for each $n$, $i$ are not in $A_{\varepsilon}$, $x$ is not an endpoint so that 
we know from ($\ref{e:l7p3psi0}$)
that we can choose an $r_x>0$ such that 
$$B_{r_x}(x) \cap A \subset T_{n_0,j_0}$$
for some choice of $n_0 \in \mathbb{N}$ and $j_0 \in \{1,...,2^{n_0}\}$ 
and such that for all $n > n_0$
\begin{eqnarray}
\theta_{\delta} & := & tan^{-1}(\delta) \nonumber \\
& > & 3\p{n-1} + 2\p{n-2} \nonumber \\
& > &\p{n}. \nonumber 
\end{eqnarray}
Since $x \in A_{\varepsilon} \cap T_{n_0,j_0}$, for each $n > n_0$, $x \in T_{n,j(n)}$ for some 
$j(n) \in \{1,...,2^n\}$.
For each $r \in (0,r_x]$ we can therefore choose an $n_1 > n_0$ and $j_1 = j(n_1)$ such that 
$$\Hm{1}(A_{n_1,j_1}) \in [r/2,r).$$
We now consider $x$ as simply being some element of $T_{n_1,j_1}$ and set $L_{r,x}$ to be the affine 
space parallel to $A_{n_1,j_1}$ containing $x$.
\newline \newline
We now check that 
$$2^{2-n_1}\varepsilon > {{\delta r}\over{2}}.$$ 
First, we note that
\begin{eqnarray}
\delta  & > & tan (\p{{n_1}}) \nonumber \\
& = & {{2^{2-n_1+n_1+1}\varepsilon}\over{(1 + n16\varepsilon^2)^{1/2}}} \nonumber \\
& = & {{8\varepsilon}\over{(1 + n_116\varepsilon^2)^{1/2}}} \nonumber
\end{eqnarray}
which we get from the selection of $n_1$. Also, from the selection of $n_1$ with respect to $r$ that we have
\begin{eqnarray}
r & > & \Hm{1}(A_{n_1},j_1) \nonumber \\
& =& {{(1 + (n_1)16\varepsilon^2)^{1/2}}\over{2^{n_1}}} \nonumber 
\end{eqnarray}
so that 
\begin{eqnarray}
\delta r & > & {{8\varepsilon (1 + (n_1)16\varepsilon^2)^{1/2}}\over{(1 + n_116\varepsilon^2)^{1/2}2^{n_1}}} \nonumber \\
&\geq & 2^{3-n_1} \nonumber 
\end{eqnarray}
giving the desired inequality. 
This gives us that the vertical height of $T_{n_1,j_1}$ is less than half the diameter of 
the neighbourhood that we need around $L_{r,x}$ (that is $L_{r,x}^{\delta r} $). Thus
$$B_r(x) \cap T_{n_1,j_1} \subset L^{\delta R}_{r,x}.$$
It only remains to show that the remainder of $A_{\varepsilon} \cap B_r(x)$ is inside of an appropriate cone around 
$L^{\delta r}_{r,x}$. 
Since from the choice of $n_1$ with respect to $r$ we have that 
$$B_r(x) \subset \pi_x \left( O_{A_{n,i}}\left(\cup_{j:|i-j|\leq 1}A_{n,j}\right) \right) \times 
[-2\Hm{1}(A_{n,i}),2\Hm{1}(A_{n,i})].$$
Thus from Lemma $\ref{lem7}$ (2) it follows that the remainder of $A$ is contained in
$$\bigcup_{i:0<|i-j|<3}T_{n_1,i}$$  
so that it suffices to prove that these four caps are in the appropriate cone around $L^{\delta r}_{r,x}$.
We note that the union of these four caps is the subset of three $T_{n_1-1},k$ caps,
$$T_{n_1,j_1-2}\cap T_{n_1,j_1-1}\cap T_{n_1,j_1+1}\cap T_{n_1,j_1+2} \subset T_{n_1-1,j_1-1}\cap T_{n_1-1,j_1}\cap 
T_{n_1-1,n_1+1}.$$
By construction the maximal angle divergence from $L^{\delta r}_{r,x}$ that an edge on a neighbouring 
triangular cap of order $n_1$ is $2\p{n-1}$ and similarly for a triangular cap of order $n_1-1$, the 
maximal angular divergence is $2\p{n-2}$. Adding these together (which is actually worse than could possibly 
occur) we find that the maximal angle {\it requirement} for a cone around $L^{\delta r}_{r, x}$ is 
\begin{eqnarray}
2\p{n-1} + 2 \p{n-2} & < & 3\p{n-1} + 2 \p{n-2} \nonumber \\
& < &\theta_{\delta}. \nonumber 
\end{eqnarray}
It follows that we now have
$$B_r(x) \cap A \subset L^{\delta r}_{r,x}.$$
Since $x$ and $\delta$ were arbitrary, this shows that $A_{\varepsilon}$ 
has the fine weak $1$-dimensional $\varepsilon$-approximation property with local $r_x$ uniformity, 
(that is, it satisfies (iv)) and thus completes the proof. 
\end{proof}
\end{lem}\noindent
Corollary $\ref{cor5}$ and Lemma $\ref{lem10}$ 
allow us to provide the answer to question (iv) (2). We present this result formally in the following Theorem.
\begin{thm}\label{thm3}\thst
The answer to (iv) (2) is no.
\pf
From Lemma $\ref{lem10}$ $\A$ is a set that satisfies (iv) (2). Since, from Corollary $\ref{cor5}$ 
we know that $\A$ is not 
weak locally $\Hm{1}$-finite it follows that $\A$ is a counter example to the answer to (iv) (2) being 
yes. 
\end{proof}
\end{thm}\noindent
\section{Approximate $j$-Dimensionality of $\G$}
As previously discussed, the remainder of the answers to our definitions are completely dependent on showing 
that $\G$ satisfies (v). We show that this is true, or at least sufficiently true in the 
following Lemma. Sufficiently true here means that we can find an appropriate $\e$ such that 
$\G$ constructed with this $\e$ satisfues (v) for any given $\delta>0$. This is sufficient since definition (v) is 
dependent on some arbitrary but fixed $\delta$ unlike (iv) which requires $\delta$ to be able to be chosen arbitrarily 
for any set satisfying (iv). We show first that $\G$ satisfies (v) and then how the remainding classification follows.
\begin{lem}\label{lem11}\thst
 For all $\delta > 0$ there exists an $\e_{\delta}= \e_{\delta}(\delta) > 0$ such that 
$\Gamma_{\e_{\delta}}$ satisfies property (v) with respect to $\delta$.
\pf
Let $0<\e < 1/100$. We show, in fact, that there exists a function 
$$\delta (\e):\R \rightarrow \R$$ 
such that 
$$\lim_{\e \rightarrow 0}\delta (\e) = 0$$
such that $\G$ satisfies (v) with respect to $\delta (\e)$. It then follows that for all 
$\delta > 0$ there is an $\e_{\delta} >0$ such that $\delta(\e_{\delta}) < \delta$; $\G$ then satisfies 
(v) with respect to $\delta (\e_{\delta})$ and therefore with respect to $\delta$. 
\newline \newline
Let $w \in \G$ and $\rho \in (0,\rho_0] (= (0,1])$. Then,
as in Lemma $\ref{lem10}$, we know that there exists an $n \in \mathbb{N}$ such that $w \in T_{n,i}$ 
for some $i$ with $\Hm{1}(T_{n,i}) \in [\rho,2\rho )$. \newline \newline
Now from Lemma $\ref{lem7}$ (1) 
$$\{\psi^{T_{n,j}}_{T_{n,j+1}}\}_{j=i-1}^{j=i} < 2\p{0} < {{\pi}\over{16}}$$
so that
\begin{eqnarray}
O_{A_{n,i}}^{-1}(R_{n,i}) & = & O_{A_{n,i}}^{-1}(\pi_x ( O_{A_{n,i}}(\cup_{j:|i-j|\leq 1}A_{n,j}) ) \times 
[-2\Hm{1}(A_{n,i}),2\Hm{1}(A_{n,i})]) \nonumber \\
& \supset & O_{A_{n,i}}^{-1}([-(0.5\Hm{1}(T_{n,\cdot}) + 0.9\Hm{1}(T_{n,\cdot})),0.5\Hm{1}(T_{n,\cdot}) + 
0.9\Hm{1}(T_{n,\cdot})] \nonumber \\
& & \times [-2\Hm{1}(A_{n,i}),2\Hm{1}(A_{n,i})]) \nonumber \\
& \supset & O_{A_{n,i}}^{-1}([-\rho,\rho] \times [-2\Hm{1}(A_{n,i}),2\Hm{1}(A_{n,i})]). \nonumber
\end{eqnarray}
This implies that 
\begin{equation}
B_{\rho}(w) \subset O_{A_{n,i}}^{-1}(R_{n,i}).
\label{e:l7p5n1}
\end{equation}
From Lemma $\ref{lem7}$ (2) it follows that
$$\G \cap B_{\rho}(w) \subset \bigcup_{j:|i-j|<2}T_{n,j} \cup O_{A_{n,i}}^{-1}(R_{n,i})^c,$$
Since, from ($\ref{e:l7p5n1}$) 
$$B_{\rho}(w) \cap O_{A_{n,i}}^{-1}(R_{n,i})^c = \emptyset,$$
\begin{equation}
\G \cap B_{\rho}(w) \subset \bigcup_{j:|i-j|<2}T_{n,j}
\label{e:l7p5n2}
\end{equation}
and more importantly, that
$$\G \cap \left(B_{\rho(w)} \sim \bigcup_{j:|i-j|<2}T_{n,j}\right)  = \emptyset.$$
Since 
\begin{eqnarray}
\sup \{\pi_{y}(x):x \in O_{A_{n,i}}(T_{n,i}) \} & \leq & \varepsilon\Hm{1}(A_{n,i}) \nonumber \\
& \leq & \varepsilon 2\rho \nonumber
\end{eqnarray}
and since from Lemma $\ref{lem7}$ (2)
$$O_{A_{n,i}}\left(\bigcup_{j:|i-j|=2}T_{n,j}\right) \subset C_{4\p{0}}((0,0)) $$ 
so that  we have
\begin{eqnarray}
\sup \left\{ |\pi_y(z)|:z  \in  O_{A_{n,i}}\left( \bigcup_{j:|i-j|=1}T_{n,j} \cap B_{\rho}(w) \right) \right\}
& \leq & \sin (4\p{0})\rho \nonumber  
\end{eqnarray}
it follows that 
\begin{eqnarray}
\sup \{|\pi_y(z)|:z \in O_{A_{n,i}}(\G \cap B_{\rho}(w))\} & < & \sup \{2\e,sin(4\p{0})\}\rho \nonumber \\
& = & sin(4\p{0})\rho \nonumber
\end{eqnarray}
and thus by choosing $L_{w, \rho } || A_{n,i}$ we have  
$$\sup \{|\pi_{L_{w, \rho}}^{\perp}(z)|:z \in \G \cap B_{\rho}(w)\} < sin(4\p{0})\rho,$$
that is 
$$\G \cap B_{\rho}(w) \subset L_{w,\rho}^{sin(4\p{0})\rho}.$$
Thus $\G$ satisfies (v) for $\delta > sin(4\p{0})$. Which, since 
$$\lim_{\e \rightarrow 0} sin(4\p{0}) = 0,$$
by setting $\delta (\e) = sin(4\p{0})$, proves the lemma. 
\end{proof}
\end{lem} \noindent
The dimension of $\G$ follows from the work of Hutchinson \cite{hutch}. The proof is not trivial and so we do not 
present the proof here. We will however apply Hutchinsons proof regularly as a fundamental theorem of dimension 
to which we can reduce all of our investigations into the dimension of the generalised Koch Sets considered 
in Chapters 7 and 8. It is therefore important to state the Theorem and to show that $\G$ satisfies the 
conditions required for the Theorem to be applied.
\newline \newline
We first mention a result of Mandelbrot \cite{mandel} required to make sense of the result in \cite{hutch} that we use.
\begin{prop}\label{prop6}\thst
Let $\{r_i\}_{i=1}^N$ be a sequence of positive real numbers, then there exists a unique $D \in \R$ such that 
$$\sum_{i=1}^Nr_i^D = 1.$$
\end{prop}
With this $D$ we can now consider the appropriate result about dimension from \cite{hutch}.
\begin{thm}\label{thm4}\thst
If 
$$K = \bigcup_{i=1}^NS_i(K)$$
where $S_i$ are contraction mappings and if there exists an open set $O$ such that 
\begin{enumerate}
\item $0 \not= \emptyset$
\item $\bigcup_{i=1}^NS_i(O) \subset O $
\item $S_i(O) \cap S_j(O) = \emptyset$ whenever $i \not= j$
\end{enumerate}
Then if Lip$S_i =: r_i$ for each $1 \leq i \leq N$ and $D$ is the unique real number for which
$$\sum_{i=1}^Nr_i^D = 1$$
$$\hbox{dim}K = D.$$
\end{thm}\noindent
We can apply this Theorem directly to our case with $\G$ by appealing to Proposition 3 as follows. \newpage
\begin{lem}\label{lem12}\thst
For each $\e>0$, dim$\G >1$.
\pf
By Proposition $\ref{prop3}$ there exist, for each $\e>0$  
contraction maps $S_1$, $S_2$ with Lip$S_i = l(\e) > 1/2$ for each $i = 1,2$ and 
an open set $O$ such that the requirements of Theorem $\ref{thm4}$ are satisfied for $K = \G$.
\newline \newline
It follows that 
$$\sum_{i=1}^2(\hbox{Lip}S_i)^{\hbox{dim}\G} = 1.$$
That is 
$$2l^{\hbox{dim}\G} = 1,$$
or
$$\hbox{dim}\G = -{{\hbox{ln}2}\over{\hbox{ln}l}}>1.$$
\end{proof}
\end{lem}\noindent
We now have the tools to, and do in the following Theorem and Corollary, give the answers to our remaining definitions.
\begin{thm}\label{thm5}\thst
The answer to (v) (1) is no.
\pf
From Lemma $\ref{lem11}$ we know that $\G$ satisfies (v). Lemma $\ref{lem12}$ shows that dim$\G > 1$ and therefore 
that $\G$ is a counter example to the answer to (v) (1) being yes. 
\end{proof}
\end{thm}\noindent
\begin{cor}\label{cor7}\thst
The answer to the following definitions is no.
\begin{list}{}{}
\item (v) (2),
\item (ii) (1), and
\item (ii) (2).
\end{list} 
\pf
Since from Lemma $\ref{lem12}$ we know that 
the dimension of $\G$ is greater than $1$, it follows that $\G$ cannot be weak locally $\Hm{1}$-finite. Since 
Lemma $\ref{lem11}$ then shows that $\G$ satisfies (v), it follows that the answer to (v) (2) 
must be no.
\newline \newline
Since Property (v) is strictly stronger thatn Property (ii). Any set that satisfies (v) 
must also satisfy (ii). It then follows that $\G$ satisfies (ii) and thus in the same 
way that the answer to (v) (1) and (2) is no it follows that the answers to (ii) (1) 
and (2) is no. 
\end{proof}
\end{cor}\noindent
This completes the classification results that were the inital motivating aim for this work. We present again here 
a summary of the classification results:
\begin{eqnarray}
(i) & &(1) \ \  no \ \ (2) \ \  no \nonumber \\
(ii) & &(1)  \ \ no \ \  (2) \ \  no \nonumber \\
(iii) & &(1) \ \  yes \ \  (2) \ \ no \nonumber \\
(iv) & &(1) \ \ yes \ \ (2) \ \ no \nonumber \\
(v) & &(1) \ \ no \ \  (2) \ \ no \nonumber \\
(vi) & & (1) \ \ yes \ \  (2) \ \  no \nonumber \\
(vii) & & (1) \ \  yes \ \  (2) \ \  yes(weak)/no(strong) \nonumber \\
(viii) & & (1) \ \  yes \ \ (2) \ \ yes. \nonumber  
\end{eqnarray}
We next continue with results related to the fitting of the counter examples to the eight properties. In particular 
we show that $\A$ does indeed spiral in a sense that will be defined and we show that the counterexamples can be 
extended to higher dimensions.
\newline \newline
We have already seen that a rich tapestry of results follows from these more complicated examples. In the interest 
of finding as much interesting mathematics as possible that could arise from these sets we then in Chapters 7 and 
8 allow for generalisation of these sets and show various measure theoretic properties of the resulting sets.

\chapter{Miscellaneous Results}
In this section we present some further interesting and relevant results found in association with the study leading 
to the classification that we have presented, but that were not directly necessary in the classification. 
In particular we show that our present counter examples would not 
be sufficient for a $\delta$-fine version of property (v) and that $\A$ does not satisfy (vii), 
showing that there is no direct strength ranking of the 8 definitions in Definition $\ref{defa}$ since $\Lambda_{\delta}$ 
which satisfies (vii) does not satisfy (iii) which is satisfied by $\A$. Further,in the proof that 
$\A$ does not satisfy (vii) we see that the sets $\A$ do infact spiral infinitely finely 
in a sense that will become clearer.
\newline \newline
We also discuss how to extend the presented counter examples into higher dimensional counter examples. We show 
one such extension since the process of extending to higher dimensions remains the same for each of the sets.
\section{The Existence of Spiralling}
We start by showing that $\A$ does not satisfy (v) for each $\delta >0$. Similarly, but 
oppositely to Lemma $\ref{lem11}$ we show that there is also a function $\delta_1 : \R \rightarrow \R$ such that 
for each $\varepsilon > 0$, for each $\delta < \delta_1(\e)$, $\A$ does not satisfy (v) 
for $\delta$. Thus, although for each $\delta$ there is an $\A$ that fites, there is no $\e$ such that 
$\A$ satisfies (v) for every $\delta$, thus showing that $\A$ and indeed $\G$ are not 
sufficient as counter examples to any $\delta$-fine version of (v).
\begin{prop}\label{prop7}
There is a function 
$$\delta_1: \R \rightarrow \R^+$$
with $\delta_1(x) > 0$ for all $x>0$ such that for each $\e>0$ and all $\delta < \delta(\e)$ 
$\A$ does not satisfy (v) with respect to $\delta$.
\pf
First, we take 
$$y \in \A \cap B_{{{\varepsilon}\over{32}}}\left(\left({{1}\over{2}},2\varepsilon\right)\right)$$
and $\rho > {{3}\over{8}}$, say $\rho = \rho_0 = {{1}\over{2}}$ ($\rho_0 = {{1}\over{2}}$ as 
$B_{1/2}((1/2,0)) \supset A$.
It is not hard to see that we then have
$$\partial B_{\rho}(y) \cap int(T_{3,1}) \not= \emptyset$$
and 
$$\partial B_{\rho}(y) \cap int(T_{3,4}) \not= \emptyset$$
so that, more particularly 
$$B_{\rho}(y) \cap T_{3,1} \cap \A \not= \emptyset$$
and
$$B_{\rho}(y) \cap T_{3,4} \cap \A \not= \emptyset.$$
We now note that 
$$\sup_{x \in T_{3,i}}\pi_y(x) < \varepsilon $$
for $i = 1,4$ and that clearly 
$$\inf_{x \in B_{{{\varepsilon}\over{32}}}(({{1}\over{2}},2\varepsilon))}\pi_y(x) > {{63\varepsilon}\over{32}}.$$
It follows that a vertical gap between $y$ and points in $A \cap B_{\rho}(y)$ of atleast ${{31\varepsilon}\over{32}}$ 
exists both "to the left" of $y$ (that is points $z \in T_{3,1}$ where we must have $\pi_x(z) < \pi_x(y)$) and 
"to the right" of $y$ (similar to above, that is points $z \in T_{3,4}$ where we must have $\pi_x(z) > \pi_x(y)$). 
\newline \newline
Similarly clearly, we know that $\pi_x(z) > 0$ for all $x \in T_{3,1}$ (and also in fact $T_{3,4}$) and conversely 
we have $\pi_x(z) < 1$ for all $z \in T_{3,4}$ (and in fact, but unimportantly $T_{1,4}$). Also we have 
\begin{eqnarray}
\pi_x(y) & \in & \left({{1}\over{2}}-{{\varepsilon}\over{32}},{{1}\over{2}}+{{\varepsilon}\over{32}}\right) \nonumber \\
& \subset & \left({{63}\over{128}},{{65}\over{128}}\right) \nonumber 
\end{eqnarray}
since $\varepsilon < {{1}\over{4}}$.
\newline \newline
This means that in the best case any cone has less than a horizontal length of ${{65}\over{128}}$ to spread out to meet 
a set of vertical distance $${{31\varepsilon}\over{32}}$$ 
away from it's center.
\newline \newline
Supposing, to begin with, that $L_{y,\rho} || \R_x$ (that is $L_{y, \rho}$ is parallel to the horizontal axis)
then the cone angle must be, to cover the "best case mentioned above" at least
$$tan^{-1}\left({{\left({{31\varepsilon}\over{32}}\right)}\over{\left({{65}\over{128}}\right)}} \right).$$
Now, should  $L_{y, \rho}$ not be horizontal, we have that it is either positively or negatively sloped, but in 
either case, it continues to go through $y$. In the former case, we have that the cone angle estimate is imporoved 
for points in $T_{3,1}$ however, continuing to observe the $y=(63/128,63\varepsilon/32)$ case with a 
$z \in T_{3,4}$, it is clear that the resultant required cone angle for this $z$ can be 
no better than the cone angle required to include $z = (1,\varepsilon)$. 
We must therefore have that the minimum cone angle is no better than 
\begin{eqnarray}
\theta & = & tan^{-1}\left({{\left({{31\varepsilon}\over{32}}\right)}\over{\left({{65}\over{128}}\right)}} \right)
+ ||(L_{y,\rho}-y)-\R_x||_{G(1,2)} \nonumber \\
& > & tan^{-1}\left({{\left({{31\varepsilon}\over{32}}\right)}\over{\left({{65}\over{128}}\right)}} \right) \nonumber 
\end{eqnarray}
where $||\cdot||_{G(1,2)}$ denotes the norm on the grassman manifold of $1$-planes in $\R^2$.
A similar argument works considering points in $T_{3,1}$ in the case that $L_{y,\rho}$ is negatively sloped possibly 
improving the estimate for points in $T_{3,4}$. We therefore have that the cone angle 
$$tan^{-1}\left({{\left({{31\varepsilon}\over{32}}\right)}\over{\left({{65}\over{128}}\right)}} \right)$$
cannot be improved on, so that for any 
$$\delta < {{\left({{31\varepsilon}\over{32}}\right)}\over{\left({{65}\over{128}}\right)}} $$ 
$\A$ cannot satisfy (v) with respect to $\delta$. 
\newline \newline
Thus the function $\delta_1$ defined by 
$$\delta_1(x) := {{\left({{31x}\over{32}}\right)}\over{\left({{65}\over{128}}\right)}}$$
satisfies the requirements for the Proposition. 
\end{proof}
\end{prop}\noindent
To prove that $\A$ (and indeed $\G$) cannot satisfy (vii) irrespective of $\delta$, we have to 
show that although no spiralling occurs in the vicinity of any given point in $\A$ at a given approximation level, 
spiralling does indeed occur. 
\newline \newline
This means that for any point $x\in \A$, any radius $r>0$ and any potential approximating affine space, there 
exists a (smaller) 
triangular cap in $B_r(x)$ whose base is arbitrarily close to perpendicular to the approximating affine space.
\newline \newline
It then follows that an approriate choice of testing points and testing 
radius smaller than or equal to $r$ in such a triangular cap will allow us to show that 
for any $\delta<1$ $\A$ and indeed $\G$ cannot possibly satisfy (vii).
\begin{prop}\label{prop8}\thst
For each $\e > 0$ and $0< \delta < 1$, $\A$ does not satisfy Property (vii) with respect to $\delta$.
\pf
Let $\delta \in (0,1)$, $\e > 0$ and
$y \in \A$; then should $\A$ satisfy the definition then for each $\rho_y >0$ there would exist 
an affine space $L_{y, \rho_y}$ such that for all $x \in \A \cap B_{\rho_y}(y)$ and all $\rho<\rho_y$ 
$B_{\rho}(x) \cap \A \subset L_{y, \rho}^{\delta \rho}+x$.
\newline \newline
Now, since we are assuming that $\A$ satisfies the definition there must be a function, 
$$\phi:(0,1) \rightarrow (0,2\pi)$$ 
dependent only on $\delta$ which describes the cone outside of which no boundary points of a ball around 
a point $x \in B_{\rho_y}(y)$ are in $\A$. That is by defining
$$E_{\phi,\rho, x} := \left\{x\in \partial B_{\rho}(x) : tan^{-1}
\left({{\pi^{\perp}_{L_{y, \rho_y}}(x)}\over{\pi_{L_{y, \rho_y}}(x)}}\right) \geq \phi(\delta) \right\}$$
we have 
$$\A \cap E_{\phi(\delta),\rho,x} = \emptyset,$$
for all $x \in A \cap B_{\rho_0}(y)$ and also that there is a $\eta(\delta)>0$ such that for all 
$x\in A \cap B_{\rho_0}(y)$, $\rho \in (0, \rho_0]$ and all 
$z \in E_{{{\pi}\over{2}} - {{\pi - \psi(\delta)}\over{2}},\rho,x}$
\begin{equation}
B_{\rho \eta(\delta)}(z) \cap A = \emptyset.
\label{e:conddef(iii)}
\end{equation}
That is, around points in the central part of $E_{\phi,\rho, x}$ we can put a ball depending only on 
$\rho$ and $\delta$ that will be completely empty of $\A$.
\newline \newline
We observe that $y$ must be in some trianglular cap of the construction of $A$, $T_{n,i}$, for some $n$ and 
$i$, also such that $T_{n,i} \subset B_{\rho_0}(y)$. 
We make the nomenclatutorial choice to call the vertices of the triangular cap $T_{m,i}$ 
$a_{m,i}, l_{m,i}$, and $r_{m,i}$ chosen such that 
\begin{list}{}{}
\item $\pi_x \circ O_{A_{m,i}}(a_{m,i}) = 0$,
\item $\pi_x \circ O_{A_{m,i}}(l_{m,i}) < 0$, and 
\item $\pi_x \circ O_{A_{m,i}}(r_{m,i}) > 0$.
\end{list}
That is $a$ denotes the "top" vertici as we have previously defined, 
and $l$ and $r$ denote the identical "left" and "right" base angles.
\newline \newline
We now note that for each $k \in \mathbb{N}$ we have  
$$\psi_{T_{n,i} + z(i,k)}^{T_{n+k,2^ki + 4^{k-1}+2}} = \sum_{i=n}^{n+k}\p{i}$$
for some appropriate point $z(i,k) \in \R^2$.
\newline \newline
We now need some properties of the sequence $\{\p{i}\}_{i=1}^{\infty}$. First of all we recall that 
\begin{equation}
\lim_{i \rightarrow \infty} \p{i} = 0.
\label{e:theta0}
\end{equation}
And that we can specifically write that 
$$\p{n} = tan^{-1} \left({{8\varepsilon}\over{(1+16n\varepsilon^2)^{1/2}}}\right).$$
so that using the facts that 
$${{dtan}\over{dx}}(x) \geq 0,$$
and 
$${{dtan}\over{dx}}(x)|_{x=0} = 1$$
(and hence for sufficiently large $n$, $tan^{-1}(1/\varepsilon n) > 1/(2 \varepsilon n)$),
we get for any $n_0 \in \mathbb{N}$
\begin{eqnarray}
\sum_{n=n_0}^{\infty}\p{i} & = & \sum_{n=n_0}^{\infty} 
tan^{-1} \left({{8\varepsilon}\over{(1+16n\varepsilon^2)^{1/2}}}\right) \nonumber \\
& \geq & \sum_{n=n_0}^{\infty}tan^{-1}\left({{8\varepsilon } \over{4(1+n\varepsilon^2 )^{1/2}}}\right) \nonumber \\
& = & \sum_{n=n_0}^{\infty}tan^{-1}\left({{8\varepsilon } \over{4\varepsilon ({{1}\over{\varepsilon^2 }}+n)^{1/2}}}\right)
\nonumber \\
& \geq &\sum_{n=n_0}^{\infty}{{8\varepsilon } \over{8\varepsilon ({{1}\over{\varepsilon^2 }}+n)^{1/2}}} \nonumber \\
& = & \sum_{n=n_0}^{\infty}{{1}\over{({{1}\over{\varepsilon^2}} + n)^{1/2}}} \nonumber \\
& \geq & \sum_{n=n_0 + E}^{\infty}{{1}\over{n^{1/2}}} \nonumber \\
& > & \sum_{n=n_0 + E}^{\infty}{{1}\over{n}} \nonumber \\
& = & \infty, \nonumber
\end{eqnarray}
where $E$ denotes the smallest integer greater than or equal to $1/ \varepsilon^2$. 
It follows that there exists a sequence, $\{n_k\}$ , such that for each $k \in \mathbb{N}$ 
$$\sum_{i=n}^{n_k-1}\p{i} < {{2k\pi -\pi}\over{2}} < \sum_{i=n}^{n_k+1}\p{i}.$$
So that there is a triangular cap $T_{n_k,i(k)}$ (for the appropriate $i$ depending on $k$) such that 
\begin{eqnarray}
tan^{-1}\left({{\pi^{\perp}_{L_{y,\rho_0}}(r_{n_k,i(k)}-l_{n_k,i(k)})}\over
{\pi_{L_{y,\rho_0}}(r_{n_k,i(k)}-l_{n_k,i(k)})}}\right) & - & {{2k\pi - \pi}\over{2}} \nonumber \\
& < & \p{n_k-1}+\p{n_k}. \nonumber 
\end{eqnarray}
Thus, by $(\ref{e:theta0})$ there exists a $k \in \mathbb{N}$ such that 
\begin{eqnarray}
tan^{-1}\left({{\pi^{\perp}_{L_{y, \rho_0}}(r_{n_k,i(k)}-l_{n_k,i(k)})}\over
{\pi_{L_{y,\rho_0}}(r_{n_k,i(k)}-l_{n_k,i(k)})}}\right) & - & {{2k\pi - \pi}\over{2}} \nonumber \\
& < & \p{n_k-1}+\p{n_k} \nonumber \\
& < & {{{{\pi}\over{2}} - \phi(\delta)}\over{2}}. \nonumber
\end{eqnarray}
That is the triangular cap, $T_{n_k, i(k)}$ has the property that 
$$r_{n_k,i(k)} \in E_{\phi(\delta)+{{{{\pi}\over{2}} - \phi(\delta)}\over{2}},|r_{n_k,i(k)}-l_{n_k,i(k)}|,l_{n_k,i(k)}}.$$
The endpoints themselves are not in $A$, however, we can choose $x_l, z_r \in A$ such that 
$$|x_l - r_{n_k,i(k)}| < {{\eta(\delta) |r_{n_k,i(k)}-l_{n_k,i(k)}|}\over{2}}$$
and
$$|z_r - l_{n_k,i(k)}| < {{\eta(\delta) |r_{n_k,i(k)}-l_{n_k,i(k)}|}\over{2}}$$
so that there is a $x_r \in \partial B_{|r_{n_k,i(k)}-l_{n_k,i(k)}|}(x_l)$ such that 
$$z_r \in B_{\eta(\delta) |r_{n_k,i(k)}-l_{n_k,i(k)}|}(x_r).$$
Since, by our choice of triangular cap, $T_{n,i}$, $x_l \in B_{\rho_0}(y)$ and $|r_{n_k,i(k)}-l_{n_k,i(k)}| < \rho_0$ 
this contradicts $(\ref{e:conddef(iii)})$, proving the proposition since $\e$ and $\delta$ were chosen 
arbitrarily. 
\end{proof}
\end{prop}\noindent
\section{Higher Dimension Analogies of $\G$, $\A$ and $\Ae$}
We now come to the higher dimensional generalisations of the counterexamples.
\newline \newline
It is unfortunately trivial - unfortunate from the view of finding interesting mathematics - 
to generalise our counter examples to higher dimensions so that we obtain no further 
insite into how the structures of sets work. In each case we simply cross each set with either an interval or simply 
the plane of the required dimension, depending on whether or not we need the set to be bounded (as we do for 
$\rho_0$ uniformity properties). We show, as an example, how $\A$ is extended, and demonstrate how it continues to 
satisfy Property (iii).
\newline \newline
Suppose that we are taking $j$-dimensional approximations in $\R^{j+k}$.
We take 
$$S_{\e} = \A \times \R^{j-1} = \subset \R_A \times \R^{j-1} \times \R^{A_c},$$
where $\R_A = \R_{Ac} = \R$ but have been given names for notational convenience. $\R_{A}$ and $\R_{A_c}$ are 
identified with $\R$ and $\R^2/\R$ as we have been considering in the preceeding sections so that 
$\A \subset \R_A \times \R_{A_c}$.
Further $S_{\e}$ is constructed inside of 
$$\R^{j+k} = \R_A \times \R^{j-1} \times R_{A,c} \times \R^{k-1}.$$
We can thus see $S_{\e}$ as 
\begin{eqnarray}
S_{\e} & = & \{y = (y_1,x_2,...,x_{n-1},y_2,0,...,0):(y_1,y_2) \in \A \subset \R_A \times \R_{Ac}, x_i \in \R \}
\nonumber \\
& \subset & \R_A \times \R_x^{j-1} \times R_{A_c} \times \R_z^{k-1}, \nonumber 
\end{eqnarray}
where $\R_x^{j-1} \cong \R^{j-1}$ and $\R_z^{k-1} \cong \R^{k-1}$ are again notational conveniences denoting the 
dimensions along which the extension of $\A$ into $S_{\e}$ exist ($\R_x^{j-1}$), and the additional codimensions 
($\R_z^{k-1}$).
\begin{prop}\label{prop}\thst
$S_{\e}$ shows that the answer to (iii) (2) is no for arbitrary $j$.
\pf
There are two properties that we need to show that $S_{\e}$ has. That it has the fine weak $j$-dimensional 
$\delta$-approximation property, and that for each $x \in S_{\e}$ and $R>0$, $B_R^{j+k}(x) = +\infty$.
\newline \newline
First, to show that $S_{\e}$ has property (iii). We take 
arbitrary $y \in S_{\e}$ and $\delta > 0$. We now need only show that there exists a $j$-dimensional affine space, 
$L_{y,\rho}$ for each $\rho > 0$, such that 
$$S_{\e} \cap B_{\rho}(y) \subset L_{y,\rho}^{\delta \rho}.$$
We note that since $\A$ has property (iii), there exists for the chosen 
$\delta$ and $y$ a $1$-dimensional affine space $L_{(y_1,y_2),\rho}$ such that 
$$\A \cap B_{\rho}^2(y_1,y_2) \subset L_{(y_1,y_2),\rho}^{\delta \rho}.$$
It is therefore reasonable to take and test $L_{y, \rho} = L_{(y_1,y_2),\rho} \times \R_x^{j-1}$ as our affine space.
Clearly 
\begin{eqnarray}
S_{\e} \cap B_{\rho}(y) & = &(\A \cap \pi_{\R_A \times \R_{A_c}}(B_{\rho}(y)))\times \R_x^{j-1} \nonumber \\
& \subset & L_{(y_1,y_2),\rho}^{\delta \rho} \times \R_x^{j-1} \nonumber \\
& = & L_{y, \rho}^{\delta \rho}, \nonumber
\end{eqnarray}
which gives us that $S_{\e}$ has the appropriate property.
\newline \newline
To show that there is infinite measure in each ball $B_{R_1}(y)$ we take an $R_1 > 0$ and a $y \in S_{\e}$. 
Let $R = d(y,\partial S_{\e})$. We then get that 
\begin{eqnarray}
\Hm{j}(S_{\e} \cap B_{R_1}(y)) & \geq & \Hm{j}(S_{\e} \cap B_{R}(y)) \nonumber \\
& \geq & \Hm{j}(S_{\e} \cap ([-R/4,R/4]^{j+k}+(y_1,0,...,y_2,0,...,0))) \nonumber \\
& = & \Hm{1}(A \cap ([-R/4,R/4]^2 + (y_1,y_2)))\Hm{j-1}(\pi_{R_x^{j-1}}) \nonumber \\
& = & \Hm{1}(A \cap ([-R/4,R/4]^2 + (y_1,y_2)))\left({{R}\over{4}}\right)^{j-1} \nonumber \\
& = & + \infty, \nonumber 
\end{eqnarray}
showing that $S_{\e}$ is not weak locally $\Hm{j}$-finite. 
\end{proof}
\end{prop}\noindent

\chapter{Generalised Koch type sets and Relative centralisation of sets}
We turn now to the generalisation of the sets $\A$ and $\G$ which in our generalisations turn out to be two 
examples of the same sort of set. As already hinted at in Definition $\ref{def18}$ the generalisation can be 
seen as increasing the freedom with which the base angles of the triangular caps $\T{n,i}^A$ for a set $A$.
We allow this freedom in two differing strengths. Firstly that $\T{n,i}=T{n,j}$, $i,j\in\{1,...,2^n\}$ as in the 
construction of $\A$. Secondly that $\T{n,j}$ are allowed to vary freely over $n$ and $j$. A common restriction 
to the two variations is that $T_{n,i} \subset T_{m,j} \Rightarrow \T{n,i}\leq \T{m,j}$. That is, as we take 
triangular caps inside of previously constructed ones, the base angles reduce. The rate of reduction in seperate 
triangular caps may of course vary.
\newline \newline
It is clear that the second variation is a direct generalisation of the first. We keep them seperate however since 
the second allows more complications than the first and so some results are able to be presented in a stronger form 
for the first variation. 
\newline \newline
The original motivation for this investigation stems from an interest in the dimension of these sets. 
$\G$ and $\A$ are both examples of the first variation where for $\G$, $\T{n,i}^{\G}$ is constant over $n$ and $i$, 
whereas $\T{n,i}^{\A}$ varies by strictly decreasing to $0$ in $n$. The question being whether higher dimensions 
than (in this case) $1$ could only be reached with constant base angle as in $\G$. The answer turns out to be no.
Along with a presentation of this answer in both variations of our generalisation we present various other results 
concerning measure and rectifiability erlating to our generalisations.
\newline \newline
In this chapter we present the two main definitions of the sets in question and show their equivalence (both 
definitions will be used as which is more convenient in proofs that we present varies). We further show another 
characterisation of these sets in terms of a bijection from $\R$. We then present some general lemmas and background 
results necessary to present the main results concerning measure, rectifiability and dimension. The main results 
are then presented in the next and final chapter. 
\section{Equivalent Constructions of Koch Type Sets}
We start, quite naturally with definitions, equivalences and characterisations. First of all with a formal definition 
of the first variation of the generalisations.
\begin{def1}\label{def19}\thst
Suppose we can construct a set $B$ as follows: \newline \newline
Let $A_{0,1}$ be a base (a line in $\R^2$) and $T_{0,1}$ be a triangular cap on $A_{0,1}$ with vertical height 
$\e\Hm{1}(A_{0,1})$ with $\e<1/100$. Let $\p{0}$ be the base angles of $T_{0,1}$ and the two shorter sides of $T_{0,1}$ 
be named $A_{1,1}$ and $A_{1,2}$. We then construct two new triangular caps $T_{1,1}$ and $T_{1,2}$ on $A_{1,1}$ and 
$A_{1,2}$ with base angles $\p{1} \leq \p{0}$. We define 
$$A_0 = T_{0,1}$$ and 
$$A_1 = \bigcup_{i=1}^2T_{1,i}.$$
Then suppose we have $A_n = \cup_{i=1}^{2^n}T_{n,i}$ a union of $2^n$ triangular caps with base angles $\p{n}$ and 
$2^{n+1}$ shorter sides labelled $A_{n,i}$, $i \in \{1,...,2^{n+1}\}$. Then construct a triangular cap $T_{n+1,i}$ on 
each $A_{n+1,i}$ such that the base angles $\p{n+1}$ satisfy $\p{n+1} \leq \p{n}$. Define 
$A_{n+1}=\cup_{i=1}^{2^{n+1}}T_{n,i}$. Finally define 
$$B=\bigcap_{n=0}^{\infty}A_n.$$
We then call a set $A$ an {\bf $A_{\e}$-type} set whenever $A \in \{B, B \sim E(B)\}$.
\end{def1}\noindent
Then immediately we define the second variation.
\begin{def1}\label{def20}\thst
 Suppose we can construct a set $B$ as follows: \newline \newline
Let $A_{0,1}$ be a base (a line in $\R^2$) 
(for our purposes, provided that the line has non-infinite, non-zero length, it's position and length have no 
effect on the properties with respect to rectifiability, dimension, etc. and so without loss of 
generality we will generally assume that $A_{0,1}=[0,1] \subset \R$)
and $T_{0,1}$ be a triangular cap on $A_{0,1}$ with vertical height 
$\e\Hm{1}(A_{0,1})$ with $\e<1/100$. Let $\theta_0$ be the base angle of $T_{0,1}$ and the two shorter sides of $T_{0,1}$ 
be denoted $A_{1,1}$ and $A_{1,2}$. We then construct two new triangular caps $T_{1,1}$ and $T_{1,2}$ on $A_{1,1}$ and 
$A_{1,2}$ with base angles $\theta_{1,1}, \theta_{1,2} \leq \theta_0$. We define 
$$A_0 = T_{0,1}$$ and 
$$A_1 = \bigcup_{i=1}^2T_{1,i}.$$
Then suppose we have $A_n = \cup_{i=1}^{2^n}T_{n,i}$ a union of $2^n$ triangular caps with base angles $\theta_{n,i}$ and 
$2^{n+1}$ "shorter sides" (two per triangular cap) 
labelled $A_{n+1,i}$, $i \in \{1,...,2^{n+1}\}$. Then construct a triangular cap $T_{n+1,i}$ on 
each $A_{n+1,i}$ such that the base angles $\{\theta_{n+1,i}\}_{i=1}^{2^{n+1}}$ 
satisfy for each $i \in \{1,...,2^n\}$ 
$$\theta_{n,i} \geq \cases{\theta_{n+1,2i-1} \cr \theta_{n+1,2i}\cr}.$$ 
(i.e. the new base angles for each triangular cap are bounded by the base angle of the nth level that the new triangular 
cap is contained in).
\newline \newline
Define 
$A_{n+1}=\cup_{i=1}^{2^{n+1}}T_{n+1,i}$. Finally define 
$$B=\bigcap_{n=0}^{\infty}A_n.$$
We then call a set $A$ a {\bf Koch type set} whenever $A \in \{B, B \sim E(B)\}$. 
We denote the set of all such sets by $\cl{K}$.
\end{def1}\noindent
{\bf Remark:} In general any notation that can be considered in relation to some set $A \in \cl{K}$, for e.g.
$\T{n,j}$, $T_{n,j}$, etc., the superscript $A$ will denote association with the set $A$ when it may be unclear 
which set we are talking about. That is $T_{n,j}^A$ will denote the triangular cap $T_{n,j}$ associated 
with the construction of $A$.
\begin{def1}\label{def21}\thst
Let $A \in \cl{K}$. Then 
$$\tilde{A}_n^A := \bigcup_{i=1}^{2^n}A_{n,i}^A$$
\end{def1}\noindent
The second round of definitions for the two variations of generalisation are directly analogous to the original 
construction of $\A$ in that we consider $\tilde{A}_{n,j}$ sets instead of the $T_{n,j}$ sets.
\begin{def1}\label{def22}\thst
Suppose we can construct a set $B$ as follows: \newline \newline
Let $A_{0,1}$ be a base (a line in $\R^2$ of finite 
length) and $T_{0,1}$ be a triangular cap on $A_{0,1}$ with vertical height 
$\e\Hm{1}(A_{0,1})$ with $\e<1/100$. Let $\p{0}$ be the base angles of $T_{0,1}$ and the two shorter sides of $T_{0,1}$ 
be named $A_{1,1}$ and $A_{1,2}$. We then construct two new triangular caps $T_{1,1}$ and $T_{1,2}$ on $A_{1,1}$ and 
$A_{1,2}$ with base angles $\p{1} \leq \p{0}$. We define 
$$A_0 = T_{0,1}$$ and 
$$A_1 = \bigcup_{i=1}^2T_{1,i}.$$
Then suppose we have $A_n = \cup_{i=1}^{2^n}T_{n,i}$ a union of $2^n$ triangular caps with base angles $\p{n}$ and 
$2^{n+1}$ shorter sides labelled $A_{n,i}$, $i \in \{1,...,2^{n+1}\}$. Then construct a triangular cap $T_{n+1,i}$ on 
each $A_{n+1,i}$ such that the base angles $\p{n+1}$ satisfy $\p{n+1} \leq \p{n}$. Define 
$\tilde{A}_{n+1}=\cup_{i=1}^{2^{n+1}}A_{n,i}$. Finally define 
$$B=\overline{\bigcup_{n=0}^{\infty}\tilde{A_n}}\sim \bigcup_{n=0}^{\infty}\tilde{A_n}.$$
We then call a set $A$ an {\bf $A_{\e}$-type} set whenever $A \in \{B, B \sim E(B)\}$.
\end{def1}\noindent
Then immediately we define the second variation.
\begin{def1}\label{def23}\thst
 Suppose we can construct a set $B$ as follows: \newline \newline
Let $A_{0,1}$ be a base (a line in $\R^2$) 
(as previously, provided that the line has non-infinite, non-zero length, it's position and length have no 
effect on the properties with respect to rectifiability, dimension, etc. and so without loss of 
generality we will generally assume that $A_{0,1}=[0,1] \subset \R$)
and $T_{0,1}$ be a triangular cap on $A_{0,1}$ with vertical height 
$\e\Hm{1}(A_{0,1})$ with $\e<1/100$. Let $\theta_0$ be the base angle of $T_{0,1}$ and the two shorter sides of $T_{0,1}$ 
be denoted $A_{1,1}$ and $A_{1,2}$. We then construct two new triangular caps $T_{1,1}$ and $T_{1,2}$ on $A_{1,1}$ and 
$A_{1,2}$ with base angles $\theta_{1,1}, \theta_{1,2} \leq \theta_0$. We define 
$$A_0 = T_{0,1}$$ and 
$$A_1 = \bigcup_{i=1}^2T_{1,i}.$$
Then suppose we have $A_n = \cup_{i=1}^{2^n}T_{n,i}$ a union of $2^n$ triangular caps with base angles $\theta_{n,i}$ and 
$2^{n+1}$ "shorter sides" (two per triangular cap) 
labelled $A_{n+1,i}$, $i \in \{1,...,2^{n+1}\}$. Then construct a triangular cap $T_{n+1,i}$ on 
each $A_{n+1,i}$ such that the base angles $\{\theta_{n+1,i}\}_{i=1}^{2^{n+1}}$ 
satisfy for each $i \in \{1,...,2^n\}$ 
$$\theta_{n,i} \geq \cases{\theta_{n+1,2i-1} \cr \theta_{n+1,2i}\cr}.$$ 
(i.e. the new base angles for each triangular cap are bounded by the base angle of the nth level that the new triangular 
cap is contained in).
\newline \newline
Define 
$\tilde{A}_{n+1}=\cup_{i=1}^{2^{n+1}}A_{n+1,i}$. Finally define 
$$B=\overline{\bigcup_{n=0}^{\infty}\tilde{A_n}}\sim \bigcup_{n=0}^{\infty}\tilde{A_n}.$$
We then call a set $A$ a {\bf Koch type set} whenever $A \in \{B, B \sim E(B)\}$.
We denote the set of all such sets by $\cl{K}$.
\end{def1}\noindent
\begin{def1}\label{def24}
Let $A \in \cl{K}$ we then define the edge points of $A$, $E(A)$ by 
$$E(A):= \bigcup_{n=1}^{\infty}\bigcup_{i=1}^{2^n}E(T_{n,i}^A)$$
where $E(T_{n,i}^A)$ is as defined in Definition $\ref{def9}$.
\end{def1}\noindent
Before going on to show that these definitions are equivalent we need the following simple but important fact.
\begin{lem}\label{lem13}\thst
Let $A \in \cl{K}$. Then for any sequence $\{n,i(n)\}_{n \in \mathbb{N}}$ such that $T_{n,i(n)} \subset T_{n-1,i(n-1)}$ 
for each $n \in \mathbb{N}$
$$\lim_{n \rightarrow \infty}\Hm{1}(A_{n,i(n)})=0.$$
\pf
Since, by assumption $\T{0,1} <\pi/32$ and by construction $\T{n,i(n)}$ is decreasing in $n$. It follows from 
the inductive definition of the $A_{n,i(n)}$'s that 
\begin{eqnarray}
\Hm{1}(A_{n,i(n)}) & = & (cos\T{n-1,i(n-1)})^{-1}\Hm{1}(A_{n-1,i(n-1)}) \nonumber \\
& \leq & (cos\T{0,1})^{-1} \Hm{1}(A_{n-1,i(n-1)}) \nonumber \\
& = & C\Hm{1}(A_{n-1,i(n-1)}) \nonumber
\end{eqnarray}
where $C= (cos\T{0,1})^{-1} < 1$. It follows inductively that 
$$\Hm{1}(A_{n,i(n)}) \leq C^n\Hm{1}(A_{0,1}).$$
Since $\Hm{1}(A_{0,1})<\infty$ by construction, the result follows.
\end{proof}
\end{lem} \noindent
We now show that these definitions are equivalent.
\begin{prop}\label{prop10}\thst
Definition $\ref{def19}$ is equivalent to Definition $\ref{def22}$.
Definition $\ref{def20}$ is equivalent to Definition $\ref{def23}$.
\pf
We show these equivalences by showing that should $\cl{A}_2$ be defined as in Definition $\ref{def20}$ and $\cl{A}_1$ 
be defined as in Definition $\ref{def23}$ with the same $T_{n,i}$, $A_{n,i}$, $\T{n,i}$ etc. then
$$\cl{A}_1=\overline{\bigcup_{n=0}^{\infty}\tilde{A_n}}\sim \bigcup_{n=0}^{\infty}\tilde{A_n} = 
\left(\bigcap_{n=0}^{\infty}\bigcup_{i=1}^{2^n}T_{n,i}\right)-E(A)=\cl{A}_2-E(\cl{A}_2)=\cl{A}_2-E(\cl{A}_1).$$
That $E(\cl{A}_1)=E(\cl{A}_2)$ follows from Definition $\ref{def24}$ and the fact that the $T_{n,i}$ used for 
$\cl{A}_1$ and $\cl{A}_2$ are 
the same. We thus denote $E(A):=E(\cl{A}_1)=E(\cl{A}_2)$.
This will complete the proof since $E(A)$ is countable and thus $\Hm{1}(E(A))=0$.
\newline \newline
As in Lemma $\ref{lem4}$ we see that 
$$\cl{A}_1 + E(A)$$ 
is closed. Let 
$$x \in \cl{A}_2-E - \cl{A}_1.$$ 
then $d_x := d(x, \cl{A}_1 + E)>0$. \newline \newline
Now, for each $n \in \mathbb{N}$, $x \in T_{n,i}$ for some $i$ so that 
$$d(x,\cl{A}_1+E)<diam(T_{n,i}) = \Hm{1}(A_{n,i}).$$
From Lemma $\ref{lem13}$ we have
$$\lim_{n \rightarrow \infty}\Hm{1}(A_{n,i})=0.$$
Hence there is an $n_0 \in \mathbb{N}$ such that $diam(T_{n_0,j})=\Hm{1}(A_{n_0,j}) < d_x$ which implies 
$$d(x, \cl{A}_1+E)<\Hm{1}(A_{n_0,i}) <d_x = d(x,\cl{A}_1+E).$$
This contradiction implies 
$$\bigcap_{n=0}^{\infty}\bigcup_{i=1}^{2^n}T_{n,i} \subset \cl{A}_1 +E$$ and thus that 
$$\left(\bigcap_{n=0}^{\infty}\bigcup_{i=1}^{2^n}T_{n,i}\right)-E(A) \subset \cl{A}_1.$$
Next, it is clear from definition that 
$$\bigcup_{i=1}^{2^n}A_{n,i} \subset \bigcup_{i=1}^{2^n}T_{n,i}$$
so that
$$\bigcap_{n=1}^{\infty}\bigcup_{i=1}^{2^n}A_{n,i} \subset \bigcap_{n=0}^{\infty}\bigcup_{i=1}^{2^n}T_{n,i}=\cl{A}_2$$
and thus, since $\cl{A}_2=\cap_{n=0}^{\infty}\cup_{i=1}^{2^n}T_{n,i}$ is closed 
$$\overline{\bigcap_{n=1}^{\infty}\bigcup_{i=1}^{2^n}A_{n,i}} \subset 
\cl{A}_2.$$
Hence 
$$\cl{A}_1 \subset \overline{\bigcap_{n=1}^{\infty}\bigcup_{i=1}^{2^n}A_{n,i}} -E(A) \subset 
\cl{A}_2 -E(A).$$
Therefore,
$$\cl{A}_1 = \cl{A}_2-E(A).$$
\end{proof}
\end{prop}\noindent
Before moving on to the further characterisations of these sets we present another useful equivalence of representation 
concerning the constructional pieces of sets in $\cl{K}$
\begin{prop}\label{prop19}\thst
For any $A \in \cl{K}$ and any $\xi \in \R^{+}$
$$\bigcup_{n=0}^{\infty}\bigcap_{i(n,x):x\in \Lambda_{\xi +}} T_{n,i(n,x)} = 
\bigcup_{x \in \Lambda_{\xi +}}\bigcap_{n=0}^{\infty}T_{n,i(n,x)}.$$
\pf
First, suppose 
$$z \in \bigcup_{n=0}^{\infty}\bigcap_{i(n,x):x\in \Lambda_{\xi +}} T_{n,i(n,x)}.$$
Then, since $\theta_{n,i(n,x)}$ is decreasing in $n$ for all $x \in A$, so that $\theta_{n,i(n,x)} \geq \xi$ for 
all $ n \in \mathbb{N}_0$ and $x \in \Lambda_{\xi +}$, and since for all $n \in \mathbb{N}_0$ $z in T_{n,i(n,x0}$ for 
some $x \in \Lambda_{\xi +}$ it follows that $\theta_{n,i(n,z)} \geq \xi$ for each $n \in \mathbb{N}_0$ and thus 
$$\lim_{n \rightarrow \infty}\theta_{n,i(n,z)}\geq \xi$$
so that $z \in \Lambda_{\xi +}$. 
\newline \newline
Since clearly $z \in T_{n,i)n,z)}$ for each $n \in \mathbb{N}_0$ we can write 
$$\bigcup_{x \in \Lambda_{\xi +}}\bigcap_{n=0}^{\infty}T_{n,i(n,x)} \supset
\bigcup_{n=0}^{\infty}\bigcap_{i(n,x):x\in \Lambda_{\xi +}} T_{n,i(n,x)} \ni z,$$
so that 
$$\bigcup_{n=0}^{\infty}\bigcap_{i(n,x):x\in \Lambda_{\xi +}} T_{n,i(n,x)} \subset
\bigcup_{x \in \Lambda_{\xi +}}\bigcap_{n=0}^{\infty}T_{n,i(n,x)}.$$
For the other direction, suppose 
$$z \in \bigcup_{x \in \Lambda_{\xi +}}\bigcap_{n=0}^{\infty}T_{n,i(n,x)}.$$
Then for some $x\in \Lambda_{\xi +}$ $z \in \cap_{n=0}^{\infty}T_{n,i(n,x)}$ and therefore 
\begin{eqnarray}
z & \in & \bigcap_{n=0}^{\infty}T_{n,i(n,x)} \nonumber \\
& \subset & \bigcap_{n=0}^{\infty}\left(T_{n,i(n,x)} \cup \bigcup_{i(n,x):x \in \Lambda_{\xi +}}T_{n,i(n,x)}\right) 
\nonumber \\
& = & \bigcap_{n=0}^{\infty}\bigcup_{i(n,x):x\in \Lambda_{\xi +}}T_{n,i(n,x)}, \nonumber 
\end{eqnarray}
so that 
$$\bigcup_{x \in \Lambda_{\xi +}}\bigcap_{n=0}^{\infty}T_{n,i(n,x)} \subset
\bigcup_{n=0}^{\infty}\bigcap_{i(n,x):x\in \Lambda_{\xi +}} T_{n,i(n,x)}.$$
Combining these two inclusions gives the result. 
\end{proof}
\end{prop}

\section{Bijective Characterisation of Koch Type Sets}
We now show that sets in $\cl{K}$ can be characterised by a bijection from $\R$ into $\R^2$. Since some sets in $\cl{K}$ 
do not have dimension $1$ it may seem odd at first glance that such a bijection exists. By quoting the fact that there is 
a bijection between $\R$ and the Cantor set, however, we see that the concept is neither new or foreign in mathematics.
\newline \newline
We show also immediately that a certain level of control of the preimage can retained. To this end we need the following 
definition.
\begin{def1}\label{def26} \thst
Let $A \in \cl{K}$ and $n \in \mathbb{N}$, then the dyadic interval of order $n$ in $A_{0,0}$ 
(or simply, dyadic intervals of order $n$) are defined as the intervals 
$D_n$ of the form 
$$D_n = [l_{0,0} + i2^n(r_{0,0}-l_{0,0}),l_{0,0} + (i+1)2^n(r_{0,0}-l_{0,0})]$$
for some $i\in\{0,...,2^n-1\}$. For some chosen $j \in \{0,...,2^n-1\}$, the particular interval $D_{n,j}^A$ is defined as 
$$D_{n,j}^A = [l_{0,0} + j2^n(r_{0,0}-l_{0,0}),l_{0,0} + (j+1)2^n(r_{0,0}-l_{0,0})].$$
As per usual the superscript $A$ is dropped when the set is understood. 
\end{def1} \noindent
{\bf Remark:} Note that should $A_{0,0}$ be, or be normed to be $[0,1]$ on the real line, then the dyadic intervals in 
$A_{0,0}$ are simply the usual dyadic intervals.
\begin{prop}\label{prop11}\thst
Let $A \in \cl{K}$. Then there exists a sequence of Lipschitz functions $F_n:\R\mapsto \R^2$ ($F_n^A$ when which 
set $F_n$ is related to is not clear from the context) such that 
$$F_n(A_{0,1})=\tilde{A}_{n-1}.$$ 
Further there exists a bijection $\cl{F}$ ($\cl{F}^A$ when which set $\cl{F}$ is related to 
is not clear from the context) such that 
$$\cl{F}(A_{0,1})=A.$$ 
Additionally, denoting the relatively dyadic points of 
$A_{0,1}$ by $D$; 
\newline that is, for $\{x_1,x_2\}=E(A_{0,1})$, $x_1<x_2$,
$$D:=\{y: y=x_1 + (x_2-x_1)j2^{-n}, n \in \mathbb{N}, j \in \{0,...,2^{n}\}\};$$
we have
$$\cl{F}(D)=E(A).$$
Finally for each dyadic interval $D_{n,i}$ in $A_{0,0}$, 
$$F_{n}(D_{n,i}) = A_{n,i}$$
and 
$$\F1(D_{n,i}) \subset T_{n,i}.$$
\pf
Since the proof is the same for any $A_{0,1}$, we assume for notational convenience that $A_{0,1}=[0,1]$. In this case 
$D$ is also exactly the set of dyadic rationals in $[0,1]$. That is 
$$D=\{j2^{-n}:n\in \mathbb{N}, j\in \{0,...,2^n\}\}.$$
We will define $\F1$ as the limit of the $F_n$ functions, and then show that it is well defined and has the required 
properties. Firstly, we define $f_0:A_0 \rightarrow \R^2$ as 
$$f_0(y) = \cases{(y,tan\theta_{0,1} y) \ \ \hbox{ }\ \  y \in [0,1/2) \cr (y,tan\T{0,1} (1-y)) \ \ y \in [1/2,1] \cr}.$$
We see clearly that $f_0$ is a Lipschitz bijection between $A_0$ and $A_1$ 
(Since the graph of the function draws out the triangular cap $T_{0,1}^A$) with Lipschitz constant (and Jacobian) 
$Lipf_0 = Jf_0 \equiv cos\T{0,1}^{-1}$. We then similarly define for each 
$n \in \mathbb{N}$ $f_{n,i}:A_{n,i} \rightarrow \R^2$ by
$$f_{n,i}(y)=\cases{O_{A_{n,i}}^{-1}(\pi_x(O_{A_{n,i}}(y)),tan\T{n,i}(\pi_x(O_{A_{n,i}}(y))+\Hm{1}(A_{n,i})/2)) \ \ y \in 
I_1\cr O_{A_{n,i}}^{-1}(\pi_x(O_{A_{n,i}}(y)),tan\T{n,i}(\Hm{1}(A_{n,i})/2-\pi_x(O_{A_{n,i}}(y)))) \ \ y \in I_2 \cr},$$
where $I_1=O_{A_{n,i}}^{-1}([-\Hm{1}(A_{n,i})/2,0))$ and $I_2=O_{A_{n,i}}^{-1}([0,\Hm{1}(A_{n,i})/2])$. (Note that 
the $(1-y)$ factor in the definition of $f_0$ would change to some other appropriate constant should $A_{0,1}\not= [0,1]$.)
We note in particular that $f_{n,i}(A_{n,i}) \subset T_{n,i}$. Noting also 
that the two end points of $A_{n,i}$ stay fixed we can define $f_n:A_n \rightarrow \R^2$ by 
$$f_n(y) = f_{n,i}(y) \ \ y \in A_{n,i}.$$
We see then that similarly to the $f_0$ situation $f_n$ is a Lipschitz bijection between $A_n$ and 
$A_{n+1}$ with Lipschitz constant (and Jacobian in the case $A$ is an $\A$ type set)
$Lipf_n (=Jf_n) =\max_{1\leq i\leq 2^n}cos \T{n,i}^{-1}$.
\newline \newline
By writing for a collection of functions $\{g_i\}_{i=0}^n$
$$\circ_{i=0}^ng_i = g_n\circ g_{n-1}\circ ... \circ g_0$$
we can then define the Lipschitz bijection between $A_0$ and $A_{n+1}$, $F_n:A_0 \rightarrow \R^2$ by 
$$F_n = \circ_{i=0}^nf_i$$
which will then have Lipschitz constant (and Jacobian in the $\A$ type set case)
$$LipF_n = JF_n = \prod_{i=1}^n(\cos\theta_i )^{-1}.$$
This demonstrates the first claim.
\newline \newline
We can then propose a definition for $\F1$ and indeed we propose the definition of $\F1 : A_0 \rightarrow \R^2$ to 
be 
$$\F1(y) = \lim_{n \rightarrow \infty}F_n(y).$$
We need first of all to show that this function is well defined. To do this we suppose first of all that 
$$F_n(y) \in A_{n+1,i} \subset T_{n+1,i},$$
for some $i \in \{1,...,2^{n-1}\}$. Then
$$F_{n+1}(y) = f_{n+1,i}(y) \subset T_{n+1,i}.$$
Thus by induction, for each $n,k \in \mathbb{N}$
$$F_n(y) \in T_{n+1,i} \Rightarrow F_{n+k}(y) \in T_{n+1,i}$$
Then, From Lemma $\ref{lem13}$,
since $diam(T_{n,i})=\Hm{1}(A_{n,i})$, $diam(T_{n,i(n)}) \rightarrow 0$ as $n \rightarrow \infty$ 
for any sequence $\{n,i(n)\}_{n \in \mathbb{N}}$ and thus by setting the sequence $\{i(y,n)\}_{n \in \mathbb{N}}$ 
to be the sequence such that $y \in T_{n,i(y,n)}$ for each $n \in \mathbb{N}$ (so that it is always well defined, 
we choose arbitrarily $i(n)$ to be chosen such that $y =l_{n,i}$ for each $n$ for which $y$ is an edge point) 
it follows that for any $\e>0$ there 
is an $n_0 > 0$ such that for all 
$n,m=n+k > n_0$, $$d(F_n(y),F_m(y))<diam T_{n_0+1,i(y,n)} < \e$$ 
so that 
$\{F_n(y)\}$ is a Cauchy sequence in $\R^2$ and thus converges. It follows that $\F1$ is well defined. 
\newline \newline
We need still to show that $\F1$ is a bijective function between $A_0$ and $A$.
\newline \newline
We note firstly that for any $y \in A_0$ $F_n(y) \in A_n$ so that $\F1(y) \in \overline{\cup_{n=0}^{\infty}A_n}$ and 
thus 
$$\F1(A_0)= \bigcup_{y \in A_0}\F1(y) \subset \overline{\bigcup_{n=0}^{\infty}A_n}.$$
Now, since new edge points $a_{n,i}$ are by the definition of triangular caps always directly over the 
center of the base of the triangular cap, it follows that for all $e \in E$, $e = a_{n,i}$ for some $n \in \mathbb{N}$ 
and $i \in \{1,...,2^n\}$ and thus $e = F_n((2i-1)/2^{n+1})$. Since the set 
$\{(2i-1)/2^{n+1})\}_{n \in \mathbb{N},i\in \{1,...,2^n\}}=D$ the set of dyadic rationals, it follows that 
$\F1(D) = E$ which is a claim in our Proposition.
\newline \newline
Further, for all $y \in A_0 \sim D$, $F_{n+k}(y) \cap A_n \subset (A_{n+1+k} \sim E) \cap A_n = \emptyset$ 
for each $k \geq 0$ so that $\F1(y) \not\in A_n$ for all $n \in \mathbb{N}$. It thus follows that 
$$\F1(A_0 \sim D) \subset \overline{\bigcup_{n=0}^{\infty}A_n} - \bigcup_{n=0}^{\infty}A_n  \in \{A,A-E(A)\} \subset A$$
and thus that 
$$\F1(A_0) = \F1(A_0 \sim D) \cup \F1(D) \subset A \cup E = A.$$
We therefore have $\F1:A_0 \rightarrow A$. We now need to show that it is bijective. We first show, however, the 
final two claims that refer to the relationship of $\F1$ to the dyadic intervals of $A_{0,0}$.
\newline \newline
We quickly mention a sketch of a proof and motivation of the last two claims which will be more rigorously proven 
in the following result.
\newline \newline
From the above comment on the image of the dyadic rationals and the definition of $F_n$ for an $n \in \mathbb{N}$ 
it follows that for each $i \in \{1,...,2^{n+1}\}$
$$F_n\left(\left[{{i-1}\over{2^{n-1}}},{{i}\over{2^n}}\right]\right) = A_{n+1,i}.$$
This proves also our second last claim. Since, we have from definition that from each $n \in \mathbb{N}$ and any 
$x \in A_{0,0}$, $F_{n+1}(x)$ is in the same triangular cap $T_{n,i}$ as $F_n(x)$. It follows from induction 
that $\F1(x) \in T_{n,i}$. Since this is true for each $x_0$ such that $F_n(x_0) \in T_{n,i}$ and from the 
above this set is equal to $D_{n,i}$. It follows that $\F1(D_{n,i}) \subset T_{n,i}$ which is our final claim in the 
Theorem.
\newline \newline
Continuing with the proof of bijective we use the above proven important facts as follows. 
\newline \newline
Firstly, that should $x,z \in A_0$ with $x \not= z$ we then have that there 
is an $n \in \mathbb{N}$ such that $2^{1-n}\geq |x-z|>2^{-n}$ and thus there exist $i,j \in \{1,...,2^{n+2}\}$ with 
$4\geq |i-j|\geq 2$ and 
the property that $x \in [(i-1)2^{-n-2},i2^{n-2}]$ and $z \in [(j-1)2^{-n-2},j2^{n-2}]$. 
\newline \newline
It then follows that 
$F_n(x) \in T_{n+2,i}$ and thus, as above, $\F1(x) \in T_{n+2,i}$. Similarly $\F1(j) \in T_{n+2,j}$. 
\newline \newline
Since from Lemma $\ref{lem7}$ 
we know that for any $n \in \mathbb{N}$ $T_{n+2,i} \cap T_{n+2,j} = \emptyset$ whenever $4\geq |i-j|\geq 2$ it follows 
that $\F1(x) \not= \F1(z)$ and therefore that $\F1$ is injective.
\newline \newline
For surjectivity, we consider an arbitrary element $y \in A$. For all $n \in \mathbb{N}$, $y \in T_{n,i(y,n)}$ for 
some $i(y,n) \in \{1,...,2^n\}$. Then, again from 
$$F_n\left(\left[{{i-1}\over{2^{n-1}}},{{i}\over{2^n}}\right]\right) = A_{n+1,i}$$
we see that it is instructive to consider the intervals 
$$\F1^{-1}(A_{n,i(y,n)}) = [(i(y,n)-1)2^{-n},i(y,n)2^{-n}] =:D_{n,i(y,n)}.$$
Since $T_{n+1,i(y,n)} \subset T_{n,i(y,n)}$ for each $n$ it follows that $D_{n+1,i(y,n+1)} \subset D_{n,i(y,n)}$ 
for each $n$. We now observe $y_0 = \cap_{n=0}^{\infty}D_{n,i(y,n)}$. For this $y_0$ 
$$F_n(y_0) \subset F_n(D_{n,i(y,n)}) \subset A_{n,i(y,n)} \subset T_{n,i(y,n)}$$
for each $n$. Thus for each $n \in \mathbb{N}$, $|F_n(y_0) - y| \leq diam T_{n,i(y,n)}$. 
Since this diameter goes to zero as $n$ approached infinity it follows that 
$$\F1(y_0) = \lim_{n \rightarrow \infty}F_n(y_0)=y.$$
From well definedness and the arbitrariness of $y$ the surjectivity and thus bijectivity of $\F1$ follows.
\end{proof}
\end{prop}\noindent
We now show some results on the structure of $\F1$ which expand on the last two points of the previous results, as well 
as embellsihing the proof somewhat. We show that the function can be looked at as a function on each dyadic interval. 
$A$ in any given triangular cap is a bijection between $A$ in this cap and a dyadic interval in $A_{0,0}$. These results 
make it much easier to track images and pre-images and thus also to track how much measure has come from, or gone to 
where.
\begin{prop}\label{prop14}\thst
Let $A \in \cl{K}$ be constructed from a base $[0,1]$. Then when $\{F_n\}_{n=0}^{\infty}$ are the Lipschitz 
functions such that 
$$\F1^A:=\lim_{n \rightarrow \infty}F_n$$
pointwise on $A_{0,0}^A$ and 
$$F_n = \circ_{i=0}^nf_i$$
and writing $l_{ni}^A, a_{ni}^A$ and $r_{ni}$ as the edge point of $A_{n,i}^A$ adjoining $A_{n,i-1}^A$ (or (0,0) should 
$i=1$), the centerpoint of $A_{n,i}^A$ and the edge point of $A_{n,i}^A$ adjoining $A_{n,i+1}^A$ (or (1,0) should $i=2^n$) 
respectively. \newline \newline
Then for $n \in \mathbb{N}$ and $i \in \{1,...,2^n\}$ we have 
$$A_{n,i}^A = F_{n-1}([(i-1)2^{-n},i2^{-n}]), F_{n-1}^{-1}(A_{n,i}^A) = [(i-1)2^{-n},i2^{-n}],$$
$$l_{ni}=F_{n-1}(2^{-n}(i-1))$$
$$r_{ni}=F_{n-1}(2^{-n}i)$$
and that $F_{n-1}$ preserves relative distances. That is  for each $x,y \in [(i-1)2^{-n},i2^{-n}]$
$$|F_{n-1}(x) - F_{n-1}(y)| = p_{n,i}|x-y|$$
for some $p_{n,i} \in \R.$
\newline \newline
{\bf Remark:} Of the claims stated we are most interested in and thus emphasise 
$$A_{n,i}^A = F_{n-1}([(i-1)2^{-n},i2^{-n}]), F_{n-1}^{-1}(A_{n,i}^A) = [(i-1)2^{-n},i2^{-n}],$$
which gives us in essence a trace of the movement of a dyadic interval as it approaches the limit set $A$. With this we 
can follow the track either forward or backwards to identify which parts of $A$ or $A_{0,0}$ have positive 
measure given information about the measure of the other of $A$ and $A_{0,0}$. The other claims are stated here as an 
aid to proving the inductive step which is the key to the proof.
\pf
We prove the statement by induction on $n$. \newline \newline
From the definition of $A_{1,1}^A$, $A_{1,2}^A$ and the definition 
$$f_0(y) = \cases{(y,tan\theta_0 y) \ \ \hbox{ }\ \  y \in [0,1/2) \cr (y,tan\theta_0 (1-y)) \ \ y \in [1/2,1] \cr}.$$
it follows that $A_{1,1}^A=F_0([0,1/2]), A_{1,2}^A=F_0([1/2,1])$, that 
$F_0^{-1}(A_{1,1}^A)=[0,1/2], F_0^{-1}(A_{1,2}^A)=[1/2,1]$, that $F_0(0) = (0,0)=l_{11}$, that $F_0(1)=(1,0)=r_{12}$, 
and hence that $F_0(1/2) = r_{11}=l_{12}$. \newline \newline
We see also that the preservation of relative distances holds with $p_{1,i} \equiv tan\theta_{0,1}^A$ for $i=1,2$. The 
claim thus holds for $n=1$.
\newline \newline
Now suppose that the claim is true for each $n \leq m$ for some $m \in \mathbb{N}$. \newline \newline
We note that for any arbitrary $i \in \{1,...,2^{m+1}\}$ there is a $j \in \{1,...,2^m\}$ such that $i \in \{2j-1,2j\}$.
Now since $A_{m,j}^A = F_{m-1}([(j-1)2^{-m},j3^{-1}])$
\begin{eqnarray}
F_m([(j-1)2^{-m},j2^{-m}]) & = & f_m \circ F_{m_1}([(j-1)2^{-m},j2^{-m}]) \nonumber \\
& = & f_m(A_{m,j}). \nonumber
\end{eqnarray}
Since $F_{m-1}((j-1)2^{-m})=l_{lmj}^A$, $F_m(j2^{-m})=l_{rmj}^A$ and $F_m$ preserves relative distances we also have 
$$F_{m-1}((j-1)2^{-m} + 2^{-m-1})=l_{mj}^A,$$
and thus $\pi_x(O_{A_{n,j}^A}(F_{m-1}((j-1)2^{-m} + 2^{-m-1})))=0$.
Thus, again from relative distance preservation 
$$\pi_x(O_{A_{n,j}^A}(F_{m-1}([(j-1)2^{-m},(j-1)2^{-m}+2^{-m-1}])) = [-\Hm{1}(A_{n,j}^A)/2,0] =:I_1$$
and
$$\pi_x(O_{A_{n,j}^A}(F_{m-1}([(j-1)2^{-m}+2^{-m-1},j2^{-m}])) = [0,\Hm{1}(A_{n,j}^A)/2] =:I_2.$$
It follows then from the definition of $f_m|_{[(j-1)2^{-m},j2^{-m}]}$
$$f_{n,i}(y)=\cases{O_{A_{n,i}}^{-1}(\pi_x(O_{A_{n,i}}(y)),tan\theta_n(\pi_x(O_{A_{n,i}}(y))+\Hm{1}(A_{n,i})/2)) \ \ y \in 
I_1\cr O_{A_{n,i}}^{-1}(\pi_x(O_{A_{n,i}}(y)),tan\theta_n(\Hm{1}(A_{n,i})/2-\pi_x(O_{A_{n,i}}(y)))) \ \ y \in I_2 \cr},$$ 
and the definition of $A_{m+1,k}^A$, $k \in \{1,...,2^{-m-1}\}$ that 
$$A_{m+1,2j-1}^A = F_m([(2j-2)2^{-m-1},(2j-1)2^{-m-1}])$$
$$A_{m+1,2j}^A = F_m([(2j-1)2^{-m-1},2j2^{-m-1}])$$ 
and since we know $F_m$ is a bijection 
$$F_m^{-1}(A_{m+1,2j-1}^A) = [(2j-2)2^{-m-1},(2j-1)2^{-m-1}]$$
$$F_m^{-1}(A_{m+1,2j}^A) = [(2j-1)2^{-m-1},2j2^{-m-1}]$$
Further: 
$$F_m((2j-1)2^{-m-1}),F_m((2j-2)2^{-m-1}) \in E(A_{m+1,2j-1}^A)$$ and 
$$F_m((2j-1)2^{-m-1}),F_m(2j2^{-m-1}) \in E(A_{m+1,2j}^A)$$
from which it must therefore follow that 
$$F_m((2j-1)2^{-m-1}) = r_{(m+1)(2j-1)}^A=l_{l(m+1)(2j)}^A$$
$$F_m((2j-2)2^{-m-1}) = l_{(m+1)(2j-1)}^A$$
$$F_m(2j2^{-m-1}) = r_{(m+1)(2j)}^A.$$
Further since $F_{m-1}|_{[(2j-1)2^{-m},2j2^{-m}]}$ preserves relative distance with 
$|F_{m-1}(x) - F_{m-1}(y)| = p_{m-1,j}|x-y|$ for all $x,y \in [(2j-1)2^{-m},2j2^{-m}]$
from the definition of $f_m(y)$ and $F_m = f_m \circ F_{m-1}$ 
it follows that $F_m$ preserves relative distances on 
$[(2j-2)2^{-m-1},(2j-1)2^{-m-1}]$ and $[(2j-1)2^{-m-1},2j2^{-m-1}]$ with 
$$p_{m,2j-1} = p_{m,2j}= (tan\theta_{n,j}^A)p_{m-1,j}.$$
By substituting in $i$ for $2j-1$ or $2j$ as necessary it follows that all required preoperties are satisfied for $m+1$ 
with the choice of $i \in \{1,...,2^{m+1}\}$. Since the choice of $i$ was arbitrary this completes the inductive step 
and thus the proof. 
\end{proof}
\end{prop}\noindent
\section{Further Characterisations and Properties of Sets in $\cl{K}$}
Equiped with these results we are able to give a list of nomenclaturial definitions that will be instrumental in 
describing our results.
\begin{def1}\label{def25}\thst
Let $A\in \cl{K}$, then we write
$$\tilde{\theta}_{x_0}^A = \tilde{\theta}_{x}^A = \lim_{n \rightarrow \infty}\theta^A_{n,i(n,x)}$$
and define the functions $\tilde{\Pi}^A,\tilde{\Pi}^A_n, \tilde{\Pi}^A_{n,i}:\R \rightarrow \R$ by
$$\tilde{\Pi}^A(x)=\prod_{i=0}^{\infty}(cos\theta_{n,i(n,x)}^A)^{-1}$$
$$\tilde{\Pi}_n^A(x)=\prod_{i=0}^{n}(cos\theta_{n,i(n,x)}^A)^{-1}$$
and 
$$\tilde{\Pi}^A_{n,i} = \tilde{\Pi}^A_n(x),$$
for any $x \in A \cap T_{n,i}^A$.
The superscript $A$ is dropped when the set $A$ is understood.
\newline \newline
Further 
$$\Lambda_m :=\{x\in A: \tilde{\Pi}(x) \leq m\}$$
$$\Lambda_m^{-1}:=\F1^{-1}(\Lambda_m)$$
$$\Lambda_{m+} :=\{x\in A: \tilde{\Pi}(x) \geq m\}$$
$$\Lambda_{m+}^{-1}:=\F1^{-1}(\Lambda_{m+})$$
$$\Lambda_{\infty}:=\{x\in A:\tilde{\Pi}(x)=\infty\}$$
$$\Lambda^{-1}_{\infty} := \F1^{-1}(\Lambda_{\infty})$$
Also we introduce 
$$i(n,x):= \mathbb{N} \times A_0 \rightarrow \mathbb{N}$$
defined by 
$$i(n,x):= \{i \in \{1,...,2^n\}:x \in T_{n,i}^A\}$$
Also, we define for each $a \in \R$
$$\Upsilon_a^{-1} :=\{x\in A_0:\tilde{\theta}^A_x \leq a\}$$
$$\Upsilon_a :=\F1(\Upsilon_a^{-1})$$
$$\Upsilon_{a+}^{-1} :=\{x\in A_0:\tilde{\theta}^A_x \geq a\}$$
and
$$\Upsilon_{a+} :=\F1(\Upsilon_{a+}^{-1})$$
As with the other notations, when which $A \in \cl{K}$ we are referring to is unclear we add a superscript $A$, for 
example $(\Lambda^{-1}_{\infty})^A$.
\end{def1}\noindent
Two further definitions relating to sets being used will now be presented. Firstly a variant of the angle between sets , 
and then a generalisation of the $i(n,x)$ notation.
\begin{def1}\label{def28}\thst
Let $L_1,L_2$ be any two stright lines in $\R^2$ and $L^1, L^2$ be the extensions of these lines to simply connected 
lines of infinite length in both directions. We then denote the smaller of the two types of angles that occur at the 
intersection of $L^1$ and $L^2$ by $\psi^{L_1}_{L_2}$.
\end{def1} \noindent
\begin{def1}\label{def29}\thst
Let $A \in \cl{K}$ and $B \subset A_{0,0}$. suppose for a $n \in \mathbb{N}$, $i(n,x)$ is uniform for all $x \in B$. 
Then we will sometimes for convenience denote this common value $i(n,B)$.
\end{def1} \noindent
Further notations will occasionally be used, but not regularly and so will be defined as they are used. We continue now with 
further definitions and properties relating to the above terms and $\F1$ which will be necessary in our main results 
concerning sets in $\cl{K}$ which will be presented in the next chapter.
\begin{def1}\label{def27}\thst
We define, for $r \in \R$, the collection $A^r$ by 
$$A^r:=\{A:A \hbox{ is an }A_{\e} \hbox{ type set and } \tilde{\theta}^A_{\cdot}=r\}.$$
\end{def1}\noindent
We now state formally, to connect to the previous work, which $A^r$ sets that our previous sets $\G$ and $\A$ are 
members of. 
\begin{prop}\label{prop12}\thst
$\G \in A^{tan^{-1}(2\e)}$ and $\A \in A^0$.
\pf
that $\G \in A^{tan^{-1}(2\e)}$ follows from the definition of $\G$ since we can calculate from the construction that 
$\theta_{n,\cdot}^{\G} \equiv tan^{-1}(2\e)$. Since $\theta^A_{n,\cdot}$ is constant and from the proof of Lemma 
$\ref{lem10}$ $\lim_{n \rightarrow \infty}\theta^A_{n,\cdot}=0$ it follows that $A \in A^0$.
\end{proof}
\end{prop} \noindent
We now wish to investigate some of the properties possesed by $\F1$ and resultant from the definitions that we have 
just made. We first look at three results concerning the $\theta_{n,i}$. We see that the stretch (and when $\F1$ has 
appropriate properties the Jacobian) that occurs to each $\tilde{A_n}$ is described by a product of the base angles. 
Secondly we consider a convergence equivalence of this stretch factor to a convergence of the sum, which can be thought 
of as a test of whether a set $A \in \cl{K}$ spirals infinitely or not. Finally we look at the first of several 
results we have concerning the density of $A$ around the image of a considered point in $A_{0,0}$.
\begin{lem}\label{lem14}\thst
For any $A \in \cl{K}$, $n \in \mathbb{N}$ and $i \in \{1,...,2^n\}$
$$\Hm{1}(A_{n,i}^A) = {{\Hm{1}(A_{0,0}^A)}\over{2^n}}\prod_{j=1}^n{{1}\over{cos(\theta_{j,D_{j,i}}^A)}}.$$
\pf
By considering the right angled triangle consisting of $A_{n,i}^A$, half of the base $A_{n-1,j}^A$ of the triangular cap 
in which $A_{n,i}^A$ arises and the line connecting the ends that don't meet, we see that 
$${{\Hm{1}(A_{n-1,j}^A)}\over{2\Hm{1}(A_{n,i}^A)}} = cos(\theta_{n,D_{n,i}}^A)$$ 
so that
$$\Hm{1}(A_{n,i}^A) = {{\Hm{1}(A_{n-1,j}^A)}\over{2cos(\theta_{n,D_{n,i}}^A)}}.$$
Thus repeating this step inductively we get 
\begin{eqnarray}
\Hm{1}(A_{n,i}^A)& =&  {{\Hm{1}(A_{n-1,j}^A)}\over{2cos(\theta_{n,D_{n,i}}^A)}} \nonumber \\
& = & (cos(\theta_{n,D_{n,i}}^A)^{-1})(cos(\theta_{n-1,D_{n-1,i}}^A))^{-1}{{1}\over{4}}\Hm{1}(A_{n-2,\cdot}^A) \nonumber \\
& = & ... \nonumber \\
& = & \left(\prod_{j=0}^{n-1}(cos(\theta_{j,D_{j,i}}^A))^{-1}\right)\Hm{1}(A_{0,1}^A) \nonumber 
\end{eqnarray}
as required. \end{proof}
\end{lem}\noindent
\begin{prop}\label{prop13}\thst
Let $A \in A^0$ and $x \in A$. Then
$$\prod_{n=0}^{\infty}(cos(\theta_{n,i(n,x)}^A))^{-1} < \infty \Longleftrightarrow 
\sum_{n=0}^{\infty}(\theta_{n,i(n,x)}^A)^2 < \infty.$$
\pf
We first show that the claim is true for a sequence $\{\theta_{n,i(n,x)}^A\}$ composed of entirely sufficiently small 
elements.
Where what sufficiently small entails will be shown in the proof. 
\newline \newline
Let $M:=\prod_{n=0}^{\infty}(cos(\theta_{n,i(n,x)}^A))^{-1}$ and note 
\begin{eqnarray}
ln(M) & = & ln\left(\prod_{n=0}^{\infty}(cos(\theta_{n,i(n,x)}^A))^{-1}\right) \nonumber \\
& = & \sum_{n=0}^{\infty}ln((cos(\theta_{n,i(n,x)}^A ))^{-1}) \nonumber \\
&&\hbox{ so that using a Taylor expansion for $ln$ around $1$ we have } \nonumber \\
&=&\sum_{n=0}^{\infty}\sum_{j=1}^{\infty}{{(-1)^{j+1}(cos\theta_{n,i(n,x)}^A)^j}\over{j}}\left({{1}\over{cos\theta_{n,i(n,x)}^A}}-1\right)^j 
\nonumber \\
& = & \sum_{n=0}^{\infty}\sum_{j=1}^{\infty}{{(-1)^{j+1}(1-cos\theta_{n,i(n,x)}^A)^j}\over{j}}. \nonumber 
\end{eqnarray}
So that since
$$\lim_{x \rightarrow \infty}xsinx = 2\lim_{x\rightarrow \infty}1-cosx$$
we have that for sufficiently small $\theta_{n,i(n,x)}^A$ that 
$$(1-cos\theta_{n,i(n,x)}^A) \in (c_1\theta_{n,i(n,x)}^A sin\theta_{n,i(n,x)}^A,c_2\theta_n^A sin\theta_{n,i(n,x)}^A)$$ 
for $0<c_1<c_2<1$ and that 
$$c_2{{(\theta_{n,i(n,x)}^A)^j sin^j \theta_{n,i(n,x)}^A}\over{j}} < c_1{{1}\over{2}}{{(\theta_{n,i(n,x)}^A)^{j-1} sin^{j-1} 
\theta_{n,i(n,x)}^A}\over{j-1}}$$
and thus that
\begin{eqnarray}
ln(M) & = & \sum_{n=0}^{\infty}\sum_{j=1}^{\infty}{{(-1)^{j+1}(1-cos\theta_{n,i(n,x)}^A)^j}\over{j}} \nonumber \\
&\in &\left(\sum_{n=0}^{\infty}c_1\theta_{n,i(n,x)} sin\theta_{n,i(n,x)}^A - {{c_2(\theta_{n,i(n,x)}^A)^2sin^2\theta_{n,i(n,x)}^A}\over{2}},
\sum_{n=0}^{\infty}c_1\theta_{n,i(n,x)}^A sin\theta_{n,i(n,x)}^A \right). \nonumber 
\end{eqnarray}
Since 
$$\lim_{x \rightarrow \infty}x^2 = \lim_{x \rightarrow \infty}xsinx$$
we have, again for sufficiently small $\theta_{n,i(n,x)}^A$ that
\begin{eqnarray}
ln(M)&\in &\left(\sum_{n=0}^{\infty}c_1\theta_{n,i(n,x)} sin\theta_n^A - {{c_2(\theta_{n,i(n,x)}^A)^2sin^2\theta_{n,i(n,x)}^A}\over{2}},
\sum_{n=0}^{\infty}c_1\theta_{n,i(n,x)}^Asin\theta_{n,i(n,x)}^A\right)\nonumber \\
& \in & \left(\sum_{n=0}^{\infty}{{c_1}\over{2}}(\theta_{n,i(n,x)}^A)^2-{{c_2(\theta_{n,i(n,x)}^A)^2sin^2\theta_{n,i(n,x)}^A}\over{2}},
\sum_{n=0}^{\infty}2c_2(\theta_{n,i(n,x)}^A)^2 \right)\nonumber \\
& = & \left(\sum_{n=0}^{\infty}c_3(\theta_{n,i(n,x)}^A)^2 ,\sum_{n=0}^{\infty}c_2(\theta_{n,i(n,x)}^A)^2\right) \nonumber
\end{eqnarray}
for an appropriate $0<c_3<c_2$. It follows that
$$M \in \left(e^{\sum_{n=0}^{\infty}c_3(\theta_{n,i(n,x)}^A)^2},e^{\sum_{n=0}^{\infty}c_2(\theta_{n,i(n,x)}^A)^2}\right)$$
and thus for $\{\theta_{n,i(n,x)}^A\}$ being comprised of sufficiently small terms we have 
$$\prod_{n=0}^{\infty}(cos(\theta_{n,i(n,x)}^A))^{-1} < \infty \Longleftrightarrow \sum_{n=0}^{\infty}(\theta_{n,i(n,x)}^A)^2 < \infty.$$
The general case follows from noting that since $A \in A^0$, $\theta_{n,i(n,x)}^a \rightarrow 0$ and thus for sufficiently 
large $n$ the tail will be a sequence of sufficiently small $\theta_{n,i(n,x)}^A$. Since $\theta_{n,i(n,x)}^A<\pi/32$ in all cases, 
the finite number of terms at the begining of the sequnce will be a finite multiplying or adding factor for both 
sequences and thus will not affect convergence. 
\end{proof}
\end{prop}\noindent
We now present the first of three results that will be presented addressing the density of points in an $A \in \cl{K}$. 
The density is important as it will be the key to the existence or non-existence of approximate tangent spaces to 
$A$, and therefore an essential ingredient in discussing the rectifiability of sets in $\cl{K}$.
\begin{cor}\label{cor7}\thst
Let $A$ be an $\A$ type set. Then
$$\Theta^1(\Hm{1},A,y) \geq {{1}\over{2}}\lim_{j \rightarrow \infty}\prod_{n=j}^{\infty}(cos\theta_{n,\cdot}^A)^{-1}.$$
In particular, for $\A$ type sets $A$ such that $\prod_{n=0}^{\infty}(cos\theta_{n,\cdot}^A)^{-1}=\infty$, 
$\Theta^1(\Hm{1},A,y) =\infty$ for all $y \in \bar{A}$.
\pf
Let $\rho>0$ and $y \in \bar{A}$, then there is a $j \in \mathbb{N}$ and $i\in \{1,...,2^j\}$ with 
$T_{j,i}^A \subset B_{\rho}(y)$ and $\Hm{1}(A_{j,i}^A) \geq \rho$. From the proof of Lemma $\ref{lem14}$ we know 
$$\Hm{1}(A \cap T_{j,i}^A)=\Hm{1}(A_{j,i}^A)\prod_{n=j}^{\infty}(cos\theta_{n,\cdot}^A)^{-1}.$$
so that 
$$\Hm{1}(A \cap B_{\rho}(y)) \geq {{\Hm{1}(A_{j,i}^A)\prod_{n=j}^{\infty}(cos\theta_{n,\cdot}^A)^{-1}}\over{2\rho}} 
\geq{{1}\over{2}}\prod_{n=j}^{\infty}(cos\theta_{n,\cdot}^A)^{-1}$$
and thus 
$$\Theta^1(\Hm{1},A,y)=\lim_{\rho \rightarrow 0}{{\Hm{1}(A \cap B_{\rho}(y))}\over{2\rho}} > 
{{1}\over{2}}\lim_{j \rightarrow \infty}\prod_{n=j}^{\infty}(cos\theta_{n,\cdot}^A)^{-1}.$$
\end{proof}
\end{cor}\noindent
\section{Properties of The Bijective Functions}
We now examine some important properties of the functions $F_n$ and the function $\F1$. 
In order to work with $\F1$ properly we must first check that it has some basic properties. 
We show that the function $\F1$ is continuous and measurable. We show that images of compact sets are compact. We show 
that positive measure is preserves. A less well behaved, but nonetheless important property, is then to show that 
under conditions on $\Lambda_{\infty}^{-1}$ sets of positive measure have images of infinite measure. We additionally 
prove the $A \in \cl{K}$ version of Corollary $\ref{cor7}$. First of all, however, we prove that parts of the limit 
function $\F1$ can be expressed as Lipschitz functions. Recalling that $\tilde{\Pi}^A_{\cdot}$ can be seen as the 
stretching (or Jacobian) factor of $\F1$ it would seem sensible that when this is bounded, we are actually looking at a 
Lipschitz function. We show that this is true after defining how we make bounds. We make bounds by simply looking 
at the restriction of the function to pre-image sets on which $\tilde{\Pi}^A$ is bounded.
\begin{def1}\label{def30}\thst
Let $A \in \cl{K}$, then we define 
$$\F1_{m}:=\F1|_{\Lambda_m^{-1}}.$$
\end{def1} \noindent
\begin{lem}\label{lem15}\thst
For $m \in \R$, $\F1_m:=\F1|_{\Lambda_m^{-1}}$ is Lipschitz with $Lip \F1_m \leq Cm^2$.
\pf
Let $x,y \in \Lambda_m^{-1}$, and without loss of generality let $y<x$. there are then two cases to consider
\begin{enumerate}
\item $\{ty + (1-t)x : t\in [0,1]\} \subset \Lambda_m^{-1}$,
\item otherwise.
\end{enumerate}
Case 1 is the simpler. In this case we have $\F1_m|_{[y,x]} = \F1|_{[y,x]}$ and $\tilde{\Pi}(z) \in [y,x]$. It follows from 
the construction of the $F_n$ from which $\F1$ is defined as a limit that 
$$d(F_n(y),F_n(x))\leq md(y,x), \hbox{ for all } n \in \mathbb{N},$$
which impples
\begin{eqnarray}
d(\F1(y),\F1(x)) & \leq & \limsup_{n \rightarrow \infty} d(F_n(y),F_n(x)) \nonumber \\
& \leq & \limsup_{n \rightarrow \infty} md(y,x) \nonumber \\
& = & md(y,x). \nonumber 
\end{eqnarray}
For case 2 we know that there must exist a $z \in (y,x)$ such that $\tilde{\Pi}(x) > m$ and therefore there is an 
$n_0 \in \mathbb{N}$ such that 
$$y,x \not\in T_{n_0,i(n_0,z)}^A$$
and indeed 
$$i(n_0,y) < i(n_0,z) < i(n_0,x).$$
It follows that we can find a minimum such $n_0$ and therefore an $n_1 \in \mathbb{N}$ such that 
$$i(n_1,x) - 3 \leq i(n_1,y) \leq i(n_1,x)-2$$
and such that for all $n <n_1$
$$i(n,y) \in \{i(n,x)-1,i(n,x).$$
From this, it follows firstly that for each $n  < n_1$
$$[y,x] \subset T_{n,i(n,y)}^A \cup T_{n,i(n,x)}^A$$
which implies that $F_{n_0}|_{[y,x]}$ has Lipschitz constant 
$$Lip F_{n_0}|_{[y,x]} \leq \max\{\tilde{\Pi}^A_{n_0}(x), \tilde{\Pi}^A_{n_0}(y)\}  \leq m$$
so that 
$$d(F_{n_0}(y),F_{n_0}(x)) \leq md(y,x).$$
It also follows from the choice of $n_1$ that 
$$d(\F1(y),\F1(x)) < 2 \max_{w \in \{x,y\}}\Hm{1}(A_{n_0-1,i(n_0-1,w)}^A).$$
Now, using Lemma $\ref{lem7}$ we know 
$$\pi_x(O_{A_{n_0,i(n_0,y)}^A}(T_{n_0,i(n_0,y)}^A \cap T_{n_0,i(n_0,x)}^A)) 
\cap O_{A_{n_0,i(n_0,y)}^A}(A_{n_0,i(n_0,y)}^A) = \emptyset$$ 
and thus 
\begin{eqnarray}
d(F_{n_0}(y),F_{n_0}(x)) & \geq & \Hm{1}(A_{n_0,i(n_0,y)+1}^A) \nonumber \\
& \geq & {{1}\over{2}}\min_{w\in{y,x}}\Hm{1}(A_{n_0-1,i_w(n_0-1)}^A) \nonumber 
\end{eqnarray}
the latter following since $A_{n_0,i(n_0,y)+1}^A$ is a shorter side of either $A_{n_0-1,i(n_0-1,y)}^A$ or 
$A_{n_0-1,i(n_0-1,x)}^A$.
\newline \newline
Since
$$y \in T_{n_0,i(n_0,y)+1}^A \Rightarrow 1 \leq \tilde{\Pi}^A_{n_0-1}(y) \leq m$$
and 
$$x \in T_{n_0,i(n_0,x)+1}^A \Rightarrow 1 \leq \tilde{\Pi}^A_{n_0-1}(x) \leq m$$
it follows that 
$$\max_{w \in \{x,y\}}\Hm{1}(A_{n_0-1,i(n_0-1,w)}^A)\leq m\min_{w \in \{x,y\}}\Hm{1}(A_{n_0-1,i(n_0-1,w)}^A).$$
Thus
$$d(F_{n_0}(y),F_{n_0}(x)) \geq {{1}\over{2m}}\max_{w\in\{x,y\}}\Hm{1}(A_{n_0-1,i(n_0-1,w)}^A).$$
Hence
\begin{eqnarray}
d(\F1(y),\F1(x)) & \leq 2 & \max_{w\in \{x,y\}}\Hm{1}(A_{n_0-1,i(n_0-1,w)}^A) \nonumber \\
& \leq & 4md(F_{n_0}(y),F_{n_0}(x)) \nonumber \\
& \leq & 4m^2d(y,x). \nonumber 
\end{eqnarray}
Combining the two cases gives us, using $m \geq 1$ 
\begin{eqnarray}
d(\F1(y),\F1(x) ) & \leq & \max\{m,4m^2\}d(x,y) \nonumber \\
& = & 4m^2d(x,y) \nonumber 
\end{eqnarray}
for each $x,y \in \Lambda_m^{-1}$. 
\end{proof}
\end{lem} \noindent
\begin{prop}\label{prop15}\thst 
Let $A \in \cl{K}$ and let $\F1$ be the function related to $A$. Then
\begin{enumerate}
\item $\F1$ is continuous, 
\item should $B \subseteq A_0$ be closed, then $\F1(B) \subseteq A$ is compact, 
\item if $B \subset A_0$ is such that $\Hm{1}(B) > 0$ then $\Hm{1}(\F1(B)) > \Hm{1}(B)/6 > 0$,
\item if $\Hm{1}(\Lambda_{\infty}^{-1})>0$ then $\Hm{1}(\Lambda_{\infty})=\infty$, and
\item if $\Theta^1(x,\Hm{1},\Lambda_{\infty}^{-1})>0$ then $\Theta^1(\F1(x),\Hm{1},\Lambda_{\infty})=\infty$
\item $\F1$ is $\Hm{1}$-measurable.
\end{enumerate}
\pf
As we are considering only one $A$ we shall omit the $A$ superscripts.
\newline \newline
For (1), 
\newline
since for all constructions $A$ that we consider we have $\theta_{0,0} \leq \pi/32$ we see that 
$$diam(T_{n,\cdot}) = \Hm{1}(A_{n,\cdot}) \leq {{(cos(\pi/32))^{-1}}\over{2}}\Hm{1}(A_{n-1,\cdot})$$
which inductively gives us
$$diam(T_{n,i}) \leq \left({{(cos(\pi/32))^{-1}}\over{2}}\right)^n\Hm{1}(A_0).$$
Since $cos(\pi/32) >1/2$, $(cos(\pi/32))^{-1}/2 <1$ so that 
$$\lim_{n \rightarrow \infty}diam(T_{n,\cdot}) = 0.$$
It follows that for all $\e>0$, $diam(T_{n,\cdot})<\e/2$ for all $n$ greater than some sufficiently large $n_0$. 
Consider $x_1,x_2 \in A_0$ such that $|x_1-x_2|<2^{-n_0}$. Then 
$$x_1,x_2 \in \left[{{i-1}\over{2^n_0}},{{i}\over{2^n}}\right] \cup \left[{{i}\over{2^n}},{{i+1}\over{2^n}}\right]$$
for some $i \in \{1,...,2^n-1\}$, so that since $\F1([(i-1)2^{-n},i2^{-n}]) \subset T_{n,i}$ for each $n \in \mathbb{N}$ 
and each $i \in \{1,...,2^n\}$
$$\F1(x_1),\F1(x_2) \in T_{n_0,i} \cup T_{n_0,i+1},$$
which implies 
$$|\F1(x_1)-\F1(x_2)| \leq diam(T_{n_0,i-1})+diam(T_{n_0,i}) < \e.$$
For (2), \newline
Since $A_0$ is bounded, so to is any closed subset of $A_0$, thus should $B$ be a closed subset of $A_0$ it is 
also compact. It then follows from the fact that $\F1$ is continuous that $\F1(B)$ is closed and indeed bounded 
since $\F1(A_0) \subset [0,1]\times [0,1]$ and thus also compact. \newline \newline
For (3) \newline
Let our set, for convenience be denoted $K$. Let $\Hm{1}(K)>0$, say $\Hm{1}(K)=:\beta$. It follows that there 
is a $\delta_0 >0$ such that for all $0<\delta <\delta_0$
$$\Hm{1}_{\delta}(K) > {{\beta}\over{2}}.$$
Now, let $\delta <\delta_0$ and $\{B_{\delta}\}$ be a $\delta$-cover of $\F1(K)$ and consider a $B \in B_{\delta}$.
By Lemma $\ref{lem7}$ we see that there is an $n(B) \in \mathbb{N}$ such that whether or not $center(B) \in \F1(K)$ 
$B \cap \F1(B)$ subset $T^A_{n(B),i(B)-1} \cup T^A_{n(B),i(B)} \cup T^A_{n(B),i(B)+1}$ for some 
$\i(B) \in \{2,...,2^{n(B)}-1\}$with 
$$diam(T^A_{n(B),\cdot})=length(A^A_{n(B),\cdot}) \in \left({{diam(B)}\over{2}},diam(B)\right)$$
so that 
$$diam(B) > \sum_{j=i(B)-1}^{i(B)+1}diam(T^A_{n(B),\cdot}).$$
In this case we also have 
\begin{eqnarray}
\F1(B \cap \F1(K)) & \subset & \bigcup_{j=i(B)-1}^{i(B)+1}\F1^{-1}(T_{n(b),j}) \nonumber \\
& = & \bigcup_{j=1(B)-1}^{i(B)+1}F_{n(b)}^{-1}(A_{n(B),j}) \nonumber 
\end{eqnarray}
which, since $F_{n(B)}$ is an expansion map, gives three intervals $I_{B,j}$, $j=1,2,3$ with 
$$diam(I_{B,j}) = length(I_{B,j}) \leq length(A_{n(B),j}) < diam(B) < \delta.$$
It follows that 
$$\sum_{j=i(B)-1}^{i(B)+1}diam(I_{B,j}) < 3diam(B).$$
Since 
$$\F1(K) \subset \bigcup_{B \in B_{\delta}}(B \cap \F1(K))$$
it follows that  
$$K \subset \bigcup_{B \in B_{\delta}}\bigcup_{j=i(B)-1}^{i(B)+1}I_{B,j}$$
which implies that $\{\{I_{B,j}\}_{B\in B_{\delta}}\}_{j=i(B)-1}^{i(b)+1}$ is a $\delta$ cover of $K$ 
and thus that 
$$\sum_{B \in B_{\delta}}\sum_{j=i(B)-1}^{i(B)+1}I_{B,j} \geq \Hm{1}_{\delta}(K) > {{\beta}\over{2}}$$
and therefore
$$\sum_{B \in B_{\delta}}diam(B) > {{1}\over{3}}\sum_{B \in B_{\delta}}\sum_{j=i(B)-1}^{i(B)+1}I_{B,j} 
>{{1}\over{3}}{{\beta}\over{2}} = {{\beta}\over{6}}.$$
Since this is true for any such $\delta$-cover of $\F1(K)$ we see that 
$$\Hm{1}_{\delta}(\F1(K)) > {{\beta}\over{6}}$$
for any $\delta < \delta_0$ and therefore that 
\begin{eqnarray}
\Hm{1}(\F1(K)) & = & \lim_{\delta \rightarrow 0}\Hm{1}_{\delta}(\F1(K)) \nonumber \\
& \geq & \lim_{\delta \rightarrow 0} {{\beta}\over{6}} \nonumber \\
& = & {{\beta}\over{6}} \nonumber \\
& > & 0.
\end{eqnarray}
\newline \newline
For (4), \newline
Let $M>0$. Then since $\Hm{1}_{A_0}$ is Radon  and 
$$\Lambda_{\infty}^{-1} = \bigcup_{n \in \mathbb{N}}\{x \in \Lambda_{\infty}^{-1}:\tilde{\Pi}_n(x)>M\}$$
it follows that there is an $n_0 \in \mathbb{N}$ with 
$$\Hm{1}(\{x \in \Lambda_{\infty}^{-1}:\tilde{\Pi}_{n_0}(x)>M\}) > {{\Hm{1}(\Lambda_{\infty}^{-1})}\over{2}}>0.$$
We set 
$$\Lambda_{\infty, n_0}^{-1}:=\{x \in \Lambda_{\infty}^{-1}:\tilde{\Pi}_{n_0}(x)>M\}.$$
It follows that with 
$$X:=\{i \in \{1,...,2^{n_0}\}:T_{n_0,i} \cap F_{n_0}(\Lambda_{\infty,n_0}^{-1} \not= \emptyset\}$$
\begin{eqnarray}
\Hm{1}(F_{n_0}(\Lambda_{\infty,n_0}^{-1})) & = & \sum_{i \in X}\Hm{1}(F_{n_0}(\Lambda_{\infty, n_0}^{-1}) \cap T_{n_0,i})
\nonumber \\
& > & M\sum_{i \in X}\Hm{1}(\Lambda_{\infty,n_0}^{-1} \cap [2^{-n_0}(i-1),2^{n_0}i]) \nonumber \\
& = & M \Hm{1}(\Lambda_{\infty,n_0}^{-1}) \nonumber \\
& > & M{{\Hm{1}(\Lambda_{\infty}^{-1})}\over{2}}. \nonumber
\end{eqnarray}
We then apply (3) to each set $A^{n_0,i} \in \cl{K}$ defined as the subconstruction (and subset) of $A$ starting 
with $A^{n_0,i}_0 = A_{n_0,i}$ to find
$$\Hm{1}(\cl{F} \circ F_{n_0}^{-1}(F_{n_0}(\Lambda_{\infty,n_0}^{-1}) \cap T_{n_0,i}))> 
{{\Hm{1}(F_{n_0}(\Lambda_{\infty ,n_0}^{-1}) \cap T_{n_0,i})}\over{6}}$$
and thus that 
\begin{eqnarray}
\Hm{1}(\F1(\Lambda_{\infty,n_0}^{-1})) & = & \sum_{i\in X}\Hm{1}(\F1 \circ F_{n_0}^{-1}(F_{n_0}(\Lambda_{\infty, n_0}^{-1}) 
\cap T_{n_0,i})) \nonumber \\
& > & {{1}\over{6}}\sum_{i \in X}\Hm{1}(F_{n_0}(\Lambda_{\infty, n_0}^{-1}) \cap T_{n_0,i}) \nonumber \\
& = & {{1}\over{6}}\Hm{1}(F_{n_0}(\Lambda_{\infty, n_0}^{-1})). \nonumber
\end{eqnarray}
We therefore now have 
$$\Hm{1}(\F1(\Lambda_{\infty,n_0}^{-1}))>{{1}\over{6}}\Hm{1}(F_{n_0}(\Lambda_{\infty,n_0}^{-1})) > 
{{M\Hm{1}(\Lambda_{\infty}^{-1})}\over{12}}.$$
Since this is true for each $M>0$ it follows that 
$$\Hm{1}(\F1(\Lambda_{\infty}^{-1})) > \Hm{1}(\F1(\Lambda_{\infty,n_0}^{-1})) = \infty.$$
For (5), \newline
Suppose $x \in \Lambda_{\infty}^{-1}$ is such that $\Theta^1(x,\Hm{1},\Lambda_{\infty}^{-1}) > 0$.
\newline \newline
Consider $\F1(x)$ and let $\rho>0$. We know firstly from definition that there is an $n_0>0$ such that 
$F_n(x) \in B_{\rho /2}(\F1(x))$ for all $n > n_0$ and thus, since 
from the proof of (1) $diam(T_{n,i}) \rightarrow 0$ as $n \rightarrow \infty$, there is an $n_1 \geq n_0$ such that 
$diam(T_{n_1,\cdot})<\rho /4$ and thus $\cup_{j=-1}^{1}T_{n_1,i(x,n_1)+j} \subset B_{\rho}(\F1(x))$.
\newline \newline
For the remainder of (5) we write $i:=i(n_1,x)$
We now, temporarily have two cases to consider, namely CASE I that $F_n(x) \in E(T_{n_1,i})$ and CASE II that 
$F_n(x) \in T_{n_1,i} - E(T_{n_1,i})$.
\newline \newline
CASE I: \newline 
In this case $F_{n_1}(x) \in T_{n_1,i} \cap T_{n_1,i-1}$ or $F_{n_1}(x) \in T_{n_1,i} \cap T_{n_1,i+1}$, without loss 
of generality let us suppose that it is the latter case. Then 
\begin{eqnarray}
x & = & i2^{-n_1} \nonumber \\ 
& \in & ((i-1)2^{-n_1},(i+1)2^{-n_1}) \nonumber \\
& \subset & [(i-1)2^{-n_1},(i+1)2^{-n_1}] \nonumber \\
& = & F_n^{-1}(A_{n,i}\cup A_{n,i+1}). \nonumber 
\end{eqnarray}
Since $\Theta^1(x,\Hm{1},\Lambda_{\infty}^{-1}) > 0$ it follows that 
$$\Hm{1}(\Lambda_{\infty}^{-1} \cap [(j-1)2^{-n_1},j2^{-n_1}] >0$$ 
for atleast one $j \in \{i,i+1\}$. Without loss of generality let us assume that $j=i$. Then 
\begin{eqnarray}
\Hm{1}(F_{n_1}^{-1}(\Lambda_{\infty}^{-1}) \cap A_{n_1,i}) & = & \tilde{\Pi}_{n_1}(x)\Hm{1}(\Lambda_{\infty}^{-1} 
\cap [(i-1)2^{-n_1},i2^{-n_1}]) \nonumber \\
& \geq & \Hm{1}(\Lambda_{\infty}^{-1} \cap [(i-1)2^{-n_1},i2^{-n_1}]) \nonumber \\
& > & 0. \nonumber
\end{eqnarray}
CASE II: \newline
In this case $F_{n_1}(x) \in T_{n_1,i} - E(T_{n_1,i})$ so that 
$$x \in ((i-1)2^{-n_1},i2^{-n_1}) \subset [(i-1)2^{-n_1},i2^{-n_1}] = F_{n_1}^{-1}(A_{n,i}).$$
Thus since $\Theta^1(x,\Hm{1},\Lambda_{\infty}^{-1}) > 0$ it follows that
$$\Hm{1}(\Lambda_{\infty}^{-1} \cap [(i-1)2^{-n_1},i2^{-n_1}] >0$$ 
and therefore 
\begin{eqnarray}
\Hm{1}(F_{n_1}^{-1}(\Lambda_{\infty}^{-1}) \cap A_{n_1,i}) & = & \tilde{\Pi}_{n_1}(x)\Hm{1}(\Lambda_{\infty}^{-1} 
\cap [(i-1)2^{-n_1},i2^{-n_1}]) \nonumber \\
& \geq & \Hm{1}(\Lambda_{\infty}^{-1} \cap [(i-1)2^{-n_1},i2^{-n_1}]) \nonumber \\
& > & 0. \nonumber
\end{eqnarray}
That is, in either case there is a $n \in \mathbb{N}$ and $i \in \{1,...,2^n\}$ such that 
$T_{n,i} \subset B_{\rho}(\F1(x))$ and $\Hm{1}(F_{n}^{-1}(\Lambda_{\infty}^{-1}) > 0$.
Applying (iv) to the $A_1\in \cl{K}$ resulting from the subconstruction of $A$ on $T_{n,i}$ it follows that 
$\Hm{1}(\Lambda_{\infty} \cap T_{n,i})=\infty$ and thus that $\Hm{1}(B_{\rho}(\F1(x)) \cap \Lambda_{\infty}) = \infty$.
Since this is true for all $\rho>0$ it follows that 
$$\Theta^1(\F1(x),\Hm{1},\Lambda_{\infty})=\infty,$$
completing the proof of (5).
\newline \newline
Proof of (6): \newline
We note that the open sets of $A$ with respect to $\Hm{1}$ measure are $U \cap A$ for $U$ open in the usual sense in 
$\R^2$. Now consider an open set in $A$, $V:=U \cap A$ for some $U$ open in $\R^2$.
\newline \newline
Let $$\cl{T}_1 := \cup\{T_{1,i}:T_{1,i} \subset U\}$$
and in general
$$\cl{T}_n:= \cup\{T_{n,i}:T_{n,i} \subset U\}.$$
We claim that 
$$V = \bigcup_{n=1}^{\infty}(\cl{T}_n \cap A).$$
Clearly $\cl{T}_n \cap A \subset U \cap A$ 
for all $n \in \mathbb{N}$ and thus 
$$\bigcup_{n=1}^{\infty}(\cl{T}_n \cap A) \subset U \cap A = V.$$
Conversely, let $x \in V$. Then $x \in A$ and there exists $\rho > 0$ such that $B_{\rho}(x) \subset U$. Since we know 
that for any $A \in \cl{K}$, and $x \in A$
$$\lim_{n \rightarrow \infty}diam(T_{n,i(n,x)})=0$$
there exists $n_{\rho} \in \mathbb{N}$ such that $diam(T_{n_{\rho},i(n_{\rho},x)} < \rho/2$. Then 
$$T_{n_{\rho},i(n_{\rho},x)} \subset B_{\rho}(x) \subset U$$
thus 
$$T_{n_{\rho},i(n_{\rho},x)} \subset \cl{T}_{n_{\rho}}$$
and thus $x \in \cl{T}_{n_{\rho}}$.
\newline \newline
Since $x \in A$ we have $x \in \cl{T}_{n_{\rho}} \cap A$ and thus 
$$x \in \bigcup_{n =1}^{\infty}(\cl{T}_{n_{rho}} \cap A).$$
It follows that $$V \subset \bigcup_{n =1}^{\infty}(\cl{T}_{n_{rho}} \cap A).$$
Now, for each $n \in \mathbb{N}$
$$\cl{T}_n \cap A = \bigcup_{i \in I_n}T_{n,i} \cap A$$
for some (possibly empty) index $I_n \subset \{0,1,...,2^n-1\}$. Thus 
$$\F1^{-1}(\cl{T}_n \cap A) = \bigcup_{i \in I_n}D_{n,i}$$
where $D_{n,i}$ is the $i$-th dyadic interval of order $n$. Thus 
\begin{eqnarray}
\F1^{-1}(V) & = & \F1^{-1}\left(\bigcup_{n=1}^{\infty}(\cl{T}_n\cap A)\right) \nonumber \\
& = & \bigcup_{n=1}^{\infty}\F1^{-1}(\cl{T}_n \cap A) \nonumber \\
& = & \bigcup_{n=1}^{\infty}\bigcup_{i \in I_n}D_{n,i} \nonumber 
\end{eqnarray}
which is a Borel set in $A_{0,0}$ and thus $\Hm{1}$-measurable. It follows that for any Borel set $B \in A$ $\F1^{-1}(B)$ 
is a Borel set in $A_{0,0}$. Thus, finally, if $B$ is a $\Hm{1}$-measurable set in $A$, $\F1^{-1}(B)$ is a 
$\Hm{1}$-measurable set in $A_{0,0}$. The fact that the measurability of the inverse images of measurable sets 
follows from the measurability of the inverse images of open sets is standard measure theory and is discussed 
in, for example, Rudin \cite{rudin} or Bartle \cite{bartle}. \hfill 
\end{proof}
\end{prop}\noindent
To complete the preliminary results required for our study of measure and rectifiability of sets in $\cl{K}$ we have one 
more lemma concerning density to consider. It is this final general density Lemma that will be applied in the 
proof of non-rectifiability of those Koch sets which are not rectifiable (which ones they are will be made clear later).
It shows the presence of infinite density almost everywhere in the image of any measurable subset of 
$\Lambda_{\infty}^{-1}$ of positive measure. In order to prove this Lemma, however, we first need a couple of general 
measure theoretic results showing that the set of points density one are sufficiently large in a set of positive measure 
in $A_{0,0}$. The second is a condition of non-rectifiability.
\begin{prop}\label{prop16}\thst
Let $B \subset A_0$ be $\Hm{1}$-measurable, then
$$\Hm{1}(\{x\in B:\Theta^1(x,\Hm{1},B)=1\})=\Hm{1}(B).$$
\pf
Since $B$ is $\Hm{1}$-measurable we know that for all $\rho>0$ 
$$1=(2\rho)^{-1}\Hm{1}(B_{\rho}(x)) = (2\rho)^{-1}(\Hm{1}(B_{\rho}(x) \cap B) + \Hm{1}(B_{\rho}(x) \cap B^c))$$
so that 
\begin{eqnarray}
1& = & \lim_{\rho \rightarrow 0}(2\rho)^{-1}\Hm{1}(B_{\rho}(x)) \nonumber \\
& = & \lim_{\rho \rightarrow 0}(2\rho)^{-1}(\Hm{1}(B_{\rho}(x) \cap B) + \Hm{1}(B_{\rho}(x) \cap B^c)) \nonumber \\
& = & \Theta^1(x,\Hm{1},B) + \Theta^1(x,\Hm{1},B^c). \nonumber
\end{eqnarray}
From standard theory (see for example [Simon3] Theorem 3.5) we know 
$\Theta^1(x,\Hm{1},C) = 0$ for $\Hm{1}$-almost all $x \in C^c$ for any $\Hm{1}$-measurable set $C$ with $\Hm{1}(C)<\infty$.
Hence
$\Theta^1(x,\Hm{1},B^c)=0$ for $\Hm{1}$-almost all $x \in B$ and thus 
$$\Theta^1(x,\Hm{1},B)=1-\Theta^1(x,\Hm{1},B^c)=1$$
for almost all $x \in B$. The result follows 
\end{proof}
\end{prop} \noindent
\begin{prop}\label{prop17}\thst
Let $A \subset \R^2$. Let $\theta$ be a $L^1(\Hm{1},\R^2, \R)$ positive function on $A$. Suppose that $B$ is a subset 
of $A$ of positive measure that satisfies $\theta(x)\geq r>0$ for all $x \in B$. Let $x \in B$ satisfy 
$$\Theta^1(\Hm{1},A,x)\geq \Theta(\Hm{1},B,x)=\infty .$$
Then $A$ does not have a $1$-dimensional approximate tangent plane for $A$ at $x$ with respect to $\theta$.
\pf
Let $P$ be any potential approximate tangent plane for $A$ at $x$ with respect to $\theta$ 
and define $\phi \in C_c^0(\R^2;\R)$ by 
$$\phi(x) := \cases{1 \ \ |x| \leq 1 \cr 2-|x| \ \ 1\leq |x| \leq 2 \cr 0 \ \ \hbox{ otherwise }\cr}.$$
We then have 
$$\int_P\phi d\Hm{1} = 3.$$
However,
\begin{eqnarray}
\lim_{n \rightarrow \infty}\lambda^{-1} \int_A\phi(\lambda^{-1}(z-x))d\Hm{1}(z) & \geq & 
\lim_{n \rightarrow \infty}\lambda^{-1} \int_B\phi(\lambda^{-1}(z-x))d\Hm{1}(z) \nonumber \\
& \geq & \lim_{n \rightarrow \infty}\lambda^{-1} \int_{B \cap B_{\lambda}(x)}rd\Hm{1} \nonumber \\
& = & r \lim_{n \rightarrow \infty}\lambda^{-1} \int_{B \cap B_{\lambda}(x)}1d\Hm{1} \nonumber \\
& = & r \lim_{\lambda \rightarrow 0} {{\Hm{1}(B \cap B_{\lambda}(x))}\over{\lambda}} \nonumber \\
& = & 2r\Theta^1(x,\Hm{1},A) \nonumber \\
& > & 3. \nonumber 
\end{eqnarray}
It is therefore impossible that $A$ have an approximatye tangent plane at $x$ with respect to $\theta$.
\end{proof}
\end{prop}\noindent
\begin{lem}\label{lem16}\thst
Let $A \in \cl{K}$ and $\Hm{1}(B \cap \Lambda_{\infty}^{-1})>0$ for some measurable subset $B \subset A_{0,0}$. 
Then there exists 
$$B_1 \subset B \cap \Lambda_{\infty}^{-1},$$
$$\Hm{1}(B_1) = \Hm{1}(B \cap \Lambda_{\infty}^{-1})$$
such that 
$$\Theta^1(\Hm{1}, \F1(B_1), \F1(x))=\infty$$
for all $x \in B_1$.
\newline \newline
In particular, if $A \in \cl{K}$ and $\Hm{1}(\Lambda_{\infty}^{-1})>0$, then for $\Hm{1}$-a.e. 
$x \in \Lambda_{\infty}^{-1}$
$$\Theta^1(\Hm{1}, A, \F1(x))\geq \Theta^1(\Hm{1},\Lambda_{\infty},\F1(x))=\infty.$$
\pf
We note from Proposition $\ref{prop16}$ that 
$$\Theta^1(\Hm{1}, B \cap \Lambda_{\infty}^{-1},x)=1$$
for $\Hm{1}$-a.e. $x \in B \cap \Lambda_{\infty}^{-1}$. We thus choose 
$$B_1:=\{x \in B \cap \Lambda_{\infty}^{-1}:\Theta^1(\Hm{1}, B \cap \Lambda_{\infty}^{-1},x)=1\},$$
noting that $\Hm{1}(B_1) = \Hm{1}(B \cap \Lambda_{\infty}^{-1})$ as required.
\newline \newline
Choose now $y \in B_1$ arbitrarily. 
\newline \newline
We then note that from the definition of $\Theta^1$ there must exist an $r_0 > 0$ so that for all $r \leq r_0$
$$(2r)^{-1}\Hm{1}(B_r(y) \cap B_1)>7/8.$$
We now claim that for any dyadic interval $D \ni y$ with $|D|:=\Hm{1}(D) < r_0/2$
$$\Hm{1}(D \cap B_1)>3/4|D|.$$
We see this by selecting 
$$\gamma : = \max \{d(y,z):z \in E(D)\}$$
(where $E(D)$ as elsewhere denotes the endpoints of $D$). Then $\gamma < |D|<r_0$ and $D \subset \overline{B_{\gamma}(y)}$. 
Thus
$${{\Hm{1}(B_1^c \cap B_{\gamma}(y))}\over{2\gamma}} = {{1-\Hm{1}(B_1 \cap B_{\gamma}(y))}\over{2\gamma}} < {{1}\over{8}}$$ 
which implies 
$${{\Hm{1}(B_1^c \cap D)}\over{|D|}} \leq {{2\Hm{1}(B_1^c \cap B_{\gamma}(y))}\over{2\gamma}} < {{1}\over{4}}$$
and thus
\begin{eqnarray}
{{\Hm{1}(B_1 \cap D) }\over{|D|}} & = & {{|D|-\Hm{1}(B_1^c\cap D)}\over{|D|}} \nonumber \\
& > & (|D|-|D|/4)|D|^{-1} \nonumber \\
& = & {{3}\over{4}}, \nonumber 
\end{eqnarray}
proving the claim.
\newline \newline
In particular, the claim holds for any dyadic interval $D_m \ni y$ of order $m \geq m_0$ where $m_0$ is chosen 
such that $2^{-m_0}\leq r_0.$ \newline \newline
Then, selecting, independently from one another, $1>p>0$ and $M \in \R$ with 
$$\Hm{1}(A_{m_0,i(y,m_0)})=diam(T_{m_0,i(y,m_0)})>p.$$
We choose $m \geq m_0$ such that $diam(T_{m,i(y,m)}) \in (p/2,2p)$ and $T_{m,i(y,m)} \subset B_p(y).$
Note that $F_m^{-1}(A_{m,i(y,m)}) = A \cap T_{m,i(y,m)}.$ That is, defining 
$\cl{B}_1:=\F1(B_1)$
$$\Hm{1}(B_p(y) \cap \cl{B}_1) \geq \Hm{1}(\cl{B}_1 \cap T_{m,i(y,m)}) = \Hm{1}(\F1(D_m) \cap \cl{B}_1).$$
Since, for all $x \in B_1 \cap D_m$, $\prod_{m=0}^{\infty}(cos\theta_{n,i(y,n)})^{-1} = \infty$ there exists a 
$q_0 \in \mathbb{N}$ such that for 
$$B^q :=\left\{x \in B_1 \cap D_m:\prod_{n=m+1}^p(cos\theta_{n,i(y,n)})^{-1}>M\right\}$$
$$\Hm{1}(B_{q_0}) > \Hm{1}(B_1 \cap D_m)/2.$$
If this were not true then since $B^q \subset B^{q+1}$ for each $q$ it would follow that 
$$\Hm{1}\left(\bigcup_{q=m+1}^{\infty}B^q\right) \leq {{\Hm{1}(B_1 \cap D_m)}\over{2}}$$
and thus there would exist $x \in B_1 \cap D_m$ such that $\prod_{n=m+1}^{\infty}(cos\theta_{n,i(x,n)})^{-1}<M<\infty.$
This contradiction confirms our claim. \newline \newline
We then note 
$$\Hm{1}(B_1) > {{1}\over{2}}\Hm{1}(B_1 \cap D_m)> {{3}\over{8}}|D_m|$$
and that since  for all $x \in D_m$, for all $x \in B_1$ 
$$\prod_{n=0}^{\infty}(cos\theta_{n,i(x,n)})^{-1} \geq \prod_{n=0}^m(cos\theta_{n,i(x,n)})^{-1} > {{p}\over{|D_m|}}.$$
It then follows that for $\tilde{y}:=\F1(y)$ 
\begin{eqnarray}
\Hm{1}(\cl{B}_1 \cap B_p(\tilde{y})) & \geq & \Hm{1}(\cl{B}_1 \cap T_{m,i(y,m)}) \nonumber \\
& \geq & \Hm{1}\left(\cl{B}_1 \cap \bigcup_{D_q \cap B_1 \not= \emptyset}T_{q,i(D_q,q)}\right) \nonumber \\
& = & \sum_{D_q \cap B_1 \not= \emptyset}\Hm{1}(\cl{B}_1 \cap T_{q,i(D_q,q)}) \nonumber \\
& \geq & \sum_{D_q \cap B_1 \not= \emptyset}\Hm{1}(F_q(D_q \cap B_1)) \nonumber \\
& = & \sum_{D_q \cap B_1 \not= \emptyset}\prod_{n=0}^q(cos\theta_{q,i(D_q,q)})^{-1}\Hm{1}(D_q \cap B_1) \nonumber \\
& > & {{p}\over{|D_m|}}\sum_{D_q \cap B_1 \not= \emptyset}\prod_{n=m+1}^q(cos\theta_{q,i(D_q,q)})^{-1}\Hm{1}(D_q \cap B_1) 
\nonumber \\
& > & {{Mp}\over{|D_m|}}\sum_{D_q \cap B_1 \not= \emptyset}\Hm{1}(D_q \cap B_1) \nonumber \\
& = & {{Mp}\over{|D_m|}}\Hm{1}\left(\bigcup_{D_q \cap B_1 \not= \emptyset}D_q \cap B_1\right) \nonumber \\
& = & {{Mp}\over{|D_m|}}\Hm{1}(B_1) \nonumber \\
& > & {{3Mp}\over{8|D_m|}}|D_m| \nonumber \\
& = & {{3Mp}\over{8}} \nonumber
\end{eqnarray}
Since this is true for any $p < diam(T_{m_0,i(y,m_0)}$ 
$$\Theta^1(\Hm{1}, \cl{B}_1,\tilde{y}) = \lim_{p \searrow 0}{{\Hm{1}(\cl{B}_1 \cap B_p(\tilde{y}))}\over{2p}} \geq 
{{3M}\over{16}}.$$
Since this is true for each $M \in \R$ we have 
$$\Theta^1(\Hm{1}, \cl{B}_1,\tilde{y})=\infty.$$
As this is true for any $y \in B_1$ it follows that 
$$\Theta^1(\Hm{1}, \F1(B_1),\F1{y})=\infty$$
for each $y \in B_1$ completing the first part of the proof. 
\newline \newline
For the final part of the proof we note that $A_{0,0}$ is itself measurable and that 
$A_{0,0} \cap \Lambda_{\infty}^{-1}=\Lambda_{\infty}^{-1}$. It follows from the above that there is a set 
$B \subset \Lambda_{\infty}^{-1}$ with $\Hm{1}(B)=\Hm{1}(\Lambda_{\infty}^{-1}$ so that 
$$\Theta^1(\Hm{1},\F1(B),\F1(x))=\infty$$
for all $x \in B$. Since $A_{0,0} \supset \Lambda_{\infty}^{-1} \supset B$, $\F1(A_{0,0}=A$, 
$\F1(\Lambda_{\infty}^{-1}=\Lambda_{\infty}$ and $\Hm{1}$-a.e. $x \in \Lambda_{\infty}^{-1}$, $x \in B$ it follows that 
for all $x \in B$ and thus $\Hm{1}$-a.e. in $\Lambda_{\infty}^{-1}$
\begin{eqnarray}
\Theta^1(\Hm{1},A,\F1(x)) \geq  \Theta^1(\Hm{1},\Lambda_{\infty},\F1(x)) \geq \Theta^1(\Hm{1},\F1(B),\F1(x))=\infty, 
\nonumber 
\end{eqnarray}
which completes the proof.
\end{proof}
\end{lem}\noindent
This completes the preliminary 
results that we need for the rectifiability and measure results on sets in $\cl{K}$. 
\section{Relative Centralisation of Semi-Self-Similar Sets}
We now look at some preliminary 
results that we will need for results on dimension. We will reduce all of our questions to an application of the 
results of Hutchinsion \cite{hutch} to get our dimension results. We do this, in essence, by a comparison principle. 
We show that sets in $\cl{K}$ depending on properties of $\tilde{\theta}^A$ can be dimension invariantly rearranged 
so that they are supersets of some sets to which Hutchinsons results apply and subsets of others. By considering 
sequences of such rearrangements we can deduce the dimension of our sets from the dimensions of the sets to which 
we are comparing.
\newline \newline
It is infact true that we could, in principal, apply Hutchinsons results directly. However, the parameters of the sets 
and "self-similarity" functions cannot be (at least not easily) extracted from sets in $\cl{K}$. Thus actually giving 
an explicit dimension directly is not possible.
\newline \newline
Our comparison principle, or rearrangement involves seperating each triangular cap in a particular approximation 
to some $A \in \cl{K}$, $T_n$ and moving each seperately by an orthogonal transformation in such a way that each 
of the newly positioned triangular caps remain disjoint. We do this by placing each inside of a triangular cap of 
another, larger, $T_n$ from some other $A^{\prime} \in \cl{K}$. Since all Hausdorff measures are translation invariant 
it follows that Hausdorff dimension is also translation invariant and thus the union of the replaced triangular caps is 
the same dimension as the original caps. We can in this way compare the dimesion of $A$ to that of each $T_n^{A^{\prime}}$ 
and thus of $A^{\prime}$. It will be by selecting appropriate $A^{\prime}$ that we will prove our dimension results.
\newline \newline
We start by defining the transformation process, which, due to the placing of one set into parts of another, we call 
centering. That is one set is centered in the bigger one.
\begin{def1}\label{def31}\thst
Let $A_1, A_2 \subset \R^2$. We say that we can center $A_1$ in $A_2$ (or that $A_1$ can be centered in $A_2$) 
written $A_1 \Cent A_2$ if for each $m \in \mathbb{N}$ there exists sets $A_{1m}$ and $A_{2m}$ such that 
$$\bigcap_{m=1}^{\infty}A_{2m}\subset A_2, \ \ A_{2m} \subset A_{2(m-1)} \ \ \hbox{ for all } \ \ m \in \mathbb{N}$$
$$A_1 \subset \bigcap_{m=1}^{\infty}A_{1m}, \ \ A_{1m} \subset A_{1(m-1)} \ \ \hbox{ for all } \ \ m \in \mathbb{N};$$
that for each $m \in \mathbb{N}$ there exists $n_1(m), n_2(m) \in \mathbb{N}$, $n_1(m) \leq n_2(m)$, disjoint sets 
$\{A_{1mj}\}_{j=1}^{n_1(m)}$ and disjoint sets $\{A_{2mj}\}_{j=1}^{n_2(m)}$ such that 
$$\bigcup_{j=1}^{n_2(m)}A_{2mj} \subseteq A_{2m}$$ and 
$$A_{1m}\subseteq \bigcup_{j=1}^{n_1(m)}A_{1mj};$$
that the sets $A_i$, $A_{im}$ and $A_{imj}$ are all $\Hm{a}$-measurable for $i=1,2$ each $a \in \R$ and appropriate 
$m,j \in \mathbb{N}$ and that
there exist orthogonal transformations $\cl{T}_{m,j}^{A_1,A_2}:\R^2 \rightarrow \R^2$ for $j = 1,...,n_1(m)$ 
such that 
$$\cl{T}_{mj}^{A_1,A_2}(A_{1mj}) \subseteq A_{2mj}.$$
If $A_1 \Cent A_2$ we write 
$$C_n^{A_1,A_2}:= \bigcup_{j=1}^{n_1(m)}{\cal{T}}_{m,j}^{A_1,A_2}(A_{1mj}).$$
\end{def1}\noindent
{\bf Remark} \newline
For any $A_1, A_2 \in \cl{K}$ we can set $n_1(m) = n_2(m)=2^m$ and for each $i \in \{1,2\}$ 
$A_{im}:=T_m:=\cup_{j=1}^{2^m}T_{m,j}$ and $A_{imj}=T_{m,j}$. In this case, as we shall see, if 
$\theta^{A_1}_{n,i} \leq \theta^{A_2}_{n,i}$ for each $n$ and $i$, we have, ignoring the negligible set of edge points $E$, 
$A_1 \Cent A_2$. \newline \newline
It would have been a simpler statement of definition to restrict to the case  $A_1, A_2 \in \cl{K}$. However, as we shall 
see we will need to apply the definition where $A_1$ and $A_2$ are subsets of elements of $\cl{K}$ where certain triangular 
caps have been simply removed in the construction of $A_1$ and $A_2$. In any case, to make the definition intuitively 
easier to understand we may always think of each $A_i$ as an element of $\cl{K}$ with triangular caps removed, each 
$A_{im}$  as a union of a subcollection of the $T_{m,j}^{A_i}$ and each $A_{imj}$ as a $T_{m,j}^{A_i}$.
\newline \newline
In the case that $A_1$ and $A_2$ are actually in $\cl{K}$ we can restate the definition as follows 
\begin{def1}\label{def32}{\bf $\cl{K}$ version} \newline
Let $A_1$ and $A_2$ be $\A$ type sets. We say that we can center $A_1$ in $A_2$ (or that $A_1$ can be 
centered in $A_2$) written $A_1 \Cent A_2$ if for each $n \in \mathbb{N}$ and $i \in \{1,..,2^n\}$ there are orthogonal 
transformations ${\cal{T}}_{n,i}^{A_1,A_2}$ such that ${\cal{T}}_{n,i}^{A_1,A_2}(T_{n,i}^{A_1}) \subset T_{n,i}^{A_2}$.
\newline \newline
If $A_1 \Cent A_2$ then we write
$$C_n^{A_1,A_2}:= \bigcup_{i=1}^{2^n}{\cal{T}}_{n,i}^{A_1,A_2}(T_{n,i}^{A_1}).$$
\end{def1}\noindent
Note that due to the fact that they are orthogonal transformations with both disjoint preimages and disjoint images 
we have both 
$$\Hm{\eta}_{\delta}({\cal{T}}_{n,i}^{A_1,A_2}(T_{n,i}^{A_1}))=\Hm{\eta}_{\delta}(T_{n,i}^{A_1})$$ 
for each $i \in \{1,...,2^n\}$ and 
$$\Hm{\eta}_{\delta}\left(\bigcup_{i=1}^{2^n}{\cal{T}}_{n,i}^{A_1,A_2}(T_{n,i}^{A_1})\right)=
\Hm{\eta}_{\delta}\left(\bigcup_{i=1}^{2^n}T_{n,i}^{A_1}\right)$$ 
for each $n \in \mathbb{N}$, for each pair $A_1 \Cent A_2$ and for each non negative $\eta, \delta \in \R$.
\newline \newline
We now look at two properties of centering. The first is more a property of $\A$ type sets that tells a condition 
allowing one $\A$ type set to be centered into another. The second is a more general result showing that the dimension 
comparison works, thus justifying the use of centering.
\begin{prop}\label{prop18}
Let $A_1$ and $A_2$ be $\A$ type sets. Let $\theta_{n,\cdot}^{A_1}$ be denoted by $\theta_{n}^{A_1}$ and 
$\theta_{n,\cdot}^{A_2}$ be denoted by $\theta_{n}^{A_2}$ for each $n \in \mathbb{N}$.
Then, if $T_{0,1}^{A_1} \subseteq T_{0,1}^{A_2}$ and $\theta_n^{A_1} \leq \theta_n^{A_2}$ for each $n \in \mathbb{N}$
then $A_1 \Cent A_2$.
\pf
We know that $T_{0,1}^{A_1} \subseteq T_{0,1}^{A_2}$ so that by denoting the identity transformation by $\iota $ we have 
${\cal{T}}_{0,1}^{A_1,A_2} \equiv \iota$ and thus 
$${\cal{T}}_{0,1}^{A_1,A_2}(T_{0,1}^{A_1}) \subset T_{0,1}^{A_2}.$$
We then continue the proof by induction on $n$. Assume that 
$${\cal{T}}_{n,i}^{A_1,A_2}(T_{n,i}^{A_1}) \subset T_{n,i}^{A_2}$$
for some $n \in \mathbb{N}_0$ and each $i \in \{1,...,2^n\}$. Consider some arbitrarily chosen $j \in \{1,...,2^n\}$
with 
$${\cal{T}}_{n,j}^{A_1,A_2}(T_{n,j}^{A_1}) \subset T_{n,j}^{A_2}$$
and there fore since 
$$diam(T_{n,i}^A) = \Hm{1}(A_{n,i}^A)$$
for each $\A$ type set $A$ it follows that 
$$\Hm{1}(A_{n,j}^{A_1})\leq\Hm{1}(A_{n,j}^{A_2}).$$
Now, $\theta_{n+1}^{A_1} \leq \theta_{n+1}^{A_2}$ by hypothesis and thus also, by Lemma 13 
\begin{eqnarray}
\Hm{1}(A_{n+1,2j+k}^{A_1}) & = & {{1}\over{2}}(cos(\theta_{n+1}^{A_1}))^{-1}\Hm{1}(A_{n,j}^{A_1}) \nonumber \\
& \leq & {{1}\over{2}}(cos(\theta_{n+1}^{A_2}))^{-1}\Hm{1}(A_{n,j}^{A_1}) \nonumber \\
& \leq & {{1}\over{2}}(cos(\theta_{n+1}^{A_2}))^{-1}\Hm{1}(A_{n,j}^{A_2}) \nonumber \\
& = & \Hm{1}(A_{n+1,2j+p}^{A_2}) \nonumber 
\end{eqnarray}
for each $k,p \in \{-1,0\}$. \newline \newline
Combining these, it follows that $T_{n+1,2j+k}^{A_1}$ can be mapped into $T_{n+1,2j+k}^{A_2}$ by placing 
$A_{n+1,2j+k}^{A_1}$ in the center of $A_{n+1,2j+k}^{A_2}$ for $k \in \{-1,0\}$. By defining 
${\cal{T}}_{n+1,2j+k}^{A_1,A_2}$ to be the orthogonal transformation that does this it follows that 
$${\cal{T}}_{n+1,2j+k}^{A_1,A_2}(T_{n+1,2j+k}^{A_1}) \subset T_{n+1,2j+k}^{A_2}$$
for $k \in \{-1,0\}$. Since $j$ was arbitrary we have ${\cal{T}}_{n+1,i}^{A_1,A_2}$ such that 
$${\cal{T}}_{n+1,i}^{A_1,A_2}(T_{n+1,i}^{A_1}) \subset T_{n+1,i}^{A_2}$$
for all $i \in \{1,...,2^{n+1}\}$, which completes the inductive step in $n$. 
\end{proof}
\end{prop} \noindent
We now prove the crucial step for the result we need to get our desired dimension results, saying that if one set 
can be centered in another then the expected result that it has a smaller dimension than the other holds.
\begin{lem}\label{lem16}\thst
$$A_1 \Cent A_2 \Rightarrow dimA_1 \leq dimA_2.$$
\pf
Let $\eta >0$ be such that $\Hm{\eta}(A_2)=0$. \newline \newline
Now, let $m \in \mathbb{N}$ then since $\Hm{\eta}$ is invariant under orthogonal transformations we have 
\begin{eqnarray}
\Hm{\eta}(A_1) & = & \Hm{\eta}(A_1 \cap A_{1m}) \nonumber \\
& = & \Hm{\eta}\left(\bigcup_{j=1}^{n_1(m)}A_1 \cap A_{1mj}\right) \nonumber \\
& = & \Hm{\eta}\left(\bigcup_{j=1}^{n_1(m)}A_{1mj}\right) \nonumber \\
& = & \sum_{j=1}^{n_1(m)}\Hm{\eta}(\cl{T}_{m,j}^{A_1,A_2}(A_{1mj})) \nonumber \\
& \leq & \sum_{j=1}^{n_1(m)}\Hm{\eta}(A_{2mj}) \nonumber \\
& = &\Hm{\eta}\left(\bigcup_{j=1}^{n_1(m)}A_{2mj}\right) \nonumber \\
& \leq & \Hm{\eta}\left(\bigcup_{j=1}^{n_2(m)}A_{2mj}\right) \nonumber \\
& \leq & \Hm{\eta}(A_{2m}). \nonumber 
\end{eqnarray}
We then find similarly for $m+1$
$$\Hm{\eta}(A_1) \leq \Hm{\eta}(A_{2(m+1)})$$
then since $A_{2(m+1)} \subset A_{2m}$ we have 
$$\Hm{\eta}(A_1) \leq \Hm{\eta}(A_{2m} \cap A_{2(m+1)}),$$
by induction it follows that 
\begin{eqnarray}
\Hm{\eta}(A_1) &\leq& \Hm{\eta}\left(\bigcap_{m=1}^{\infty}A_{2m}\right) \nonumber \\
& \leq & \Hm{\eta}(A_2) \nonumber \\
& = & 0. \nonumber
\end{eqnarray}
As this is true for any $\eta \in \mathbb{R}$ for which $\Hm{\eta}(A_2)=0$ we have 
\begin{eqnarray}
dimA_1 & = & \inf\{\eta :\Hm{\eta}(A_1)=0\} \nonumber \\
& \leq & \inf\{\eta : \Hm{\eta}(A_2)=0\} \nonumber \\
& = & dimA_2. \nonumber 
\end{eqnarray}
\end{proof}
\end{lem} \noindent
This completes the presentation of the necessary preliminary results and thus the chapter. In the following chapter 
we look at the theorems proving various results about the actual measure, rectifiability and dimension of $\A$ type 
sets and Koch type sets.

\chapter{Dimension, Rectifiability and Measure of Generalised Koch type Sets}
We now consider the main results for Koch type sets. That is under what conditions do we have finite, or 
weak locally finite measure. Under what conditions are Koch type sets rectifiable, or not rectifiable, and 
under what conditions can we determine the dimension of a set in $\cl{K}$. The results are all determined from 
the constrtuction parameters. All of the relevant parameters can be expressed in terms of the angles $\theta_{n,i}^A$.
In the case of $\A$ type sets we can exactly categorise the sets with respect to the above questions, for the Koch 
type sets it is not possible. The difference being that in the case of Koch type sets we could be generating measure 
from a pre-image set of measure zero in an otherwise well behaved set. The question of whether or not measure can 
indeed be generated remains at this time unanswered, the important point for us, is that it cannot be ruled out. 
\newline \newline
For this reason some of the results will continue to be stated seperately. 
In the general case we find, with respect to rectifiability, that,
$$A\in \cl{K} \hbox{ is countably 1-rectifiable }\Leftrightarrow \Hm{1}(\{x:\tilde{\Pi}^A(x)=\infty\})=0.$$
With respect to measure, we find that for each $A \in \cl{K}$
$$\Hm{1}(A)=\int_{A_0 \sim \Lambda_{\infty}^{-1}}\tilde{\Pi}^Ad\Hm{1} + \Hm{1}(\Lambda_{\infty}).$$
and that $\Hm{1}(\Lambda_{\infty}^{-1})>0 \Rightarrow \Hm{1}(\Lambda_{\infty})=\infty$.
In general we would also expect $\Hm{1}(\Lambda_{\infty}^{-1})=0 \Rightarrow \Hm{1}(\Lambda_{\infty})=0$ 
(that is the nongeneration of measure condition) so that we would then have 
$$\Hm{1}(A)=\int_{A_0}\tilde{\Pi}^Ad\Hm{1}.$$
While in certain cases (e.g. $\Lambda_{\infty}^{-1}$ is countable)
it is certainly true, it may not be true in general. Note that this result 
holds also for $A \in \cl{K}$ with $dimA > 1$, in which case we get the uninformative result $\Hm{1}(A)=\infty$.
\newline \newline
Finally, with respect to dimension we define
$$\gamma^A_1:=\sup\{a:\Hm{1}(\{x:\tilde{\theta}^A_x \geq a\})>0\},$$
$$\gamma^A_2:=\sup_{x \in A_0}\tilde{\theta}^A_x$$
and find 
$$dim\Gamma_{f(\gamma^A_1)} = f_1(\gamma^A_1) \leq dimA \leq f_1(\gamma^A_2) =dim\Gamma_{f(\gamma^A_2)}$$
where
$$f(\gamma):= (1/2)(tan\gamma )$$
and therefore
$$f_1(\gamma)= -{{ln2}\over{ln((1/2)(1+(tan\gamma )^2)^{1/2})}}.$$
Again, we find simplification under the hypothesis that for $B \subset A_0$ $\Hm{1}(B) = 0 \Rightarrow \Hm{1}(\F1(B))=0$
in that we can then state
$$dim A \equiv f_1(\gamma^A_1).$$
It is in the $\A$ type set case that we can ignore the possibilty of generalisation of measure and thus the "nicer" results 
can be stated for these sets.
\section{Lipschitz Representation and Rectifiability}
We start by showing that in some cases 
an $\A$ type set is actually a Lipschitz graph, where $\F1$ would pass as a Lipschitz function.
\begin{lem}\label{lem17}\thst 
Suppose $A \in A^0$ and $\sum_{n=0}^{\infty}\theta_n^A<\infty$. Then for each $l>0$ there is an $n_0 \in \mathbb{N}$ 
such that $A \cap T_{n_0,i}^A$ can be expressed as the graph of a Lipschitz function with Lipschitz constant less than 
or equal to $l$ over $A_{n_0,i}^A$ for each $i \in \{1,...,2^{n_0}\}$.
\pf
Let $n_0$ be such that 
$$\sum_{n=n_0}^{\infty}\theta_n^A < {{tan^{-1}(l)}\over{5}}.$$
Then let $x,y \in A \cap T_{n_0,i}^A$ for some $i \in \{1,...,2^{n_0}\}$ with $x \not= y$. We then know that there 
exists a $n_1 > n_0$  such that for each $n_0 \leq n < n_1$ $x,y \in T_{n,k}^A$ for some $k$ and that 
$x \in T_{n_1,j}^A$ and $y \in T_{n_1,j\pm 1}^A$ for some integer $j$. Without loss of generality let 
$x \in T_{n_1,j}^A$ and $y \in T_{n_1,j+1}^A$. \newline \newline
By choice of $n_0$ we know that 
$$\psi^{A_{n_1,j}^A}_{A_{n_0,i}^A} < {{tan^{-1}(l)}\over{5}}$$
and by Lemma $\ref{lem7}$ 
$$\psi^{T_{n_1,j}^A}_{T_{n_1,j+1}^A} < 2\theta_{n_1,j}^A < 2{{tan^{-1}(l)}\over{5}}$$
so that when writing $X = \{z \in \R^2:z=x+ty, t \in \R\}$
$$\psi^X_{A_{n_1,j}^A}<2\psi^{T_{n_1,j}^A}_{T_{n_1,j+1}^A} < 4{{tan^{-1}(l)}\over{5}}.$$
Thus
$$\psi^X_{A_{n_0,i}^A} < {{tan^{-1}(l)}\over{5}}+4{{tan^{-1}(l)}\over{5}} = tan^{-1}(l)$$
and hence
$${{|\pi_{(A_{n_0,i}^A)^{\perp}}(x) - \pi_{(A_{n_0,i}^A)^{\perp}}(y)|}\over{|\pi_{A_{n_0,i}^A}(x)-\pi_{A_{n_0,i}^A}(y)|}} 
<tan(tan^{-1}(l)) = l.$$
Noting that $(x,y)$ was an arbitrarily chosen pair of distinct points completes the proof. 
\end{proof}
\end{lem}\noindent
Combining this lipschitz result with Lemma $\ref{lem15}$ we are now able to present the rectifiability results. We 
first prove, both by Lipschitz graphs and the existence of approximate tangent spaces, the rectifiability under 
particular conditions of $\A$ type sets. We present concurring with the philosophy that multiple proof methods 
allow further insight and understanding of the objects involved and are in any case interesting in their own right, as 
well as for comparative purposes.
\newline \newline
We first prove the rectifiability using the Lipschitz lemmmas to show that certain $\A$ type sets can then be expressed as 
$\Hm{1}$-almost everywhere subsets of a countable union of Lipschitz graphs.
\begin{thm}\label{thm6}\thst
Whenever $A \in A^0$ satisfies $\sum_{n=0}^{\infty}\theta_n^A<\infty$, $A$ is countably $1$ rectifiable.
\pf
Since $\sum_{n=0}^{\infty}\theta_n^A<\infty$ there is, by Lemma $\ref{lem17}$, an $n_0 \in\mathbb{N}$ such that for each 
$i \in \{1,...,2^{n_0}\}$ $A \cap T_{n_0,i}^A$ can be expressed as the graph of a Lipschitz graph over $A_{n_0,i}^A$. 
That is there is a Lipschitz function $f_i:\R^2 \rightarrow \R^2$ such that
$$A \cap T_{n_0,i}^A \subset f_i(A_{n_0,i}^A).$$
Then when $S_{n_0,i}^A:\R \rightarrow \R^2$ is a transformation satisfying
$$S_{n_0,i}^A([0,\Hm{1}(A_{n_0,i}^A)])=A_{n_0,i}^A$$
we can define $F_i:\R \rightarrow \R^2$ as $F_i = f_i \circ S_{n_o,i}^A$ to write
\begin{eqnarray}
A & = & \bigcup_{i=1}^{2^{n_0}}A \cap T_{n_0,i}^A \nonumber \\
& \subseteq & \bigcup_{i=1}^{2^{n_0}}f_i(A_{n_0,i}^A) \nonumber \\
& \subset & \bigcup_{i=1}^{2^{n_0}}F_i(\R). \nonumber
\end{eqnarray}
Since this is a subset of a form of expression of a set that is defined as being countably $1$-rectifiable, the 
proof is complete. 
\end{proof}
\end{thm}\noindent
The second proof applies to sets with converging sums of base angles. In this case "potential" approximate tangent 
spaces eventually stop rotating and we can then use the approximate $j$-dimensionality to say that the set will be 
arbitrarily close to the limit of the rotating bases of the triangular caps containing a point and will thus have an 
approximate tangent space there.
\begin{thm}\label{thm7}\thst
Any $A \in A^0$ satisfying $\sum_{n=0}^{\infty}\theta_n^A<\infty$ has an approximate tangent space 
with multiplicity one almost eveywhere and is thus countably $1$-rectifiable.
\pf
We first prove that $A-E$ is countably $1$-rectifiable. Let $y \in A-E$, write $H:=\Hm{1}(A)$ 
and let $f \in C_C^0(\R^2)$. It follows in particular
that $f$ is Lipschitz with Lipschitz constant $F_1$ and that there is an $M$ such that
$$sptf \subset B_{M}(0).$$
Let $F = max\{1,F_1\}$. Since the other case is trivial we assume $M >0$.\newline \newline
Let $\e >0$ and define $\delta = \e/(MF)$.
Since $A \in A^0$ we know that $A-E$ satisfies property (iv), we therefore know that there is a $\rho_y>0$ such that 
for all $\rho \in (0,\rho_y]$ there is a $L_{y,\rho}$ such that 
$$A \cap B_{\rho}(y) \subset L_{\rho,y}^{\delta \rho/2}$$ and we know in fact from the proof that $A-E$ satisfies 
(iv) that we may take $L_{y,\rho}||A_{n_{\rho},i(y,n_{\rho})}^A$ where $A_{n_{\rho},i(y,n_{\rho})}^A$ is taken 
such that $\Hm{1}(A_{n_{\rho},i(y,n_{\rho})}^A) \in (\rho/2,\rho]$ and $y \in T_{n_{\rho},i(y,n_{\rho})}^A$.
\newline \newline
Since $\sum_{n=0}^{\infty}\theta_n^A <\infty$ we know that 
$\{\psi_{\R}^{A_{n,i(n,y)}^A}\}$ is a convergent sequence
and thus there is an affine space $L$ such that 
$$\psi^L_{\R} = \lim_{n \rightarrow \infty}\psi^{A_{n,i(n,y)}^A}_{\R}.$$
We then choose $\rho_1$ such that $\rho_1 < \rho_y$, so that for all $ \rho < \rho_1$ the 
$A_{n_{\rho},i(y,n_{\rho})}^A$ taken as described above is such that 
\begin{equation}
tan^{-1}(\psi^L_{A_{n_{\rho},i(y,n_{\rho})}^A}) < {{\delta}\over{2}} \label{e:ats1}
\end{equation}
with $n_{\rho}$ large enough for Lemma $\ref{lem17}$ 
to gaurantee that $A \cap T_{n_{\rho},i(y,n_{\rho})}^A$ can be expressed 
as the graph of a Lipschitz function with Lipschitz constant $\delta$, 
and since $\sum_{n=0}^{\infty}\theta_n^A<\infty \Rightarrow \prod_{n=0}^{\infty}(cos\theta_n^A)^{-1}<\infty$
we take $\rho_1$ such that $n_{\rho_1}$is such that $\prod_{n=n_{\rho_1}}^{\infty}(cos\theta_n^A)^{-1}<1+\e$.
\newline \newline
Now let $\lambda < {{\rho_1}\over{M}}$. Then we have that 
$A \cap B_{\lambda M}(y) \subset (A_{n_{\lambda},i(y,n_{\lambda})}^A)^{\delta \lambda M/2}$ 
so that by ($\ref{e:ats1}$)
$$tan(\psi^{A_{n_{\lambda},i(y,n_{\lambda})}^A}_L) < {{\delta}\over{2}}$$
so that 
$$A \cap B_{M\lambda}(y) \subset L^{M\delta \lambda}$$
and thus 
$$\eta_{y,\lambda}A \cap B_M(0) \subset L-y)^{M\delta}.$$
On this set we also have 
$$|f(x) - f(\pi_L(x))| < Lipf\cdot \delta M \leq {{MF\e}\over{MF}} = \e $$
for all $x \in \eta_{y,\lambda}(A-E)$.
\newline \newline
By otherwise considering the positive and negative parts of $f$ we may assume that $f \geq 0$.
We then note
$$\int_{\eta_{y,\lambda}(A-E)}f(y)d\Hm{1}(y) \leq \int_{\eta_{y,\lambda}(A-E)}\e d\Hm{1} + 
\int_{\eta_{y,\lambda}(A-E)}f(\pi_L(y))d\Hm{1}(y).$$
Then by Lemma $\ref{lem17}$ and Lemma $\ref{lem8}$ 
we know that we can apply the area formula with Jacobian calculated by taking 
the maximal vertical variation per unit along $L$ as $\delta$ plus $2(2\theta_{n_{\lambda}}^A)$. That is, with 
the Jacobian factor bounded above by $(1+9\delta^2)^{1/2}$ so that we have 
\begin{eqnarray}
\int_{\eta_{y,\lambda}(A-E)}f(y)d\Hm{1}(y) & \leq & \int_{\eta_{y,\lambda}(A-E)}\e d\Hm{1} + 
(1+9\delta^2)^{1/2}\int_Lf(y)d\Hm{1}(y) \nonumber \\
& < &\e\Hm{1}(\eta_{n,\lambda}(A-E))+ (1+9\e ) \int_Lf(y)d\Hm{1}(y) \nonumber 
\end{eqnarray}
which implies
\begin{eqnarray}
\left|\int_{\eta_{y,\lambda}(A-E)}fd\Hm{1} - \int_Lfd\Hm{1}\right| & \leq & \e(1+\e)2M + (1+9\e -1)
\left|\int_Lfd\Hm{1}\right| \nonumber \\
& = & \e(1+\e)2M + (9\e )\int_Lfd\Hm{1}. \nonumber
\end{eqnarray}
Since this is true for all $\e >0$ it follows that
$$\lim_{\lambda \rightarrow 0}\left|\int_{\eta_{y,\lambda}(A-E)}fd\Hm{1} - \int_L fd\Hm{1}\right|=0$$
so that 
$$\lim_{\lambda \rightarrow 0}\int_{\eta_{y,\lambda}(A-E)}fd\Hm{1} = \int_L fd\Hm{1}.$$
That is there is an approximate tangent space for $y$. Since this is true for all $y \in A-E$ and $\Hm{1}(E)=0$ we have 
\begin{eqnarray}
\lim_{\lambda \rightarrow 0}\int_{\eta_{y,\lambda}A}fd\Hm{1} & = &\lim_{\lambda \rightarrow 0}\int_{\eta_{y,\lambda}(A-E)} 
fd\Hm{1} \nonumber \\
& = & \int_L fd\Hm{1} \nonumber 
\end{eqnarray}
for all $y \in A-E$. That is, $A$ has an approximate tangent space for all $y \in A-E$, and therefore $\Hm{1}$-almost 
everywhere which implies that $A$ is countably $1$-rectifiable. 
\end{proof}
\end{thm}\noindent
Although these results are not for the entirety of $\A$ type sets, the completion of the proofs of rectifiability 
falls under the proof for general $\cl{K}$ sets. We thus prove the more general result, stating the cleaner result 
for $\A$ type sets as a Corollary.
\begin{thm}\label{thm8}\thst
Let $A \in \cl{K}$.
\newline
If $\Hm{1}(\Lambda_{\infty})=0$ then $A$ is countably $1$-rectifiable. \newline \newline 
{\bf Remark:} \newline
It would clearly be desireable to be able to show that 
$$\Hm{1}(\Lambda_{\infty}^{-1})=0 \Rightarrow \Hm{1}(\Lambda_{\infty})$$
which would be an a better situation since we have better understanding, perception and control of sets in $A_0$ than 
sets in $A$. It is however not necessarily in general true (though it may be). We do in some limited 
cases have control from $A_0$. For example if $\Lambda_{\infty}^{-1}$ is countable then $\Hm{1}(\Lambda_{\infty})=0$ 
and so the above Theorem would then state that with such a $\Lambda_{\infty}^{-1}$, $A$ is countably $1$-rectifiable.
\pf
We note that 
\begin{eqnarray}
A & = & \Lambda_{\infty} \cup \bigcup_{m=1}^{\infty}\Lambda_m \nonumber \\
& = &  \Lambda_{\infty} \cup \bigcup_{m=1}^{\infty}\F1(\Lambda_m^{-1}) \nonumber \\
& = & \Lambda_{\infty} \cup \bigcup_{m=1}^{\infty}\F1|_{\Lambda_m^{-1}}(\Lambda_m^{-1}). \nonumber 
\end{eqnarray}
Since from Lemma $\ref{lem15}$ we know that $\F1_{\Lambda_m^{-1}}$ is Lipschitz for each $m \in \mathbb{N}$ it follows that 
$A$ is countably $1$-rectifiable should $\Hm{1}(\Lambda_{\infty})=0$.
\end{proof}
\end{thm}\noindent
Before stating the corollary of rectifiability for $\A$ sets, we prove the non-rectifiability result. In this way 
we will be able to demonstrate necessary and sufficient, that is, an equivalence of conditions for sets in $\A$ 
to countable $1$-rectifiability.
\begin{thm}\label{thm9}\thst
Let $A \in \cl{K}$ and $\Hm{1}(\Lambda_{\infty}^{-1})>0$. Then $A$ is not countably $1$-rectifiable.
\pf
Let $\theta$ be any potential multiplicity function for $A$. Then $\theta \in L^1(\Hm{1},A, \R)$ and thus $\theta$ is 
$\Hm{1}$-measurable.
\newline \newline
We then claim that there is an $r > 0$ such that 
$$\Hm{1}(\F1^{-1}(\{x \in A:\theta(x)>r\}) \cap \Lambda_{\infty}^{-1})>0.$$
This is true for otherwise 
$$\Hm{1}(\{x \in A_{0,0}:\theta \circ \F1(x)=0\})>0$$ 
and thus 
$$\Hm{1}(\{x \in A:\theta(x)=0\})>0$$
contradicting $\theta$ being a positive function on $A$. Set 
$$B : = \F1^{-1}(\{x \in A:\theta(x)>r\}).$$
Since $\theta$ is measurable, $\{x \in A:\theta(x)>r\}$ is measurable and thus, since from Proposition $\ref{prop15}$ 
we know $\F1$ is 
measurable, $B$ is $\Hm{1}$-measurable in $A_{0,0}$.
\newline \newline
It then follows from Lemma $\ref{lem16}$ that there exists a $B_1 \subset B$ with $\Hm{1}(B_1)=\Hm{1}(B)>0$ such that 
$$\Theta^1(\Hm{1},\F1(B_1),\F1(x))=\infty$$
for each $x \in B_1$.
\newline \newline
Consider now $f \in C_C^0(\R^2,\R)$ such that $\chi_{B_1(0)} \leq f \leq \chi_{B_2(0)}$ where $\chi$ is the 
characteristic function. Then for any tangent space, $P$, to $A$ that may exist with respect to $\theta$ at $\F1(x)$ for 
some $x \in B_1$
$$\theta(\F1(x))\int_Pf(y)d\Hm{1}(y) \leq \theta(\F1(x))\int_P\chi_{B_2(0)}d\Hm{1}(y) = 2\theta(\F1(x))< \infty.$$
However
\begin{eqnarray}
\lim_{\lambda\searrow 0}\int_{\eta_{x,\lambda}A}f(y)\theta(x+\lambda y)d\Hm{1}(y) & \geq & 
\lim_{\lambda\searrow 0}\int_{\eta_{x,\lambda}\F1(B_1)}f(y)\theta(x+\lambda y)d\Hm{1}(y)\nonumber \\
& \geq &\lim_{\lambda\searrow 0}\int_{\eta_{x,\lambda}\F1(B_1)}\chi_{B_1(0)}\theta(x+\lambda y)d\Hm{1}(y)\nonumber \\
& > & r \lim_{\lambda\searrow 0}\int_{\eta_{x,\lambda}\F1(B_1)}\chi_{B_1(0)}d\Hm{1}(y)\nonumber \\
& \geq & r \Theta^1(\Hm{1},\F1(B_1), x) \nonumber \\
& = & \infty. \nonumber
\end{eqnarray}
Thus
$$\lim_{\lambda\searrow 0}\int_{\eta_{x,\lambda}A}f(y)\theta(x+\lambda y)d\Hm{1}(y) \not=
\theta(\F1(x))\int_Pf(y)d\Hm{1}(y).$$
Since this is true for any $x \in \F1(B_1)$ and $\Hm{1}(\F1(B_1)) \geq \Hm{1}(B_1)>0$ it follows that $A$ does not 
have an approximate tangent space with respect to $\theta$ at $x$ on a set of $x$ of positive measure. 
\newline \newline
Since this holds for any allowed selection of $\theta$ it follows from the definition of rectifiable sets and 
Theorem $\ref{thm1}$ that $A$ is not countably $1$-rectifiable. 
\end{proof}
\end{thm}\noindent
We can now state the cleaner result for $\A$ type sets from which the particular results for $\A$ and $\Ae$ follow.
\begin{cor}\label{cor9}\thst
For an $\A$ type set $A$, $A$ is countably $1$-rectifiable if and only if 
$$\Hm{1}((\Lambda_{\infty}^{-1})^A)=0.$$
\pf
We note that $A$ being $\A$ type set implies $A \in \cl{K}$. Thus from Theorem $\ref{thm9}$, if 
$\Hm{1}((\Lambda_{\infty}^{-1})^A)>0$ then $A$ is not countably $1$-rectifiable. 
\newline \newline
Conversely, Should 
$\Hm{1}((\Lambda_{\infty}^{-1})^A)=0$ then there must exist at least one point, $x$, 
for which $\tilde{\Pi}^A_x\not= \infty$. Since $\tilde{\Pi}^A_x$ is constant for all $x \in A$ for an $\A$ type set 
it follows that $\tilde{\Pi}^A_y \not=\infty$ for each $y \in A_{0,0}$ and thus for each $y \in A$. It follows that 
$\Lambda_{\infty}^A=\emptyset$ and therefore that $\Hm{1}(\Lambda_{\infty}^A)=0$. It thus follows from 
Theorem~$\ref{thm8}$ that $A$ is countably $1$-rectifiable.
\end{proof}
\end{cor}\noindent
\begin{thm}\label{thm10}\thst
Let $\e > 0$ and $A$ be constructed as in Construction $\ref{cons2}$ with this $\e$. Then 
$$\tilde{\Pi}^A_x\equiv \infty$$
and thus $A$ is not $1$-countably $1$-rectifiable.
\pf
From Lemma $\ref{lem14}$ we know that for any $\A$ type set $A_1$,
$$\Hm{1}(\tilde{A_n}^{A_1}) = \Hm{1}(A_{0,0}^{A_1})\prod_{j=0}^n(cos\theta_{j,\cdot}^{A_1})^{-1} = 
\prod_{j=0}^n(cos\theta_{j,\cdot}^{A_1})^{-1}.$$
Since from Lemma $\ref{lem5}$
$$\Hm{1}(\tilde{A_n}^{A})=(1+n16\e^2)^{1/2}$$
it follows that 
\begin{eqnarray}
\tilde{\Pi}^A & = & \lim_{n \rightarrow \infty}\prod_{j=0}^n(cos\theta_{j,\cdot}^{A})^{-1} \nonumber \\
& = &\lim_{n \rightarrow \infty} \Hm{1}(\tilde{A_n}^{A}) \nonumber \\
& = &\lim_{n \rightarrow \infty} (1+n16\e^2)^{1/2} \nonumber \\
& = & \infty.\nonumber 
\end{eqnarray}
Thus $x \in (\Lambda_{\infty}^{-1})^A$ for each $x \in A_{0,0}$. 
This completes the first part of the proof. 
\newline \newline
It thus follows that $\Hm{1}((\Lambda_{\infty}^{-1})^A)>0$.
From Proposition $\ref{prop15}$ (3) it then follows that $\Hm{1}(\Lambda_{\infty}^A)>0$. Therefore,  
from Corollary $\ref{cor9}$, $A$ is not countably $1$-rectifiable.
\end{proof}
\end{thm}\noindent
The proof then that $\Ae$ is not countably rectifiable that we present is an indirect proof, assuming that $\Ae$ 
is countably $1$-rectifiable, which then implies that $\A$ is countably $1$-rectifiable. This contradiction completes 
the proof and the rectifiability results.
\begin{thm}\label{thm11}\thst
For any appropriate $\e > 0$ for $\Ae$ to be defined, $\Ae$ is not countably $1$-rectifiable.
\pf
We prove the Theorem by contradiction. So, suppose that $\Ae$ is countably $1$-rectifiable and so can be written in 
the form
$$\Ae \subset A_0 \cup \bigcup_{n=1}^{\infty}F_n(\R)$$
where $\Hm{1}(A_0) = 0$ and $F_n:\R \rightarrow \R^2$ is a Lipschitz function for each $n \in \mathbb{N}$.
\newline \newline
We now consider that by the construction of $\A$ we know that $\A \cap T_{i,j}$ is $A_{2^{1-i}\e}$ constructed on a base 
of length $\Hm{1}(A_{i,\cdot})$ (which we note importantly is greater than $2^{1-i}$ so that should $\A$ be well defined, 
then so too is the new $\A$). 
\newline \newline
It thus follows that by contradicting $\A$ by $2^{1-i}$ in the  vertical direction and by $\Hm{1}(A_{i,\cdot})$ in the 
horizontal direction we have that the result $C(\A)$ is a copy of any $\A \cap T_{i,j}$ (where $C$ is the contraction 
map satisfying the said conditions). 
\newline \newline
We thus know that there exists contraction maps for each $i \in \mathbb{N}$ and $j \in \{1,...,2^i\}$, 
$O_{ij}:\R^2 \rightarrow \R^2$, such that  
$$O_{ij}(\A) = \A \cap T_{i,j}$$
which implies 
$$O_{ij}(\Ae) \subset \A \cap T_{i,j}$$
and also that 
$$O_{ij}(E) = E \cap T_{i,j}.$$
Define 
$$M_{\A}:= \A \cap \bigcup_{i \in \mathbb{N}} \bigcup_{j=1}^{2^i}O_{ij}(\Ae),$$
and 
$$R_{\A}:= \A \sim \left(\Ae \cup \bigcup_{i \in \mathbb{N}} \bigcup_{j=1}^{2^i}O_{ij}(\A)\right) = \A \sim M_{\A}.$$
It follows that 
$$L_{ijn}:= O_{ij}(F(x)) \ \ \hbox{ } \ \ i,n \in \mathbb{N},j \in \{1,...,2^i\}$$
are Lipschitz functions $L_{ijn}:\R \rightarrow \R^2$. We note that $\{\{L_{ijn}\}_{i,n \in \mathbb{N}}\}_{j=1}^{2^i}$
is countable. Also that $R_{\A}$ ia a subset of the union of balls (or deformed balls) around points in $E$. Also 
that by taking the further addition to $\Ae$, $O_{ij}(\Ae)$, we infinitely reduce this area by continually refining the 
deformed ball around each $e_n$, that is 
$$R_{\A} \subset \bigcup_{n \in \mathbb{N}} \bigcap_{\{i,j:O_{ij}((0,0))=e_n\}}O_{ij}(B_{r_1}((0,0))).$$
With this set up we can then attack the proof.
\newline \newline
We first note that 
\begin{eqnarray}
M_{\A} & = & \Ae \cup \bigcup_{i=1}^{\infty}\bigcup_{j=1}^{2^i}O_{ij}(\Ae) \nonumber \\
& \subset & A_0 \cup \bigcup_{n=1}^{\infty}F_n(\R) \cup \bigcup_{i=1}^{\infty}\bigcup_{j=1}^{2^i}
O_{ij}\left(A_0 \cup \bigcup_{n=1}^{\infty}F_n(\R)\right) \nonumber \\
& = & A_0 \cup \bigcup_{n=1}^{\infty}F_n(\R) \cup \bigcup_{i=1}^{\infty}\bigcup_{j=1}^{2^i}O_{ij}(A_0) \cup 
\bigcup_{i=1}^{\infty}\bigcup_{j=1}^{2^i}\bigcup_{n=1}^{\infty}L_{ijn}(\R) \nonumber \\ 
& = & A_0 \cup \bigcup_{i=1}^{\infty}\bigcup_{j=1}^{2^i}O_{ij}(A_0) \cup \bigcup_{n=1}^{\infty}F_n(\R) \cup 
\bigcup_{i=1}^{\infty}\bigcup_{j=1}^{2^i}\bigcup_{n=1}^{\infty}L_{ijn}(\R), \nonumber 
\end{eqnarray}
where 
\begin{eqnarray}
\Hm{1}\left(A_0 \cup \bigcup_{i=1}^{\infty}\bigcup_{j=1}^{2^i}O_{ij}(A_0)\right) & \leq & \Hm{1}(A_0) 
+ \sum_{i=1}^{\infty}\sum_{j=1}^{2^i}\Hm{1}(O_{ij}(A_0)) \nonumber \\
& \leq & \Hm{1}(A_0)+\sum_{i=1}^{\infty}\sum_{j=1}^{2^i}\Hm{1}(A_0) \nonumber \\
& = & 0 + \sum_{i=1}^{\infty}\sum_{j=1}^{2^i}0 \nonumber \\
& = & 0. \nonumber 
\end{eqnarray}
and $\bigcup_{n=1}^{\infty}F_n(\R) \cup \bigcup_{i=1}^{\infty}\bigcup_{j=1}^{2^i}\bigcup_{n=1}^{\infty}L_{ijn}(\R)$ is a 
countable collection of Lipschitz images. 
\newline \newline
It thus follows that $M_{\A}$ is a countably $1$-rectifiable. That is 
$$M_{\A} = M_0 \cup \bigcup_{n=1}^{\infty}M_n(\R)$$ 
where
$$M_0 = A_0 \cup \bigcup_{i=1}^{\infty}\bigcup_{j=1}^{2^i}O_{ij}(A_0)$$ 
is a set of measure zero and $\{M_n\}_{n=1}^{\infty}$ is a reordering of 
$\{F_n\}_{n=1}^{\infty} \cup \{\{L_{ijn}\}_{i,n=1}^{\infty}\}_{j=1}^{2^i}$.
\newline \newline
We now show that $\Hm{1}(R_{\A})=0$.
\newline \newline
Let $\eta > 0$. For each $i,n \in \mathbb{N}$ there exists $j_n = j_n(i,n) \in \{1,...,2^i\}$ such that 
$O_{ij}((0,0)) = e_n$. That is, $O_{ij}(B_{r_1}((0,0)))$ covers the part of $R_{\A}$ centered on $e_n$, so that since 
$$\lim_{i \rightarrow \infty}\Hm{1}(A_{i,\cdot}) = 0$$
for each $n$ we can choose an $i_n \in \mathbb{N}$ such that $diam(O_{ij}(B_{r_1}((0,0)))) < \eta 2^{-n}.$ Then, since 
$$R_{\A} \subset \bigcup_{n=1}^{\infty}O_{i_nj_n}(B_{r_1}((0,0)))$$
and since $diam(O_{ij}(B_{r_1}((0,0)))) < \eta 2^{-n} < \eta$ for each $n \in \mathbb{N}$ we then have that 
$\{O_{i_nj_n}(B_{r_1}((0,0)))\}_{n=1}^{\infty}$ is an appropriate covering set to estimate $\Hm{1}_{\eta}$ and in fact 
we have 
\begin{eqnarray}
\Hm{1}_{\eta}(R_{\A}) & \leq & \Hm{1}_{\eta}\left(\bigcup_{n=1}^{\infty}O_{i_nj_n}(B_{r_1}((0,0)))\right) \nonumber \\
& \leq & \sum_{n=1}^{\infty}diam(O_{i_nj_n}(B_{r_1}((0,0)))) \nonumber \\
& < & \sum_{n=1}^{\infty}\eta 2^{-n} \nonumber \\
& = & \eta. \nonumber
\end{eqnarray}
Thus 
\begin{eqnarray}
\Hm{1}(R_{\A}) & = & \lim_{\eta \rightarrow 0}\Hm{1}_{\eta}(R_{\A}) \nonumber \\
& < & \lim_{\eta \rightarrow 0} \eta \nonumber \\
& = & 0. \nonumber
\end{eqnarray}
now since $\Ae = M_{\A} \cup R_{\A}$ we have 
$$\Ae = R_{\A} \cup M_0 \cup \bigcup_{n=1}^{\infty}M_n(\R).$$
Since $\Hm{1}(R_{\A})=0$, $$\Hm{1}(R_{\A} \cup M_0) = 0$$
and it follows that $\A$ is countably $1$-rectifiable. This contradicts Theorem~$\ref{thm10}$, thus $\Ae$ is not 
countably $1$-rectifiable. 
\end{proof}
\end{thm}\noindent
This completes our study of rectifiability, we move on to the measure results before finally considering the dimension 
of Koch type sets.
\section{Measure Formulae for Koch Type Sets}
For our measure result we present, as previously seen, a formula that resembles the Area Formula. We could also have 
applied the Area Formula (for more information on the Area Formula see for example Simon \cite{simon1}) but not without 
some difficulty. We therefore present a self contained direct proof of the result.
\begin{thm}\label{thm12}\thst
Let $A \in \cl{K}$. Then, for all measurable $B \subset A_{0,0}$ the following holds
$$\Hm{1}(\F1(B)) = \int_{B \sim \Lambda_{\infty}^{-1}} \tilde{\Pi}d\Hm{1} + \Hm{1}(\F1(B) \cap \Lambda_{\infty}).$$
\newline \newline
{\bf Remark:} \newline
As with the rectifiability theorem, the statement of this theorem would be simplified should it be true that  
$$\Hm{1}(\Lambda_{\infty}^{-1})=0 \Rightarrow \Hm{1}(\Lambda_{\infty})=0$$
in which case we could write 
$$\Hm{1}(\F1(B)) = \int_B\tilde{\Pi}d\Hm{1},$$
since, should $\Hm{1}(\Lambda_{\infty}^{-1})>0$, both sides would then be $\infty$ so that they could in this case 
also be reconciled with one another.
\newline \newline
It seems as though an application of the area formula for rectifiable sets is all that is necessary, which is likely to be 
true, however, since the convergence of $\tilde{\Pi}_n(x)$ is equivalent to the convergence of $\sum_n\theta_{n,i}^A(x)^2$ 
and thus not necessarily of $\sum_n\theta_{n,i}^A(x)$, the Jacobian is by no means a trivial quantity to calculate or show 
that it is equal to $\tilde{\Pi}$ on $A_{0,0} \sim \Lambda_{\infty}^{-1}$. 
\pf
We note that for any measurable $C \subset A_{0,0}$ 
$$\F1(D) = \bigcap_{n=1}^{\infty}\bigcup_{i \in X_n}T_{n,i}^A$$
where 
$$X_n:=\{i \in \{1,...,2^n\}:i=i(n,x) \hbox{ for some } x \in D\}$$
and so can be constructed from countable unions and intersections of $\Hm{1}$-measurable subsets of $\R^2$ and is therefore 
measurable. Also, since from Lemma $\ref{lem15}$ 
$F_n$ is Lipschitz for each $n \in \mathbb{N}$ these sets are also measurable. 
\newline \newline
Further, since $F_n$ is a Lipschitz map for each $n \in \mathbb{N}$, if $D \subset A_{0,0}$ so to is $F_n(D)$ for 
each $n \in \mathbb{N}$.
\newline \newline
It follows then that 
$$\Hm{1}(\F1(B)) = \Hm{1}(\F1(B) \cap \Lambda_{\infty})+\Hm{1}(\F1(B) \sim \Lambda_{\infty}).$$
We consider the second term. 
\newline \newline
Let $q \in \mathbb{N}$ and define
$$H_{nq}:=\left\{x \in A_0:{{n-1}\over{q}} < \tilde{\Pi}(x) \leq {{n}\over{q}} \right\}.$$
We see 
$$B \sim \Lambda_{\infty}^{-1} = \bigcup_{n=1}^{\infty}H_{nq}.$$
we now estimate $\Hm{1}(\F1(B \cap H_{nq}))$. Firstly  $H_{nq} \subset \Lambda_{n/q}$ so that 
$$\F1(B \cap H_{nq}) = \F1|_{\Lambda_{n/q}}(B \cap H_{nq})$$ 
is a Lipschitz graph with Lip$\F1|_{\Lambda_{n/q}} \leq n/q$ 
so that 
$$\Hm{1}(\F1(B \cap H_{nq})) \leq {{n}\over{q}}\Hm{1}(H_{nq}).$$
It is now necessary to establish a lower estimate. To do this we define 
$$H_{nqj}:=\{x \in H_{nq}:\tilde{\Pi}_j(x) > (n-1)q^{-1} \geq \tilde{\Pi}_{j-1}(x)\}$$
and note that $H_{nqi} \cap H_{nqj} = \emptyset $ whenever $i \not= j$.
We also define 
$$J_{nq}:=\{i \in \{1,...,2^j\}:i=i(n,x) \hbox{ for some } x \in H_{nqj}\}$$
We note that $\F1|_{\Lambda_{n/q}} \circ F_j^{-1}$ is a Lipschitz expansion map on $F_j(\Lambda_{n/q})$.
It follows that
\begin{eqnarray}
\Hm{1}(\F1 (H_{nqj})) & = & \Hm{1}(\F1 \circ F_j^{-1} \circ F_j(H_{nqj})) \nonumber \\
& = & \Hm{1} \left(\bigcup_{i \in J_{nq}}\F1|_{\Lambda_{n/q}} \circ F_j^{-1}(F_{j}(H_{nqj}) \cap A_{j,i}) \right) 
\nonumber \\
& = & \sum_{i \in J_{nq}}\Hm{1}(\F1|_{\Lambda_{n/q}} \circ F_j^{-1}(F_{j}(H_{nqj}) \cap A_{j,i})) \nonumber \\
& \geq & \sum_{i \in J_{nq}} \Hm{1}(F_{j}(H_{nqj}) \cap A_{j,i}) \nonumber \\
& = & \sum_{i \in J_{nq}} \tilde{\Pi}_{j,i}\Hm{1}(H_{nqj} \cap [(i-2)2^{-j},i2^{-j}]) \nonumber \\
& > & \sum_{i \in J_{nq}}{{n-1}\over{q}}\Hm{1}(H_{nqj} \cap[(i-1)2^{-j},i2^{-j}]) \nonumber \\
& = & {{n-1}\over{q}}\Hm{1}(H_{nqj}). \nonumber 
\end{eqnarray}
Since  $H_{nq}$ is the disjoint union of $\{H_{nqj}\}_{j=1}^{\infty}$ it follows that 
\begin{eqnarray}
\Hm{1}(\F1 (H_{nq})) & = & \sum_{j=1}^{\infty}\Hm{1}(\F1(H_{nqj})) \nonumber \\
& > & {{n-1}\over{q}}\sum_{j=1}^{\infty}\Hm{1}(H_{nqj}) \nonumber \\
& = & {{n-1}\over{q}}\Hm{1}(H_{nq}). \nonumber
\end{eqnarray}
It then follows that 
$${{n-1}\over{q}}\Hm{1}(H_{nq}) \leq \Hm{1}(\F1 (H_{nq}) \leq {{n}\over{q}}\Hm{1}(H_{nq}).$$
Correspondingly we have direct from the definition of $H_{nq}$ that 
$${{n-1}\over{q}}\Hm{1}(H_{nq}) < \int_{H_{nq}}\tilde{\Pi}d\Hm{1} \leq {{n}\over{q}}\Hm{1}(H_{nq})$$
so that 
$$\left|\int_{H_{nq}}\tilde{\Pi}d\Hm{1}-\Hm{1}(\F1 (H_{nq})\right| < {{1}\over{q}}\Hm{1}(H_{nq})$$
and therefore 
\begin{eqnarray}
\left|\int_{B \sim \Lambda_{\infty}^{-1}}\tilde{\Pi}d\Hm{1} - \Hm{1}(\F1 (B \sim \Lambda_{\infty}^{-1}))\right| & = & 
\left|\sum_{n=1}^{\infty}\int_{H_{nq}}\tilde{\Pi}d\Hm{1} - \Hm{1}(\F1 (H_{nq}))\right| \nonumber \\
& \leq & \sum_{n=1}^{\infty}\left|\int_{H_{nq}}\tilde{\Pi}d\Hm{1} - \Hm{1}(\F1(H_{nq}))\right| \nonumber \\
& < & \sum_{n=1}^{\infty}{{1}\over{q}}\Hm{1}(H_{nq}) \nonumber \\
& = & {{1}\over{q}}\Hm{1}(B \sim \Lambda_{\infty}^{-1}) \nonumber \\
& \leq & {{1}\over{q}}. \nonumber
\end{eqnarray}
Since this is true for all $q \in \mathbb{N}$ it follows that 
$$\left|\int_{B \sim \Lambda_{\infty}^{-1}}\tilde{\Pi}d\Hm{1} - \Hm{1}(\F1 (B \sim \Lambda_{\infty}^{-1}))\right|=0$$
and thus that 
$$\int_{B \sim \Lambda_{\infty}^{-1}}\tilde{\Pi}d\Hm{1} = \Hm{1}(\F1 (B \sim \Lambda_{\infty}^{-1})).$$
This gives us 
\begin{eqnarray}
\Hm{1}(\F1(B)) & = & \Hm{1}(\F1(B) \cap \Lambda_{\infty})+\Hm{1}(\F1(B) \sim \Lambda_{\infty}) \nonumber \\
& = & \int_{B \sim \Lambda_{\infty}^{-1}}\tilde{\Pi}d\Hm{1} + \Hm{1}(\F1(B) \cap \Lambda_{\infty}). \nonumber 
\end{eqnarray}
\end{proof}
\end{thm}\noindent
As we mentioned at the beginning of this chapter, we present the simplified result for $\A$ type sets. In this case, 
however, the result does not simplify. This is because, should $\tilde{\Pi}^A_{\cdot}\equiv\infty$ for some 
$\A$ type set $A$ then it could be this very $A$ that allows for creation of measure. Then for any set $B \subset A_{0,0}$ 
with $\Hm{1}(B)>0$ we get 
$$\int_B\tilde{\Pi}d\Hm{1} = \Hm{1}(\F1(B)) = \infty.$$
However, for a measurable set $B \subset A_{0,0}$ with $\Hm{1}(B)=0$ from which measure is created we would have 
$$\int_B\tilde{\Pi}d\Hm{1} = 0$$
but 
$$\Hm{1}(\F1(B))>0$$
preventing the simplified version of Theorem $\ref{thm12}$ 
$$\int_B\tilde{\Pi}d\Hm{1}=\Hm{1}(\F1(B))$$
holding as desired.
\newline \newline
This, therefore, concludes our discussion of measure formulae and we now conclude with the results on dimension.
\section{A Full Spectrum of Dimension}
We complete this work with a discussion of the Dimension of $\A$ and Koch type sets. 
As we discussed earlier in this Chapter, in 
order to gather results about dimension we essentially want to place sets either inside of or around sets that we 
know the dimension of. Unfortunately, generally with different $\A$ type sets they do not generally stay neatly 
inside of one another. We therefore need to use our centralisation results to rearrange each stage of construction 
to ensure that strict containment is retained by the necessary sets.
\newline \newline
As with the rectifiability results, the $\A$ type sets allow for a more cleanly stated result than the Koch type sets. 
Unlike some of the previous result, we shall not prove the asthetically more pleasing results of the $\A$ type sets 
as a corollary of the more general Koch type sets but shall rather prove the result directly. This is mainly because the 
proof attached to the $\A$ type sets is much cleaner allowing the essential ingredients to be more clearly seen. 
The proof associated with the Koch type sets is then presented afterwards where the difficulties of allowing full 
variation of base angles require a much more technical proof.
\newline \newline
As we will see from the results, a complete closed interval in $\R$ represents the possible dimensions of sets in $\cl{K}$. 
This shows the rich variation of the sets, which could otherwise perhaps have been of a dimension from a finite set of 
values.
\newline \newline
Following the proof of the dimension of the $\A$ type sets, we present a Corollary showing how the dimension of $\A$ 
(which we directly proved to be $1$ in Theorem $\ref{thm1}$) follows easily from the more general result.
\begin{thm}\label{thm13}
For $r \geq 0$ and $A \in A^r$
$$dimA = -{{ln2}\over{ln({{1}\over{2}}(1+(tan(r))^2)^{1/2})}}.$$
\pf
The proof is dependent on the dimension of $\G$. We thererfore first note that for any scaling $\lambda \in \R$ 
$$dim\lambda\G = dim\G.$$
We also note that 
$$\Gamma_{1/2(tanr)} \in A^r$$
and finally, recalling $dim\G = -ln2/(lnl)$ where $l$ is the shrinking factor per approximation stage, we 
calculate that $l_r$, the approriate $l$ for $\G \in A^r$ is 
$$l_r = -{{ln2}\over{ln({{1}\over{2}}(1+(tan(r))^2)^{1/2})}}.$$
Now, since for $A \in A^r$ $\theta_n^A \searrow r$ we have $\theta_n^A \geq r$ for all $n \in \mathbb{N}$. Thus, since 
$\theta_n^{\Gamma_{1/2(tanr)}} \equiv r$ for all $n \in \mathbb{N}$ 
and thus also $\Hm{1}(A_{0,1}^A)T_{0,1}^{\Gamma_{1/2(tanr)}} \subset T_{0,1}^A$ Proposition $\ref{prop18}$ 
then gives us that 
$$\Hm{1}(A_{0,1}^A)\Gamma_{1/2(tanr)} \Cent A.$$
Lemma $\ref{lem16}$ then gives 
\begin{equation}
dimA \geq dim\Hm{1}(A_{0,1}^A)\Gamma_{1/2(tanr)} = dim\Gamma_{1/2(tanr)}. \label{e:gdim1}
\end{equation}
Then, for any $r_1 > r$ there is an $n_0 \in \mathbb{N}$ such that for all $n > n_0$ $\theta_n^A \leq r_1$. It follows that 
by choosing arbitrarily and $j \in \{1,...,2^{n_0}\}$
$$\Hm{1}(A_{n_0,j}^A)T_{0,1}^{\Gamma_{1/2(tanr_1)}} \supset T_{n_0,j}^A.$$
Now taking $T_j \in A^r$ to be the set generated by starting with $T_{n_0,j}^A$ and 
$\theta_n^{T_j} \equiv \theta_{n+n_0}^A$, we have by Proposition $\ref{prop18}$ that 
$$T_j \Cent \Hm{1}(A_{n_0,j}^A)\Gamma_{1/2(tanr_1)}.$$
It then follows from Lemma $\ref{lem16}$ that 
$$dimT_j \leq dim \Hm{1}(A_{n_0,j}^A)\Gamma_{1/2(tanr_1)} = dim\Gamma_{1/2(tanr_1)}.$$
Taking a finite union of such ses will not alter the dimension, thus 
\begin{eqnarray}
dim A & = & dim\bigcup_{j=1}^{2^{n_0}}T_j \nonumber \\
& = & dim T_j \nonumber \\
& \leq & dim\Gamma_{1/2(tanr_1)} \nonumber \\
& = & -{{ln2}\over{ln({{1}\over{2}}(1+(tan(r))^2)^{1/2})}}. \nonumber
\end{eqnarray}
Since this is true for all $r_1>r$ it follows that 
$$dimA \leq -{{ln2}\over{ln({{1}\over{2}}(1+(tan(r))^2)^{1/2})}} = dim\Gamma_{1/2(tanr)}.$$
Combining this with ($\ref{e:gdim1}$) gives the result 
\end{proof}
\end{thm}\noindent
\begin{cor}\label{cor10}\thst
dim$\A =1$.
\pf
Since from Proposition $\ref{prop12}$ we know $\A \in A^0$ for any given $\e$, we can directly apply Theorem $\ref{thm13}$ 
to calculate
\begin{eqnarray}
dim\A & = & -{{ln2}\over{ln({{1}\over{2}}(1+(tan(0))^2)^{1/2})}} \nonumber \\
& = & -{{ln2}\over{ln(1/2)}} \nonumber \\
& = & 1. \nonumber 
\end{eqnarray}
\end{proof}
\end{cor}\noindent
Our final result is then the characterisation of dimension for the more general Koch type sets. As we see, the 
basic principle is the same as that used for $\A$ type sets, the difference being the need to adjust for individually 
varying rates of change of base angle in the more general set up. We slowly eliminate those more rapidly 
decreasing, leaving those with a base measure enough to make a difference that reduce base angle slowly 
and would then, in the sense of Theorem $\ref{thm13}$ have higher dimension. It is these sets that dictate the 
dimension of the general whole set.
\begin{thm}\label{thm14}\thst
Let $A \in \cl{K}$ and 
$$\gamma_1^A = \sup\{a : \Hm{1}(\{x \in A_0:\lim_{n \rightarrow \infty}\theta_{n,i(n,x)}^A \geq a\})>0\}$$
and
$$\gamma_2^A = \sup_{x \in A_0}\tilde{\theta}^A_x.$$
Then
$$dim\Gamma_{f(\gamma^A_1)} = f_1(\gamma^A_1) \leq dimA \leq f_1(\gamma^A_2) =dim\Gamma_{f(\gamma^A_2)}$$
where
$$f(\gamma):= (1/2)(tan\gamma )$$
and therefore
$$f_1(\gamma):= -{{ln2}\over{ln((1/2)(1+(tan\gamma )^2)^{1/2})}}.$$
Should the hypothesis that for $B \subset A_0$ $\Hm{1}(B) = 0 \Rightarrow \Hm{1}(\F1(B))=0$ hold, or should for a given 
$A \in \cl{K}$ we have $\Hm{1}(\Upsilon_{\gamma_1^A +})=0$ then
$$dim A \equiv f_1(\gamma^A_1).$$
\pf
We start by proving that $dim A \leq f_1(\gamma_1^A)$.
\newline \newline
Let $\xi < \gamma_1^A$. Then $\Hm{1}(\Upsilon_{\xi +}^{-1})>0$. There is therefore an $n_0 \in \mathbb{N}$ such that 
$\Hm{1}(\Upsilon_{\xi +}^{-1}) > 2^{-n_0}$. It follows that 
$$\Upsilon_{\xi +}^{-1} \cap [(i-1)2^{-n_0},i2^{-n_0}] \not= \emptyset$$
for at least $2^{n-n_0}$ $i \in \{1,...,2^n\}$. It follows that 
$$T_{n,i}^A \cap F_n(\Upsilon_{\xi +}^{-1}) \not= \emptyset$$
for at least $2^{n-n_0}$ $i \in \{1,...,2^n\}$. \newline \newline
In particular, this is true for all $n \geq n_0$. We set 
$$A_{2m}:= \cup \{T_{n_0+m,i}^A:T_{n_0+m,i}^A \cap F_{n_0 +m}(\Upsilon_{\xi +}^{-1}) \not= \emptyset\}$$ 
and $n_2(m)$ to be the number of $T_{n_0 +m,i}^A$ that are included in $A_{2m}$.
\newline \newline
Note that $n_2(m) \geq 2^m.$ Further we order these $T_{n_0+m,i}^A \subset A_{2m}$ as $\{A_{2mj}\}_{j=1}^{n_2(m)}$.
We consider the set $\Gamma_{f(\xi)}$ constructed on a base $A_{0,0}$ of length 
$$\Hm{1}(A_0) = 2^{-n_0}\prod_{i=0}^{n_0-1}(cos \xi)^{-1}.$$ We denote this set by $\tilde{\Gamma}$. \newline \newline
We now want to show that 
$$\tilde{\Gamma} \Cent A_2 := \bigcap_{m=1}^{\infty}A_{2m}.$$
Clearly, for any $T_{n_0+m+1,i}^A \subset A_{2(m+1)}$, 
$T_{n_0+m+1,i}^A \cap F_{n_0+m+1}(\Upsilon_{\xi +}^{-1}) \not= \emptyset$, also 
$T_{n_0+m+1}^A \subset T_{n_0+m,int(i/2)+1}^A$ so that 
$$T_{n_0+m,int(i/2)+1}^A \cap F_{n_0+m+1}(\Upsilon_{\xi +}^{-1}) \not= \emptyset$$
and thus
$$T_{n_0+m,int(i/2)+1}^A \cap F_{n_0+m}(\Upsilon_{\xi +}^{-1}) \not= \emptyset,$$
so that $T_{n_0+m,\int(i/2)+1}^A \subset A_{2m}$ and hence we have $A_{2(m+1)} \subset A_{2m}$ for any $m \in \mathbb{N}$.
\newline \newline
We see that in putting $\tilde{\Gamma}$ into the required form for Definition $\ref{def31}$ 
$$\tilde{\Gamma} = A_1, T_n^{\tilde{\Gamma}}:=\bigcup_{i=1}^{2^n}T_{n,i}^{\tilde{\Gamma}} = A_{1n}, 
T_{n,i}^{\tilde{\Gamma}} = A_{1ni}, \hbox{ and } n_1(m) = 2^m.$$
So that $A_1$ and $A_2$ individually satisfy the requirements of $A_1$ and $A_2$. Also, $n_1(m) = n_2(m)$. We therefore 
only need to show the existence of the transformations $\cl{T}_{n,i}^{A_1,A_2}$.
\newline \newline
We note that each $A_{1mi} = T_{m,i}^{\tilde{\Gamma}}$ is a triangular cap of base length 
$$2^{-m}\left(\prod_{i=0}^{m-1}(cos  \xi)^{-1}\right) \times (\hbox{base length }T_{0,1}^{\tilde{\Gamma}})$$
which equals
$$2^{-m-n_0}\prod_{i=0}^{m+n_0-1}(cos \xi)^{-1}$$
and of base angle $\xi$.
\newline \newline
We also note that for each $i \in \{1,...,n_1(m)\}$, $i \in \{1,...,n_2(m)\}$ so that $A_{2mi}$ exists and is a 
triangular cap $T_{n_0+m,k}^A$ for some $k\in \{1,...,2^{n_0+m}\}$ 
with base angle $\theta_{n_0+m,k}^A \geq \xi$ and base length 
$$2^{-n_0-m}\tilde{\Pi}_{n_0+m,k}^A.$$
Since a sequence $\{\theta_{n,i(n)}\}$ of angles in the construction of $A$ is decreasing and $\theta_{n_0+m,k}^A \geq \xi$
it follows that 
$$2^{-n_0-m}\tilde{\Pi}_{n_0+m,k}^A \geq 2^{-m-n_0}\prod_{i=0}^{m+n_0-1}(cos \xi)^{-1} = \Hm{1}(A_{1mi}).$$
It follows, since $A_{1mi}$ and $A_{2mi}$ are isoceles triangles where $A_{2mi}$ has a longer base and larger base angles 
that $A_{2mi}$ is strictly larger than $A_{1mi}$ in the sense that $A_{1mi}$ could be placed inside of $A_{2mi}$ and 
thus there must exist an orthogonal transformation $\cl{T}_{m,i}^{\tilde{\Gamma},A_2}$ such that 
$$\cl{T}_{m,i}^{\tilde{\Gamma},A_2}(A_{1mi}) \subset A_{2mi}.$$
Since this is true for any $m \in \mathbb{N}$ and $i \in \{1,...,n_1(m)\}$ it follows that $\tilde{\Gamma} \Cent A_2$.
\newline \newline
Thus, using Lemma $\ref{lem16}$ and the fact that $A_2 \subset A$ we have 
$$dim \Gamma_{f(\xi)} = dim\tilde{\Gamma} \leq dimA_2 \leq dimA.$$
Since this is true for each $\xi < \gamma_1^A$ it follows that 
$$dimA \geq dim\Gamma_{f(\gamma_1^A)} = f_1(\gamma_1^A).$$
For the $\leq $ inequalities, we let $B \subset A_0$ be $\Hm{a}$-measurable for each $a \in \R$ and show that for 
$$\gamma = \sup_{x \in B}\tilde{\theta}^A_x$$
$$dim\F1(B) \leq dim\Gamma_{f(\gamma)}=f_1(\gamma).$$
Let $\xi > \gamma$ and for each $n \in \mathbb{N}$ define 
$$\chi_n:=\cup\{T_{n,i}^A:\theta_{n,i}^A \geq \xi\}.$$
Then $\Psi_n:= T_n - \chi_n$ is the finite union of triangular caps $T_{n,j}^A$ with $\theta_{n,j}^A \leq \xi$.
\newline \newline
We see that for each such triangular cap $T_{n,j}^A \subset \Phi_n$, 
$$\Hm{1}(A_{n,j}^A) \leq \Hm{1}(A_{0,1}^A) = \Hm{1}(A_{0,1}^{\Gamma_{f(\xi)}})$$
and that for each later triangular cap $T_{n,m,k}^A \subset T_{n,j}^A$
$$\theta_{n+m,k}^A \leq \theta_{n,j}^A \leq \xi = \theta_{n+m,\cdot}^{\Gamma_{f(\xi)}}.$$
It therefore follows from Proposition $\ref{prop18}$ that for each $T_{n,j}^A \subset \Psi_n$
$$T_{n,j}^A \Cent \Gamma_{f(\xi)}$$
and hence, since $A \cap T_{n,j}^A$ equals the final set resulting from the Koch set construction starting fom $T_{n,j}^A$, 
Lemma $\ref{lem16}$ gives 
$$dim(A \cap T_{n,j}^A) \leq dim\Gamma_{f(\xi)}$$
and therefore, since this is true for any such triangular cap, that 
$$dim(A \cap \Psi_n) = dim(A\cap T_{n,j}^A) \leq dim\Gamma_{f(\xi)}.$$
Now, suppose that there exists a $y \in \F1(B)$ with 
$$y \not\in \bigcup_{n=1}^{\infty}\Psi_n.$$
Then for each $n \in \mathbb{N}$ $\theta_{n,i(n,y)}^A \geq \xi$ and therefore 
$$\tilde{\theta}^A_y=\lim_{n \rightarrow \infty}\theta_{n,i(n,y)}^A \geq \xi > \gamma.$$
Since this is impossible it follows that $\F1(B) \subset \cup_{n=1}^{\infty}(\Psi_n \cap A)$ and therefore that 
$$dim\F1(B) \leq dim\Gamma_{f(\xi)}.$$
Since this is true for each $\xi > \gamma$ we have 
$$dim\F1(B) \leq dim\Gamma_{f(B)} = f_1(B).$$
To finish the proof we note that $f_1(\gamma)\geq 1$ for each $\gamma \geq 0$, and consider firstly 
that for each $x \in A_0$, $\tilde{\theta}^A_x \leq \gamma^A_2$ so that immediately from the above we have 
$$dimaA \leq dim\Gamma_{f(\gamma^A_2)} = f_1(\gamma_2^A).$$
For the second conclusion we consider 
$$B= \Upsilon_{\gamma^A_1}^{-1}.$$
It follows then that 
$$dim\Upsilon_{\gamma^A_1 +} \leq dim\Gamma_{f(\gamma^A_1)} =f_1(\gamma^A_1).$$
Should the hypothesis hold that for all $D \subset A_0$, $\Hm{1}(D)= 0\Rightarrow \Hm{1}(\F1(D))=0$, or should we directly 
have $\Hm{1}(\Upsilon_{\gamma^A_1 +})=0$, then we have $\Hm{1}(\Upsilon_{\gamma^A_1 +})=0$ and therefore 
$$dim\Upsilon_{\gamma^A_1 +} \leq dim\Gamma_{f(\gamma^A_1)} = f_1(\gamma_1^A).$$
We therefore have
\begin{eqnarray}
dimA & \leq & max\{dim \Upsilon_{\gamma^A_1}, \Upsilon_{\gamma^A_1 +}\} \nonumber \\
& \leq & dim\Gamma_{f(\gamma^A_1)} \nonumber \\
& = & f_1(\gamma_1^A), \nonumber  
\end{eqnarray}
which completes the proof 
\end{proof}
\end{thm}

\end{document}